\documentclass[11pt,reqno]{amsart}
\usepackage{amsmath} 
\usepackage{amssymb}
\usepackage{amsfonts}

\begin{document}

\title[Analytic continuation of the resolvent]
{Analytic continuation of the resolvent of the Laplacian and the dynamical zeta function}
\author[V. Petkov]{Vesselin Petkov}
\address{Universit\'e Bordeaux I, Institut de Math\'ematiques de Bordeaux, 351,
Cours de la Lib\'eration, 
33405  Talence, France}
\email{petkov@math.u-bordeaux1.fr}
\author[L. Stoyanov]{Luchezar Stoyanov}
\address{University of Western Australia, School of Mathematics 
and Statistics, Perth, WA 6009, Australia}
\email{stoyanov@maths.uwa.edu.au}

\maketitle

\def\nexto{\kern -0.54em}
\newcommand{\C}{\protect\mathbb{C}}
\newcommand{\R}{\protect\mathbb{R}}
\newcommand{\Q}{\protect\mathbb{Q}}
\newcommand{\Z}{\protect\mathbb{Z}}
\newcommand{\N}{\protect\mathbb{N}}
\newtheorem{thm}{Theorem}
\newtheorem{prop}{Proposition}
\newtheorem{lem}{Lemma}
\newtheorem{col} {Collorary}
\newtheorem{deff}{Definition}
\def\DF{{\rm {I\ \nexto F}}}
\def\DR{{\rm {I\ \nexto R}}}
\def\DN{{\rm {I\ \nexto N}}}
\def\DZ{{\rm {Z \kern -0.45em Z}}}
\def\DC{{\rm\hbox{C \kern-0.83em\raise0.08ex\hbox{\vrule height5.8pt width0.5pt} \kern0.13em}}}
\def\DQ{{\rm\hbox{Q \kern-0.92em\raise0.1ex\hbox{\vrule height5.8pt width0.5pt}\kern0.17em}}}
\def\T{{\mathbb T}}
\def\S{{\mathbb S}}
\def\sn{{\S}^{n-1}}
\def\sN{{\S}^{N-1}}

\def\ep{\epsilon}
\def\e{\emptyset}
\def\di{\displaystyle}
\def\sk{\smallskip}
\def\bs{\bigskip}
\def\ms{\medskip}
\def\dk{\partial K}
\def\kamin{\kappa_{\min}}
\def\kamax{\kappa_{\max}}
\def\saa{\Sigma_A^+}
\def\sa{\Sigma_A}
\def\scc{\Sigma^+_C}
\def\sbb{\Sigma_B^+}
\def\san{\Sigma^-_A}

\def\Prf{\mbox{\footnotesize\rm Pr}}
\def\Pr{\mbox{\rm Pr}}

\def\be{\begin{equation}}
\def\ee{\end{equation}}
\def\beqn{\begin{eqnarray}}
\def\eeqn{\end{eqnarray}}
\def\beqn*{\begin{eqnarray*}}
\def\eeqn*{\end{eqnarray*}}
\def\endofproof{{\rule{6pt}{6pt}}}

\def\i{{\bf i}}
\def\dist{\mbox{\rm dist}}
\def\diam{\mbox{\rm diam}}
\def\pr{\mbox{\rm pr}}
\def\supp{\mbox{\rm supp}}
\def\Arg{\mbox{\rm Arg}}
\def\In{\mbox{\rm Int}}
\def\Im{\mbox{\rm Im}}
\def\span{\mbox{\rm span}}
\def\con{c}
\def\const{\mbox{\rm c}}
\def\Con{C}
\def\spec{\mbox{\rm spec}\,}
\def\Re{\mbox{\rm Re}}
\def\ecc{\mbox{\rm ecc}}
\def\var{\mbox{\rm var}}
\def\conf{{\footnotesize c}}
\def\Conf{{\footnotesize C}}

\def\do{{\mathcal D}_0}
\def\d1{{\mathcal D}_1}
\def\ff{{\mathcal F}}
\def\kk{{\mathcal K}}
\def\cc{{\mathcal C}}
\def\mm{{\mathcal M}}
\def\gg{{\mathcal G}}
\def\T{{\mathcal T}}
\def\rr{{\mathcal R}}
\def\ll{{\mathcal L}}
\def\uu{{\mathcal U}}
\def\ss{{\mathcal S}}
\def\oo{{\mathcal O}}
\def\rr{{\mathcal R}}
\def\pp{{\mathcal P}}
\def\gg{{\mathcal G}}
\def\vv{{\mathcal V}}

\def\iii{{\imath }}
\def\hy{\hat{y}}
\def\hla{\hat{\lambda}}
\def\hx{\hat{x}}
\def\hxi{\hat{\xi}}
\def\hu{\hat{u}}
\def\hf{\hat{f}}
\def\hg{\hat{g}}
\def\tf{\tilde{f}}
\def\ho{\hat{\omega}}
\def\hD{\hat{\Delta}}
\def\hxi{\hat{\xi}}
\def\heta{\hat{\eta}}
\def\hQ{\hat{Q}}
\def\nv{\nabla \varphi}
\def\Od{\mbox{\rm Od}}
\def\vixi{\varphi^{(\infty)}_{\xi,0}}
\def\Ref{\mbox{\footnotesize\rm Re}}
\def\h{h_{\mbox{\rm\footnotesize top}}}
\def\tc{\tilde{c}}
\def\tx{\tilde{x}}
\def\ty{\tilde{y}}
\def\tz{\tilde{z}}
\def\tu{\tilde{u}}
\def\tv{\tilde{v}}
\def\ta{\tilde{a}}
\def\td{\tilde{d}}
\def\tf{\tilde{f}}
\def\tg{\tilde{g}}
\def\tl{\tilde{\ell}}
\def\tr{\tilde{r}}
\def\tt{\tilde{t}}
\def\tw{\tilde{w}}
\def\tth{\tilde{\theta}}
\def\tka{\tilde{\kappa}}
\def\tla{\tilde{\lambda}}
\def\tmu{\tilde{\mu}}
\def\tga{\tilde{\gamma}}  
\def\trho{\tilde{\rho}}
\def\tvar{\tilde{\varphi}}
\def\tF{\tilde{F}}
\def\ii{{\bf \i}}
\def\tU{\tilde{U}}
\def\tV{\tilde{V}}
\def\tS{\tilde{S}}
\def\tC{\widetilde{C}}
\def\tP{\widetilde{P}}
\def\tQ{\widetilde{Q}}
\def\tR{\tilde{R}}
\def\v{{\sf v}}
\def\wuloc{W^u_{\mbox{\footnotesize\rm loc}}}
\def\jj{{\bf j}}
\def\psim{\psi^{(m)}}
\def\psimm{\psi^{(m)}_m}
\def\psime{\psi^{(m,\eta)}}
\def\psimme{\psi^{(m,\eta)}_m}
\def\phim{\varphi^{(m)}}
\def\phimm{\varphi^{(m)}_m}
\def\ovo{\overset{\circ}\Omega}
\def\Wn{W^{(n+2)}}
\def\tWn{\widetilde{W}^{(n+2)}}
\def\cu{{\mathcal U}}
\def\mt{\Lambda}
\def\hU{\widehat{U}}
\def\clip{C^{\mbox{\footnotesize \rm Lip}}}
\def\Lip{\mbox{\rm Lip}}
\def\lip{\mbox{{\footnotesize\rm Lip}}}
\def\hz{\hat{z}}
\def\tz{\tilde{z}}
\def\mtb{ \mt_{\dk}}
\def\cdo{{\mathcal D}_0}
\def\d1{{\mathcal D}_1}
\def\cd2{{\mathcal D}_2}
\def\wuloc{W^u_{\mbox{\footnotesize\rm loc}}}
\def\wsloc{W^s_{\mbox{\footnotesize\rm loc}}}
\def\la{\langle}
\def\ra{\rangle}
\def\Lip{\mbox{\rm Lip}}
\def\ka{\kappa}
\def\bv{\big\vert}
\def\wj{w_{0,j}}
\def\gj{\Gamma_j}
\def\d2{{\mathcal D}_2}
\def\d1{{\mathcal D}_1}
\def\cd1{{\mathcal D}_1}

\begin{abstract}
Let $s_0 < 0$ be the abscissa of absolute convergence of the dynamical zeta function $Z(s)$ for several 
disjoint strictly convex compact obstacles $K_i \subset \R^N, i = 1,\ldots, \kappa_0,\: \ka_0 \geq 3,$ and let 
$R_{\chi}(z) = \chi (-\Delta_D - z^2)^{-1}\chi,\: \chi \in C_0^{\infty}(\R^N),$ be the cut-off resolvent of the 
Dirichlet Laplacian $-\Delta_D$ in $\Omega = \overline{\R^N \setminus \cup_{i = 1}^{k_0} K_i}$. We prove 
that there exists $\sigma_1 < s_0$ such that the cut-off  resolvent $R_{\chi}(z)$ has an analytic continuation for $\Im (z) < - \sigma_1,\: |\Re (z)| \geq J_1 > 0.$
\end{abstract}

\section{Introduction}
\renewcommand{\theequation}{\arabic{section}.\arabic{equation}}
\setcounter{equation}{0}
\def\d1{{\mathcal D}_1}

Let $K$ be a subset of ${\R}^{N}$ ($N\geq 2$) of the form $K = K_1 \cup K_2 \cup \ldots \cup K_{\ka_0},$where $K_i$ are compact
strictly convex disjoint domains in $\R^{N}$ with $C^\infty$ {\it boundaries} $\Gamma_i = \dk_i$ and $\ka_0 \geq 3$. Set 
$\Omega = \overline{{\R}^N \setminus K}$ and $\Gamma = \partial K$. We assume that $K$ satisfies the following (no-eclipse) condition: 
$${\rm (H)} \quad \quad\qquad  
\begin{cases}
\mbox{\rm for every pair $K_i$, $K_j$ of different connected components  of $K$ the convex hull of}\cr
\mbox{\rm $K_i\cup K_j$ has no common points with any other connected component of $K$. }\cr
\end{cases}$$
With this condition, the {\it billiard flow} $\phi_t$ defined on the {\it cosphere bundle} $S^*(\Omega)$ 
in the standard way is called an open billiard flow.It has singularities, however its restriction to the {\it non-wandering set} $\Lambda$ has only 
simple discontinuities at reflection points.  Moreover, $\Lambda$  is compact, $\phi_t$ is hyperbolic and transitive on $\Lambda$, and  it follows from 
\cite{kn:St1} that $\phi_t$ is  non-lattice and therefore by  a result of Bowen  \cite{kn:Bo1}, it is topologically weak-mixing on $\Lambda$.

Given a periodic reflecting ray $\gamma \subset \Omega$ with $m_{\gamma}$ 
reflections, denote by $d_{\gamma}$ the period (return time) of
$\gamma$, by $T_{\gamma}$ the  primitive period (length) of $\gamma$
and by $P_{\gamma}$ the linear  Poincar\'{e} map associated to $\gamma$.
Denote by $\Pi$ the set of all periodic rays in  $\Omega$ and let 
$\lambda_{i, \gamma}, i = 1,\ldots,N -1,$ be the eigenvalues of 
$P_{\gamma}$ with $|\lambda_{i, \gamma}| > 1$ (see \cite{kn:PS1}). 

Let ${\mathcal P}$ be the set of primitive periodic rays. Set
$$\delta_{\gamma} = -\frac{1}{2} \log (\lambda_{1,\gamma}\ldots 
\lambda_{N-1,\gamma}),\: \gamma \in {\mathcal P},$$
$$r_{\gamma} = \begin{cases} 0 \:\; {\rm if}\:\: m_{\gamma}\:\text {is   even}, \cr
1 \:\: {\rm if} \:\; m_{\gamma}\:\text {is odd} \;,\end{cases}$$
and consider the {\it dynamical zeta function}
$$Z(s) = \sum_{m = 1}^{\infty} \frac{1}{m}\sum_{\gamma \in {\mathcal P}} 
(-1)^{mr_{\gamma}}  e^{m(-s T_{\gamma} + \delta_{\gamma})}.$$
It is easy to show that there exists $s_0 \in \R$ such that  for $\Re(s) > s_0$ the series $Z(s)$ is
absolutely convergent and $s_0$ is minimal with this property. The number $s_0$ is called {\it abscissa of absolute convergence}. 
On the other hand, using symbolic dynamics and the results of \cite{kn:PP}, it follows that $Z(s)$ is 
meromorphic for $\Re( s) > s_0 - a, \: a > 0$ (see \cite{kn:I4}) and $Z(s)$ is analytic for for $\Re (s) \geq s_0.$  According to some recent 
results (\cite{kn:St2} for $N = 2$, \cite{kn:St4} for $N \geq 3$ under some additional conditions) there exists $0 <\epsilon < a$ so that the dynamical zeta function 
$Z(s)$ admits an analytic continuation for $\Re( s) \geq s_0 - \epsilon.$

The {\it cut-off resolvent} defined by
$$R_{\chi}(z)=\chi (-\Delta_K - z^2)^{-1} \chi:\: L^2(\Omega) \longrightarrow L^2(\Omega)$$
for $\Im (z) < 0$, where $\chi \in C_0^{\infty}(\R^N),\: \chi = 1$ on $K$, and 
$\Delta_K$ is the Dirichlet Laplacian in $\Omega$,  has a meromorphic continuation 
in $\C$ for $N$ odd with poles $z_j$ such that $\Im (z_j) > 0$ and in $\C \setminus \{i\R^{+}\}$ for $N$ even. The analytic properties and the estimates of 
$R_{\chi}(z)$ play a crucial role in many problems related to the local energy  decay, distribution of the resonances etc. In the physical  literature and in many 
works concerning numerical calculation of resonances (see \cite{kn:CE}, \cite{kn:W}, \cite{kn:L}, \cite{kn:LZ}, \cite{kn:LSZ}) the following conjecture is often made.

\bs

\noindent
{\bf Conjecture}:\: {\it The poles $\mu_j$ (with $\Re (\mu_j) < 0$) of $Z(s)$ and the poles $z_j$ of $R_{\chi}(z)$ are related by $\ii z_j = \mu_j$.}

\bs

At least one would expect that the poles $z_j$ of $R_{\chi}(z)$ lie in  sufficiently small neighborhoods of $-\ii \mu_j$. Presumably for
this reason the numbers $-\ii \mu_j$ are called {\it pseudo-poles} of $R_{\chi}(z)$. 

The case of several disjoint disks has been treated in many works (see \cite{kn:W} for a comprehensive list of references), and a certain
method for numerical computation of the resonances has been used. Although it is not rigorously known whether the numerically found resonances
approximate the (true) resonances in the exterior of the discs, and whether the dynamical zeta function has an analytic continuation  
to the left of the line of absolute convergence, this way of computation is widely accepted in the physical literature.

In the case of two strictly convex disjoint domains it was proved  (\cite{kn:I1}, 
\cite{kn:G}) that the poles of $R_{\chi}(\lambda)$  are contained in small neighborhoods of the pseudo-poles
$$ m \frac{\pi}{d} + \ii \alpha_k ,\: m \in \Z, \:k \in N\; .$$
Here $d > 0$ is the distance between the obstacles and $\alpha_k > 0$ 
are determined by the eigenvalues $\lambda_j$ of the Poincar\'e map 
related to the unique primitive periodic ray.

It is  known  that the above conjecture is true for convex co-compact 
hyperbolic manifolds $X = {\bf \Gamma} \setminus {\mathbb H}^{n+1}$, where 
${\bf \Gamma}$ is a discrete group of isometries with only hyperbolic
elements admitting a finite fundamental domain (then $X$ is a manifold
of constant negative curvature).  More precisely, the zeros of the 
corresponding Selberg's zeta function coincide with the poles 
(resonances) of the Laplacian $\Delta_g$ on $X$ \cite{kn:PPe}.


The case of several convex obstacles is generally speaking much more complicated. 
However the case $s_0 > 0$ is easier, since we know  that for $-s_0 \leq \Im( z) \leq 0$ the cut-off
resolvent $R_{\chi}(z)$ is analytic (see \cite{kn:I6}). \\

In the following {\bf we assume that $s_0 < 0.$} The first problem is to 
examine the link between the analyticity of $Z(s)$ for 
$\Re(s) > s_0$ and the behavior of $R_{\chi}(z)$ for $0 \leq \Im( z) < -s_0.$ (The parameters $z$ and $s$ are 
connected by the equality $s = \ii z$). In this direction Ikawa established the following

\begin{thm} {\rm (\cite{kn:I3})} Assume $s_0 < 0$ and $N = 3$. Then for every $\epsilon > 0$ there
exists $C_{\epsilon} > 0$ so that the cut-off resolvent
$R_{\chi}(z)$ is analytic for $\Im (z) < -(s_0 + \epsilon),\: |\Re (z) | \geq C_{\epsilon}\, .$
\end{thm}

A similar result for a control problem has been established by Burq  \cite{kn:B}. The proofs in \cite{kn:I3} and \cite{kn:B} are
based on the construction of an asymptotic solution $U_M(x, s; k)$ with boundary data 
$m(x; k) = e^{\ii k \psi(x)}h(x), k \in \R, \: k \geq 1$, where $\psi$ is a phase function and $h \in C^{\infty} (\Gamma)$ has a small support. 
More precisely, $U_M(., s; k)$ {\text is}\: $C^{\infty}(\Omega)$-\text valued function for $\Re( s) > s_0,$ and we have
\begin{eqnarray}
(\Delta_x - s^2) U_M(., s; k) = 0 \:{\rm for}\: x \in \ovo,\: \Re (s) > s_0,\\
U_M(., s; k) \in L^2(\ovo) \:{\rm if} \:\Re (s) > 0,\\
U_M(x,s; k) = m(x; k) + r_M(x,s; k)\: {\rm on}\:\: \Gamma,\: \Re (s) > s_0,
\end{eqnarray}
where, for $r_M(x, s;k)$ and $\Re (s) > s_0 + d > s_0,\: |s + \ii k| \leq c,$ we have the estimates
\begin{equation}
\|r_M(., s; k)\|_{L^{\infty}(\Gamma)} \leq C_{d, \psi, h} k^{-M}.
\end{equation}
To obtain the leading term of $U_M(x,s;k)$ it is necessary to justify the convergence of series having the form
\begin{equation} \label{eq:1.5}
\sum_{n = 0}^{\infty} \sum_{|\jj| = n+3, j_{n+2} = l} 
\ e^{- s\varphi_{\jj}(x)} 
a_{\jj}(x, s; k),
\end{equation}
where $\jj = (j_0,\ldots,j_{n+2})$ is a configuration (word) of length $|\jj| = n+3$, $\varphi_{\jj}(x)$ are
phase functions and the amplitudes $a_{\jj}(x, s; k)$ depend on the complex parameter
$s \in \C$ and a real parameter $k \geq 1$ (see Sections 3 and 5 for the notation and
more details). These parameters are not connected but to have (1.4) we must take $|s + \ii k| \leq c.$ The main difficulty is to establish the summability of 
above series and to obtain suitable $C^p$ estimates of their traces on $\Gamma$ for $\Re (s) > s_0$. 
The absolute convergence of $Z(s)$ makes it possible to study the absolute
convergence of these series and to get  estimates which  lead to the properties (1.1)-(1.4). 
This might seem a bit surprising since the dynamical zeta function $Z(s)$ is determined by  the periods of periodic
rays and the corresponding Poincar\'e maps, and formally from $Z(s)$ one gets almost no information about the 
dynamics of the rays in a whole neighborhood of the non-wandering set. As it turns out, the absolute convergence of 
$Z(s)$ is a strong condition  which enables us to justify the absolute convergence of (1.5).

The existence of a  domain $\{z \in \C:\: \Re z \in [E- \delta, E + \delta],\:0 \leq \Im\: z \leq h_{\delta}\}$ free of resonances has been 
proved by S. Nonnenmacher and M. Zworski in \cite{kn:NZ} for the operator $-h^2\Delta + V(x),\: V(x) \in C_0^{\infty}(\R^n)$,
assuming that the trapping set of  the Hamiltonian flow $\Phi^t$ of $|\xi|^2 + V(x)$ has a hyperbolic dynamics similar to that of 
the billiard flow in the exterior of $K$.  The existence of a resonance free domain in \cite{kn:NZ} is established under the 
hypothesis Pr$(1/2) < 0$, where Pr$(s)$ is the topological pressure  associated with the (negative infinitesimal) unstable 
Jacobian of the flow $\Phi^t$. In our situation this condition is equivalent to ${\rm Pr}\:(g) < 0,$ where ${\rm Pr}\:(g)$ is the pressure of the function 
$g$ associated with the symbolic dynamics related to the flow (see Sect. 3 for the definition of $g$ and its pressure).
It is shown in Sect. 3 below that  $C_1{\rm Pr}\:(g) \leq s_0 \leq C_2{\rm Pr}\:(g)$ for some constants $C_1 > 0, \: C_2 > 0$, so
${\rm Pr}\:(g) < 0$ if and  only if $s_0 < 0.$ It should be mentioned that the techniques and tools in \cite{kn:NZ} are different 
from those in \cite{kn:I3}, \cite{kn:B} and the present work.

In the case $\Re( s) < s_0$, it is an interesting problem to examine the link between the analytic continuation of $R_{\chi}(z)$ for 
$\Im (z) \geq -s_0$ and that of the dynamical zeta function $Z(s)$.  Several years ago, Ikawa \cite{kn:I5} announced a result 
concerning a {\bf local} analytic continuation of $R_{\chi}(z)$ in a {\bf neighborhood} of a point $z_0$ in the region
$$ {\mathcal D}_{\alpha, \epsilon} =\{z \in \C: \Im (z) \leq -s_0 + |\Re (z)|^{-\alpha}, |\Re (z)| \geq C_{\epsilon}\},\: 0 < \alpha < 1,$$
assuming the following conditions:\\
(i) $Z(s)$ is analytic in a neighborhood of $\ii z_0$ and
\begin{equation} \label{eq:1.6}
|Z(\ii z_0)| \leq |z_0|^{1-\epsilon},\: 0 <\epsilon < 1,
\end{equation}
(ii) if $w(\eta) > 0$ is an eigenfunction of the Ruelle operator
$ L_{-s_0 \tilde{f} + \tilde{g}}$  with eigenvalue 1, then the constants
$$M = \max_{\xi, \:\eta \in \Sigma^{+}_A} \frac{w(\xi)}{w(\eta)}
\quad , \quad m = \min_{\xi \in \Sigma^{+}_A} e^{-s_0 \tilde{f}(\xi) + \tilde{g}(\xi)}$$
satisfy the inequality $\frac{M}{m}\sqrt{\theta} < 1$ with a global constant $0 < \theta < 1$
depending on the expanding properties of the billiard flow (see \cite{kn:I3}, \cite{kn:I4}).
We refer to Sect. 3 for the notation $\Sigma_A^{+}$, $\tilde{f},\: \tilde{g}$.

Ikawa announced in \cite{kn:I5} that $(ii)$ holds in the case of three 
balls centered at the vertices of an equilateral triangle, provided the radii of the balls are sufficiently small. 
In general the condition $(ii)$ is rather restrictive. On the other hand, it is  difficult to check the condition 
$(i)$ if we have no precise information about the spectral properties of $\tilde{L}_s = L_{-s\tilde{f} + \tilde{g}}$ for $\Re ( s)$ close to $s_0$. In \cite{kn:I5} 
there are no comments when (i) holds and whether this happens at all.
As we show in Sect. 5, the estimate (1.6) for $z \in D_{\alpha, \epsilon}$ is related to the behavior of the iterations of  the Ruelle operator  
$\tilde{L}_s$ introduced in  Sect. 3. It does not look like the tools required to do this were available at the time \cite{kn:I5} was written. 
To our knowledge a proof of the result  announced in \cite{kn:I5} has not been published anywhere.

Starting with the work of Dolgopyat \cite{kn:D}, there has been a considerable progress in the analysis of 
the spectral properties of the Ruelle transfer operators $\tilde{L}_s$ related to hyperbolic  systems. The so 
called Dolgopyat type estimates for the norms  of the iterations $\tilde{L}^n_s$ (see \cite{kn:D}, \cite{kn:St2}, 
\cite{kn:St4}) imply an estimate for the zeta function $Z(s)$ in a strip 
$s_0 -\epsilon \leq \Re (s) \leq s_0,\: \epsilon > 0$ (see Sect. 3 and Appendix C below for details).
On the other hand, it is important to note that the information given by the estimates of the iterations 
and the behavior of the spectrum of $\tilde{L}_s$ is richer than that related to the zeta function $Z(s)$. 

Assuming certain regularity of the family of local unstable manifolds $W^u_{\ep}(x)$ of the billiard flow over the non-wandering set 
$\Lambda$  (see Appendix C) and that the  Dolgopyat type estimates (3.3) hold for the related operator $\tilde{L}_s$ for some 
class of functions, in this paper we prove the following main result:

\begin{thm} Let $s_0 < 0$. Suppose that the estimates {\rm (3.3)} for the operator $\tilde{L}_s$ hold and that
the map $\Lambda \ni x \mapsto W^u_{\ep}(x)$  is Lipschitz.
Then there exist $\sigma_1 < s_0$ and $J_1 > 0$ such that the cut-off resolvent $R_{\chi}(z)$ is analytic in
$${\mathcal S} = \{ z \in \C:\: \Im (z) < - \sigma_1,\:  |\Re(z) | \geq J_1\}.$$
Moreover, there exists an integer $m \in \N$ such that 
\begin{equation} \label{eq:1.7}
\|R_{\chi}(z)\|_{L^2(\ovo) \to L^2(\ovo)} \leq C (1 + |z|)^m, \: z \in {\mathcal S}.
\end{equation}
\end{thm}

The geometric assumptions in the above theorem are always satisfied for $N =2$. In particular,
the Dolgopyat type estimates (3.3) stated in Sect. 3 below always hold when $N = 2$ (\cite{kn:St2}).
For $N \geq 3$ it follows from some general results in \cite{kn:St4} that (3.3) hold under
certain assumptions about the flow on $\mt$. These assumptions are listed in detail at the beginning of Appendix C. 
It seems likely that most of these assumptions are either always satisfied or not really necessary in 
proving the  estimates (3.3) for open billiard flows. In fact, it was shown very recently in  \cite{kn:St5} that one of the 
conditions\footnote{This is the non-degeneracy of the symplectic form over the non-wandering
set $\mt$ -- see the condition (ND) in Appendix C below.} imposed in \cite{kn:St4} (and in \cite{kn:PS2} as well) is always 
satisfied for pinched open billiard flows. Apart from that in \cite{kn:St5} a class of examples with $N \geq 3$ is described 
for which the results in this paper can be applied.

Our argument in Sects. 7-8 shows that the integer $m$ in (\ref{eq:1.7}) depends on $\sigma_1$ and  $N$, however
we have not tried to get precise information about $m$.  It seems that to obtain an optimal growth in (\ref{eq:1.7}) is a 
difficult problem.

One should stress that the Dolgopyat type estimates only apply to a special class of  functions on $\mt$, namely
to Lipschitz functions on $\mt$ which are constant on any local stable manifold $\wsloc(x)$ of the billiard flow 
$\phi_t$ (see Sect. 3 below for details). Notice that the estimates for the iterations of the Ruelle operator were originally obtained 
for the Ruelle operator ${\mathcal L}_s$ related to a coding  given by a Markov family of rectangles  
(see \cite{kn:PS2}, \cite{kn:St4} and Appendix C for the notation).
For the proof of Theorem 2 we need Dolgopyat type estimates for the iterations of the  Ruelle 
operator $\tilde{L}_s$ related to the symbolic coding using the connected components of $K$. 
The link between the operators $\tilde{L}_s$ and ${\mathcal L}_s$ and the estimates leading to 
(3.3) are given in Section 3 in \cite{kn:PS2} (see also Proposition 5 in Appendix C).

We should  mention that our result implies the existence of  an analytic continuation of  $R_{\chi}(z)$ in 
a strip $0 \leq \Im(z) \leq -\sigma_1,\: |\Re (z)| > J_1$, without any restrictions on the eigenfunction 
$w(\eta)$ and the behavior of $Z(s)$ for $\sigma_1 \leq \Re(s) \leq s_0.$ 
The estimate \eqref{eq:1.7} enables us to obtain a scattering expansion with an exponential 
decay rate of the remainder for  the solutions of the Dirichlet problem
\begin{equation} \label{eq:1.8}
\begin{cases} (\partial_t^2 - \Delta)u(t,x) = 0,\: x \in \ovo,\: u\vert_{\R \times \Gamma} = 0,\\
u\vert_{t = 0} = f \in C_0^{\infty}(\ovo),\: \partial_t u\vert_{t = 0} = g \in C_0^{\infty}(\ovo).\end{cases}
\end{equation} 
Set ${\mathcal H} = \dot{H}(\ovo)\oplus L^2(\ovo),\: {\mathcal D}^j = H^{j}(\ovo) \oplus H^{j-1} (\ovo),\: j \geq 2,$ where the space $\dot{H}(\ovo)$ is the closure of 
$C_0^{\infty}(\ovo)$ with respect to the norm
$$\|v\|_{\dot{H}(\ovo)} = \Bigl(\int_{\Omega} |\nabla v(x)|^2 dx\Bigr)^{1/2}.$$

\ms

\begin{col}  Let $N$ be odd and let $\chi \in C_0^{\infty}(\R^n)$ be equal to 1 in a neighborhood of $K$. Let $u(t,x)$ be the 
solution of \eqref{eq:1.8} with initial data 
$(\chi f, \chi g).$ Then under the assumptions of Theorem $2$ there exists $L \in \N$ such that for every 
$\epsilon > 0$ and for $t > 0$ sufficiently large we have
$$\chi u(t,x) = \sum_{\Im\: (z_l) \leq -\sigma_1} \sum_{j = 1}^{m(z_l)} w_{z_l, j}(x) e^{\ii t z_l} t^{j-1} + E(t) (f, g),$$
where
$$\|E(t)(f,g)\|_{{\mathcal H}} \leq C_{\epsilon} e^{(\sigma_1 + \epsilon)t}\|(f, g)\|_{{\mathcal D}^L}.$$
Here $\sigma_1 < s_0$ is as in Theorem $2$, $z_l$ are the resonances with $\Im (z_l) \leq -\sigma_1$, $m_l(z_l)$ are the multiplicities of $z_l$ and  
$w_{z_l, j}$ are related to the cut-off resonances states corresponding to $z_l.$
\end{col}

A similar result was established by Ikawa \cite{kn:I3} with $\sigma_1$ replaced  by $s_0 < 0$. Recently, a 
local decay result for the solutions of the wave equation related to hyperbolic convex co-compact manifolds 
${\bf \Gamma} \setminus {\mathbb H}^{n+1}$ was proved by C. Guillarmou and F. Naud \cite{kn:GN}. They 
obtain an exponentially decreasing remainder related to the abscissa $\delta$ of absolute convergence of the Poincar\'e series
$$P_s(m, m') = \sum_{\gamma \in {\bf \Gamma}} e^{-s d_h(m, \gamma m')} ,\: m, m' \in {\mathbb H}^{n +1},$$
$d_h$ being the hyperbolic distance. To improve this result, one would have to establish a polynomial growth of the corresponding 
cut-off resolvent for  $\delta- \epsilon \leq \Re (s) \leq \delta,\: |\Im (s)| \geq C_{\epsilon}$ and small $\epsilon > 0$, and an 
analog of Corollary 1 can be conjectured for  convex co-compact manifolds (for which Dolgopyat type estimates are known). 
For other results concerning scattering expansions for trapping obstacles the reader  could consult \cite{kn:TZ} and the 
references given there.

The proof of Theorem 2 is long and technical. The reason for this is that we are trying to exploit some quite weak information coming from the Dolgopyat 
type estimates for some {\bf restrictive} class of functions defined on a symbolic model to build approximations of the resolvent of a boundary value problem 
based on infinite series which are not absolutely convergent.  This reflects the  geometric situation and we have to deal with infinite series related to reflections 
of trapping rays. In this direction it appears the present work is the first one where infinite series of this kind are used for a WKB construction. 

Below we discuss the main steps in the proof of Theorem 2.

As in \cite{kn:I3}, \cite{kn:I5},  the idea is to construct an approximative solution $U_M(x, s; k)$ for $\sigma_1
\leq \Re (s) \leq s_0,\: |\Im (s)| \geq J_1,\:k \geq 1$, so that $U_M(x, s; k)$ satisfies the conditions (1.1) - (1.3).  For our analysis in 
Sect. 8 we need to study the Dirichlet  problem for $(\Delta_x - s^2)$ with initial data 
$m(x; k) = G(x) e^{\ii k \la x, \eta \ra}\bv_{x \in \Gamma_j} = G(x) e^{\ii k \varphi(x)}\bv_{x \in \gj}$ coming from a representation by using the Fourier transform. On the other hand, it is convenient to pass to data  
$m(x, s; k) = e^{-s \varphi(x)} b_1(x, s; k)$ with $b_1(x, s; k) = e^{(s + \ii k)\varphi(x)} G(x)$ and to work with two parameters $s \in \C$ and $k \geq 1$. After the preparation in 
Sects. 3-5, we construct in Sect. 6  the first approximation $V^{(0)}(x, s; k)$. The first step in the construction of 
$V^{(0)}(x, s; k)$ is the analysis of the series
$$w_{0,j}(x, s; k) = \sum_{n = -2}^{\infty} \sum_{|\jj| = n+3, j_{n+2} = j} e^{- s\varphi_{\jj}(x)} a_{\jj}(x, s; k) 
= \sum_{n = -2}^{\infty}U_{n +2, j}(x, s; k),\: x \in \Gamma_j,$$
where $\jj = (j_0,\ldots,j_n, j_{n+1}, j_{n+2})$ are configurations of length $|\jj| = n +3$, 
$\varphi_{\jj}(x)$ are phase functions and $a_{\jj}(x, s; k)$ are amplitudes 
determined by a recurrent procedure  starting with $m(x, s; k).$ This series corresponds to the sum of the leading terms of the 
asymptotic solutions constructed after an infinite number of reflections. The analysis of $\wj(x,s; k)$ is given in Sects. 3-5. The main 
goal there is to justify the existence of $\wj(x, s; k)$ and to  obtain an analytic continuation of $w_{0, j}(x, s; k)$  from $\Re (s) > s_0$ 
to a strip $\sigma_0 \leq \Re (s) \leq s_0$ with $\sigma_0 < s_0.$ 
To do this, as in the analysis of Dirichlet series with complex parameter, the strategy is to establish suitable estimates for 
$U_{n+ 2, j}(x, s; k)$ and to apply a {\it summation by packages.}  The structure of $U_{n + 2, j}$ is rather complicated since 
the phases  $\varphi_{\jj}(x)$ and the amplitudes $a_{\jj}(x, s; k)$ are related to the 
dynamics of the reflecting rays having $|\jj|$ reflections and issued from the convex front
$\{(x, \nabla \varphi(x)) : \: x \in {\rm supp}\: h\}.$ It seems unlikely that an  explicit relationship exists between 
$U_{n+2, j}(x, s; k)$ and the iterations $L^n_{-s\tf + \tg}$ of  the Ruelle operator $L_{-s\tf + \tg}$ (see Sects. 3 and 5). Consequently,
one would not expect a particular relationship between  $\sum_{n = -2}^{\infty} U_{n +2, j}(x, s; k)$ and the zeta function $Z(s).$
Thus, it appears the situation considered here is rather different from the case of convex co-compact surfaces where it is known 
that the singularities  of the Selberg zeta function coincide with the singularities of  the corresponding Poincar\'e series which in 
turn is related to the resolvent of the Laplacian \cite{kn:PPe}.

It was observed by Ikawa \cite{kn:I5} that $U_{n + 2, j}(x, s; k)$ can be compared with
$L_{-s\tf + \tg}^n {\mathcal M}_{n, s}(x) {\mathcal G}_s \tilde{v}_s(\xi),$
where ${\mathcal M}_{n, s}(x)$ and ${\mathcal G}_s$ are suitable operators defined by means of
billiard trajectories  issued from appropriate unstable or stable manifolds, while 
$\tilde{v}_s(\xi)$ is a function related to the boundary data 
$m(x, s; k) = e^{-s \varphi(x)} h.$ The precise definitions with some small but essential
differences\footnote{In fact, it is difficult to see how the original definitions of the operators
$ {\mathcal M}_{n, s}$ and ${\mathcal G}_s$ in \cite{kn:I5} would work without the changes we have made
in Sect. 3 below.} are given in Sect. 3.  

The crucial step in this direction  is Theorem 3 in Sects. 3-4 below which provides an estimate of the form 
$$\|L_{-s\tf + \tg}^n {\mathcal M}_{n, s}(x) {\mathcal G}_s \tilde{v}_s(\xi) - 
U_{n+2, l}(x, s; k)\|_{C^p(\Gamma)} \leq C_p(s, \varphi, h)(\theta + ca)^n, \: \forall p \in \N,\:\forall n \in \N,$$ 
where $a = s_0 - \Re (s)$ and $c > 0$, $0 < \theta < 1$, $C_p > 0$ are  global constants. The assumption concerning the 
Dolgopyat type estimates (3.3) of  $\tilde{L}_s$  is not required for the proof of Theorem 3. A statement similar to part (a) of 
Theorem 3 (corresponding to $p = 0$) was announced by Ikawa in \cite{kn:I5}, however as far as we know no proof has ever 
been published. The proof of Theorem 3 is long and technical, however we  consider it in detail since it is of fundamental importance 
for the considerations later on.  It is essential to notice that the link between $U_{n+2, j}$ and the iterations of the Ruelle operator 
$L_{_s\tf + \tg}$ is crucial and allows us to find suitable estimates and deduce the convergence of $\wj(x,s; k)$. This could be considered as a 
mathematical interpretation of the interaction between the  terms with complex phases in $U_{n+2, j}$. Sect. 3 
contains the proof of Theorem 3 in the case $p = 0$, while Sect. 4 deals with $p \geq 1$. 

In Sect. 5 we obtain estimates for $w_{0, j}(x, s; k)$  applying Theorem 3. The convergence of $\wj(x, s; k)$ is reduced to the convergence of the series  
$\sum_{n=0}^{\infty}L_{-s\tf + \tg}^n {\mathcal M}_{n, s}(x) {\mathcal G}_s \tilde{v}_s(\xi)$. Here the 
Dolgopyat type estimates (3.3) for the iterations $L^n_{-s\tf + \tg}$  play a crucial role and we can justify the analyticity of 
$w_{0, j}(x, s; k)$ for $\Re (s) \geq \sigma_0$ with $\sigma_0 < s_0.$ 
The estimates of $w_{0, j}(x, s; k)$ for  $\sigma_0 \leq \Re (s) \leq s_0$ are different from those in 
the domain of absolute convergence $\Re (s) > s_0$. 

In Sect. 6 we construct outgoing parametrix $P_h, P_g, P_e$ respectively for the hyperbolic, glancing and elliptic sets of $T^*(\gj)$ related to a fixed strictly convex obstacle $K_j$. We set ${\mathcal S}_j(s) = P_h + P_g + P_e$ and define the first approximation
$$V^{(0)}(x, s; k) = \sum_{j=1}^{\kappa_0} \Bigl({\mathcal S}_j(s) \wj\Bigr)(x,s;k), x \in \Omega,$$
which is an analytic function for $s \in {\mathcal D}_0 = \{s \in \C:\:\sigma_0 \leq \Re (s) \leq 1,\:|\Im (s)| \geq J \geq 2\}.$
Here the estimates for $U_{n+2, j}(x, s; k)$ obtained in Sect. 5 are crucial for the convergence of the series ${\mathcal S}_j(s) \wj$.  Next, we need to examine the leading terms of the traces of $V^{(0)}$ on $\Gamma_{\ell}, \ell \neq j,$ 
and for this purpose we use a microlocal analysis based on the frequency set introduced in \cite{kn:GS} and 
\cite{kn:G} as well as a global construction of asymptotic solution with oscillatory boundary data $e^{-\ii s\varphi_{\jj}(x)} b(x, s; k)$ with frequency set in the  hyperbolic domain given by Ikawa \cite{kn:I3}. Thus, we show that $V^{(0)}(x, s; k)$ satisfies the conditions:
$$\begin{cases}
(\Delta_x - s^2)V^{(0)}(x, s; k) = 0,\: x \in \ovo,\: s \in \do,\\ V^{(0)}(x, s ; k) \in L^2(\ovo) \: {\rm for}\: \Re (s) > 0, \\
V^{(0)}(x, s; k) =  m(x, s; k) + s^{-1} R_1(x, s; k)\: {\rm on}\: \Gamma,\: s \in \do\;,
\end{cases}$$
with estimates
$$\|R_1(x, s; k)\|_{C^p(\Gamma)} \leq C_p \la s +\ii k \ra^{p + 2}\, |s|^{p +(N + 3)/2 + \beta_0}, \: 0 < \beta_0 < 1, \: \forall p \in \N,$$
where $\la z \ra = (1 + |z|)$. The main point here is that $R_1(x, s; k)$ is analytic for $s \in \do$.
Finite higher order approximations $V^{(j)}(x, s; k),\: j = 0,...,M-1,$ are examined in Sect. 7, and we show that
$$\sum_{j = 0}^{M-1} V^{(j)}(x,s; k) = m(x, s; k) + s^{-M} {\mathcal Q}_M(x,s; k),\: x \in \Gamma, \: s \in \do\;,$$
with  estimates
$$\|{\mathcal Q}_M(x,s; k)\|_{C^0(\Gamma)} \leq C_M |s|^{N(M)}\la s + \ii k\ra^{L(M)},\: s \in \do\;,$$
where $N(M) > M$ depends on $M$ and $ L(M) \to \infty$ as $M \to \infty$ and $Q_M(x, s; k)$ is analytic for $s \in \do.$
The situation here is quite different from the absolutely convergent case treated in \cite{kn:I3}, \cite{kn:B}, where we have 
$N(M) = 0$ for $\Re (s) > s_0 + d> s_0$. We need a finite number $M-1 > (N - 3)/2$ of higher order 
approximations, so we fix $M$ and, applying a version of the three lines theorem, we choose $\sigma_1 < s_0$ close to $s_0$ so that for
$$s \in  \{s \in \C:\: \sigma_1 \leq \Re (s) \leq s_0 + c, \: |\Im (s)| \geq J, |s + \ii k| \leq |\sigma_0| + c\},\: s_0 + c \geq 1$$
we get an estimate
$$\|{\mathcal Q}_M(x, s; k) \|_{C^0(\Gamma)} \leq B_M k^{\alpha}$$
with $0 < \alpha < M - \frac{N -1}{2}.$ The final step of our argument 
is in Sect. 8, where we solve  an integral equation on the boundary $\Gamma.$ To do this,
we invert in $L^2(\Gamma)$ an operator $I + Q(s; k)$ and we apply the last estimate to show 
that $Q(s; k)$ has a small  $L^2(\Gamma)$ norm for $k \geq k_1.$\\

Depending on how much details the reader is prepared to see in trying to understand the proof of our main result, we
would suggest {\bf three} different ways to proceed. The first ({\it shortest}) one is to start by reading Sect. 2 and only the 
beginning of Sect. 3 concerning the definitions of $u_{\jj}(x, s)$ and the statement of Theorem 3, however omitting
the proof of this theorem  in Sects. 3-4. Then one should read the definition of $\wj(x, s)$ in Sect. 5, and skipping the proof 
of the estimates (\ref{eq:5.8}) of $\wj$ in Sect. 5, one could go directly to the  constructions in Sect. 6, followed by 
Sects. 7 and 8. The arguments in  Sect. 6-8 use only the estimates (\ref{eq:5.8}) and some geometrical facts from Sect. 2 and 
Appendix B, so the reader should be able to understand the proof of Theorem 2 in Sect. 8 modulo the omitted technical details. 
The second way to proceed is to read Sect. 2 and then to follow the dynamical proofs in Sect. 3, assuming the estimate (3.3).
One could then proceed as above up to Sect. 8. In this way at a first reading Sect. 4 could be skipped, if the reader is not interested 
in the details of the estimates of the derivatives of $U_{n+2, j}$.  
Finally, the third ({\it complete}) way is to read Sect. 2 and then Appendix A and Appendix C in order to understand the estimates (3.3) and the restrictions 
on the class of functions for which we have Dolgopyat type estimates based on \cite{kn:St4} and \cite{kn:PS2}.
Then one could proceed  as in the second way. 

\ms

\noindent
{\it Acknowledgement.} The authors  are very grateful to the referees for their thorough and careful reading of the paper. Their remarks and suggestions lead 
to a significant improvement of the  first version of this paper.

\section{Preliminaries}
\renewcommand{\theequation}{\arabic{section}.\arabic{equation}}
\setcounter{equation}{0}

This section contains some basic facts about the dynamics of the  billiard flow in the exterior $\Omega$ of $K$.
Our main reference is \cite{kn:I3}; see also \cite{kn:B} and \cite{kn:PS1}. The notation follows mainly \cite{kn:I3}.

Throughout the whole paper we use the notation $\con$ and $\Con$ to denote positive {\it global constants} depending only on $K$. These constants might be 
different in different expressions. Notation of the form $\Con_p$, $\con_p$ will be used to denote global constants that  depend on $K$ and possibly the number $p$.

Here and in the rest of the paper we assume that  $K$ is as in Sect. 1.
Denote by $A$ the $\ka_0\times \ka_0$ matrix with entries $A(i,j) = 1$ if $i \neq j$ and
$A(i,i)= 0$ for all $i$, and set
$$\sa = \{ (\ldots, \eta_{-m}, \ldots, \eta_{-1}, \eta_0,\eta_1,  \ldots, \eta_{m},  \ldots ) : 
1\leq \eta_j \leq \ka_0,\:\eta_j \in \N,\:\:\: \eta_{j} \neq \eta_{j+1}\:\: \mbox{\rm for all }\: j \in \Z\}\;,$$
$$\saa = \{  (\eta_0, \eta_1, \ldots, \eta_{m},  \ldots ) : 1\leq \eta_j \leq \ka_0,\:\eta_j \in \N,\:\:\:
\eta_{j} \neq \eta_{j+1}\:\: \mbox{\rm for all }\: j \geq 0\}\;,$$
$$\san = \{ (\ldots, \eta_{-m},  \ldots, \eta_{-1}, \eta_0) : 1\leq \eta_j \leq \ka_0,\:\eta_j \in \N,\:\:\:
\eta_{j-1} \neq \eta_j\:\: \mbox{\rm for all }\: j \leq 0\}\;.$$
Let $\pr_1 : S^*(\Omega) = \Omega\times \sN \longrightarrow \Omega$ and
$\pr_2 : S^*(\Omega)  \longrightarrow \sN$ be the natural {\it projections}. 
Introduce the shift operator $\sigma: \: \sa \longrightarrow \sa$ and $\sigma: \saa \longrightarrow \saa$ by
$(\sigma(\xi))_i = \xi_{i+1}, \: i \in \Z, \: \xi \in \sa$ and $(\sigma(\xi))_i = \xi_{i +1},\: i \in \N,\: \xi \in \saa.$

Fix a {\it large ball} $B_0$ containing $K$ in its interior. For any $x\in  \Gamma = \partial K$ we will denote by
$\nu(x)$ the {\it outward unit normal} to $\Gamma$ at $x$.

For any $\delta > 0$ and $V \subset \Omega$ denote by $S^*_\delta(V)$
the set of those $(x,u)\in S^*(\Omega)$ such
that $x \in V$ and there exist $y\in \Gamma$ and $t \geq 0$ with 
$y+t u = x$, $y+ s u \in \R^N \setminus K$ for all
$s\in (0,t)$ and $\langle u, \nu (y) \rangle \geq \delta$. 

The condition (H) implies the following (see Lemma 3.1 in \cite{kn:I3})
\begin{lem} There exist  constants $\delta_0 > 0$ and $d_0 > 0$ such that for all $i, j = 1,\ldots ,\kappa_0$, if a ray issued from  $x \in \Gamma_i$ with direction $u$ hits 
$\Gamma_j$ at a point $y \in \Gamma_j$ such that $\la u, \nu(y) \ra \geq -\delta_0$, then the forward ray issued from $(y, v)$ with 
$v = u - 2\la u, \nu(y) \ra \nu(y)$ does not meet a $d_0$ neighborhood of $\cup_{\ell \not= j} K_{\ell}.$

\end{lem}

That is, there exists a constant $\delta' > 0$ such that if for some $(y,v) \in S^*(\Omega)$ with
$y \in \Gamma$, both its forward and backward billiard
trajectories have common points with $\Gamma$, then $\delta' \leq \la v, \nu(y)\ra$.

Let $z_0 = (x_0,u_0)\in S^*(\Omega)$. Denote by $X_1(z_0), X_2(z_0), \ldots, X_{m}(z_0), \ldots\;$
the successive {\it reflection points} (if any) of the {\it forward trajectory}  $\gamma_+(z_0) = \{ \pr_1(\phi_{t}(z_0)) : 0 \leq t \}\;.$
If $\gamma_+(z_0)$ is bounded (i.e. it has infinitely many reflection points), we will  say that it {\it has a forward itinerary} $\eta = (\eta_1, \eta_2, \ldots) $ 
(or that it follows the {\it configuration} $\eta$) if $X_j(z_0)\in \partial K_{\eta_j}$ for all $j \geq 1$. Similarly, we will denote by
$\gamma_-(z_0)$ the {\it backward trajectory} determined by $z_0$ and by  $\ldots, X_{-m}(z_0),  \ldots, X_{-1}(z_0), X_0(z_0)$
its backward reflection points (if any). For any $j\in \Z$ for which $X_j(z_0)$ exists denote by
$\Xi_j(z_0)$ the {\it direction} of $\gamma(z_0) = \gamma_-(z_0) \cup \gamma_+(z_0)$ at $X_j(z_0) = \pr_1(\phi_{t_j}(z_0))$, i.e.
$\Xi_j(z_0) = \lim_{t\searrow t_j} \pr_2(\phi_t(z_0))$. Thus, $\phi_{t_j}(z_0) = (X_j(z_0), \Xi_j(z_0))$.
A finite string $\jj =  (j_0,j_1, j_2, \ldots, j_m)$ of numbers $j_i = 1, 2, \ldots, \ka_0$ will be called
an {\it admissible  configuration} (of {\it length} $|\jj| = m+1$) if $j_i \neq j_{i+1}$ for all $i = 0,1,\ldots,m-1$.  We will say that a
billiard trajectory $\gamma$ with successive reflection points $x_0,x_1, \ldots,x_m$
{\it follows the configuration} $\jj$ if $x_i \in \Gamma_{j_i}$ for all $i =0,1,\ldots,m$.

A {\it phase function} on an open set $\uu$ in $\R^N$ is a smooth ($C^\infty$) function $\varphi : \uu \longrightarrow \R$ such that $\| \nabla \varphi \| = 1$ 
everywhere in $\uu$. For $x\in \uu$ the level surface
$$\cc_\varphi(x) = \{ y\in \uu: \varphi(y) = \varphi(x)\}$$
has a unit normal field $\pm \nabla \varphi(y)$. 

\bs

\noindent
{\bf Remark 1.} It should be mentioned that in  Sects. 2-4 the $C^\infty$ smoothness assumption  
can be replaced by $C^k$ for any $k \geq 1$.

\begin{deff} The phase function $\varphi$ defined on $\uu$ is said to satisfy the  {\it condition \rm{(${\mathcal P}$)}} on ${\mathcal V}$ if:

(i) the normal curvatures of $\cc_\varphi$ with respect to the normal field $-\nabla \varphi$ are
non-negative at every point of $\cc_\varphi$;

(ii) $\di \uu^+(\varphi) =  \{ y+t\nv(y) : t\geq 0, y\in \uu\cap {\mathcal V}\} \supset \cup_{i\neq j} K_i$.
\end{deff}

A natural extension of $\varphi$ on $\uu^+(\varphi)$ is obtained by setting
$\varphi( y+t\nv(y)) = \varphi(y) + t$ for $ t\geq 0$ and $y\in \uu\cap {\mathcal V}$.

Given a phase function $\varphi$ satisfying $({\mathcal P})$ on $\Gamma_j$ and $i \neq j$, denote by $\uu_i(\varphi)$ {\it the set of all points $x$ of the form } 
$x = X_1(y,\nv(y)) + t \,\Xi_1(y, \nv(y))$, where $y\in \uu\cap \Gamma_j$ and $t \geq 0$ are such that $X_1(y, \nv(y))\in \Gamma_{i,(j)}$, where
$$\Gamma_{i,(j)} = \{ x\in \Gamma_i : \langle \nu(x), (y-x)/\|y-x\|\rangle  \geq \delta_0 \;\; \mbox{\rm for all }\:\; y\in \Gamma_j\}\;.$$
Then setting $\varphi_i(x) = \varphi(X_1(y,\nv(y))) + t$, one gets a phase function $\varphi_i$ satisfying the
Condition (${\mathcal P}$) on $\Gamma_i$ (\cite{kn:I3}). The operator sending $\varphi$ to $\varphi_i$ is
denoted by $\Phi^i_j$, i.e. $\Phi_j^i (\varphi) = \varphi_i$.

Given an admissible configuration $\jj =  (j_0,j_1,  \ldots, j_m)$ and a phase function $\varphi$
satisfying the Condition (${\mathcal P}$) on $\Gamma_{j_0}$, define
$$\varphi_{\jj} = \Phi_{j_{m-1}}^{j_m} \circ \Phi_{j_{m-2}}^{j_{m-1}} \circ \ldots
\Phi_{j_{1}}^{j_2} \circ \Phi_{j_{0}}^{j_1}\, (\varphi) \;.$$
Notice that for any $z$ in the domain $\uu_{\jj}(\varphi)$ of $\varphi_{\jj}$ there exists 
$(x,u)\in S^*_+ (\Gamma_{j_0})$ such that $x\in \uu$ and $\gamma_+(x,u)$
{\it follows the configuration} $\jj$, i.e. it has at least $m$ reflection points and
$X_i(x,u) \in \Gamma_{j_i}$ for all $i = 1,\ldots,m$, and $z = X_m(x,u) + t\,\Xi_m(x,u)$ 
for some $t \geq 0$. Denote
$$X^{-\ell}(z, \varphi_{\jj}) = X_{m-\ell}(x,u) \:\:\:\:, \:\: 0 \leq \ell \leq m\;.$$

Several well-known facts about the dynamics of the billiard in  $\Omega$, phase functions and related objects will be 
frequently used throughout the paper and for convenience of the reader we state them here.

The following is a consequence of the hyperbolicity of the billiard flow in the exterior of $K$ and can be derived from
the works of Sinai on general dispersing billiards (\cite{kn:Si1}, \cite{kn:Si2}) and from Ikawa's papers on open billiards
(\cite{kn:I3}; see also \cite{kn:B}). In this particular form it can be found in \cite{kn:Sj} (see also Ch. 10 in  \cite{kn:PS1}).\\

\begin{prop}  There exist global constants $\Con > 0$ and $\alpha\in (0,1)$ such that for any
admissible  configuration $\jj =  (j_0,j_1, \ldots, j_m)$ and any two billiard trajectories in $\Omega$
with successive reflection points $x_0,x_1, \ldots,x_m$ and $y_0,y_1, \ldots,y_m$, both  following the configuration $\jj$, we have
$$\| x_i - y_i \| \leq \Con \, (\alpha^i + \alpha^{m-i})  \quad , \quad 0 \leq i \leq m\;.$$
Moreover, $\Con$ and $\alpha$ can be chosen so that if 
there exists a phase function $\varphi$ satisfying the condition {\rm (${\mathcal P}$)} on some open set $\uu$ containing $x_0$ and $y_0$ and such that
$\nabla \varphi(x_0) = (x_1-x_0)/\|x_1-x_0\|$ and $\nabla \varphi(y_0) = (y_1-y_0)/\|y_1-y_0\|$, then $\| x_i - y_i \| \leq \Con \, \alpha^{m-i}$ for $0 \leq i \leq m\;.$
\end{prop}

\ms

Next, given a vector $a = (a_1, \ldots,a_N)\in \R^N$, denote
$\di D_a = a_1 \, \frac{\partial}{\partial x_1} + \ldots + a_N\, \frac{\partial}{\partial x_N}\;,$
and for any $C^1$ vector field $f: U\longrightarrow \R^N$ ($U \subset \R^N$) and any $V\subset U$ set
$\|f\|_{0} (V) = \sup_{x\in V} \| f(x)\|$ and $\|f\|_0 = \|f\|_0(U)$.
Assuming $f$ has continuous derivatives of all orders $\leq p$ ($p \geq 1$), set
$$\|f\|_{p} (x) = \max_{a^{(1)}, \ldots, a^{(p)}\in \S^{N-1}} 
\| (D_{a^{(1)}} \ldots D_{a^{(p)}} f)(x)\| \:\:\:, \:\:\: 
\|f\|_{p}(V) = \sup_{x\in V} \|f\|_{p}(x)\:\:\:, \:\:\: \|f\|_p = \|f\|_p(U)\;,$$
$$\|f\|_{(p)} (x) = \max_{0 \leq j\leq p} \|f\|_j (x) \:\:\:, \:\:\: 
\|f\|_{(p)}(V) = \sup_{x\in V} \|f\|_{(p)}(x)\:\:\:, \:\:\: \|f\|_{(p)} = \|f\|_{(p)}(U)\;.$$
Similarly, for $x \in \Gamma$ and $V \subset \Gamma$ set
$$\|f\|_{\Gamma, p} (x) = \max_{a^{(1)}, \ldots, a^{(p)}\in S_x\Gamma} \| (D_{a^{(1)}} \ldots D_{a^{(p)}} f)(x)\| \:\:\:, \:\:\: 
\|f\|_{\Gamma, p}(V) = \sup_{x\in V} \|f\|_{\Gamma,p}(x)\:\:\:, \:\:\: \|f\|_{\Gamma, p} = \|f\|_{\Gamma, p}(U)\;,$$
where $S_x \Gamma$ is the {\it unit sphere} in the {\it tangent plane}
$T_x\Gamma$ to $\Gamma$ at $x$. Finally, set
$$\|f\|_{\Gamma,(p)} (x) = \max_{0 \leq j\leq p} \|f\|_{\Gamma,j} (x) \:\:\:, \:\:\: 
\|f\|_{\Gamma,(p)}(V) = \sup_{x\in V} \|f\|_{(p)}(x)\:\:\: , \:\:\: \|f\|_{\Gamma, (p)} = \|f\|_{\Gamma, (p)}(U)\;.$$

\noindent
{\bf Remark 2.} It follows easily from the definitions that for 
any $\delta > 0$ and any  integer $p \geq 1$ there exists
a constant $A_p = A_p(\delta,K) > $ such that if $\psi$ is a 
phase function which is at least $C^{p+1}$-smooth on 
some subset $V$ of $\Omega$ and $x\in V\cap \Gamma$ with 
$(x,\nabla \psi(x)) \in S^*_\delta(V)$, then 
$\|\nabla \psi\|_p(x) \leq A_p\, \|\nabla \psi\|_{\Gamma,p}(x)\;.$

\bs

The  following comprises Proposition 5.4 in \cite{kn:I1},
Propositions 3.11 and 3.12 in \cite{kn:I3} and Lemma 4.1 in \cite{kn:I2} 
(see also the proof of the estimate (3.64) in \cite{kn:B}).\\

\begin{prop} For every integer $p \geq 1$ there exist global constants $\Con_p > 0$ and $\alpha\in (0,1)$ such that for any 
admissible configuration $\jj =  (j_0,j_1, \ldots, j_m)$ and any phase functions $\varphi$ and $\psi$ satisfying the 
Condition {\rm (${\mathcal P}$)} on $\Gamma_{j_0}$ on some open set $\uu$, we have
\be
\|\nabla \varphi_{\jj}\|_{p}(x) \leq \Con_p\, \|\nabla \varphi\|_{(p)} (\uu\cap B_0)
\ee
for any $x \in \uu_{\jj}(\varphi)\cap B_0$, and
\be
\|\nabla \varphi_{\jj} - \nabla \psi_{\jj}\|_{p}(x) \leq 
\Con_p\, \alpha^m\, \|\nabla \varphi - \nabla \psi \|_p (\uu\cap B_0)\;,
\ee
\be
\|X^{-\ell}(\cdot, \nabla \varphi_{\jj}) - X^{-\ell}(\cdot, \nabla \psi_{\jj})\|_{\Gamma,p}(x)
\leq \Con_p\, \alpha^{m-\ell}\, \|\nabla \varphi - \nabla \psi \|_{(p)} (\uu\cap B_0)
\ee
for any $x \in \uu_{\jj}(\varphi)\cap \uu_{\jj}(\psi)\cap B_0$ and $0\leq \ell < m$.
Finally, we can choose  $\Con_p > 0$ so that
\be
\|X^{-\ell}(\cdot, \nabla \varphi_{\jj}) \|_{\Gamma,p}(x)
\leq \Con_p\, \alpha^{\ell}
\ee
for all $x \in \uu_{\jj}(\varphi)\cap B_0$ and $0\leq \ell < m$.
\end{prop}
\bs

Given $x$ in the domain $\uu$ of a phase function $\varphi$, denote
$$\Lambda_\varphi (x) = \left( \frac{G_\varphi(x)}{G_\varphi(X^{-1}(x, \nabla \varphi))} \right)^{1/(N-1)} \;,$$
where $G_\varphi(y)$ is the {\it Gauss curvature} of $C_\varphi(y)$ at $y$.  It follows from \cite{kn:I3} (or \cite{kn:B})
that there exist global constants $0 < \alpha_1 < \alpha < 1$ such that
\be
0 < \alpha_1 \leq \Lambda_\varphi(y) \leq \alpha < 1 
\ee
for any phase function $\varphi$ and any  $y \in \uu(\varphi)$.

Now for any $\jj = (j_0= 1, j_1, \ldots, j_m)$ and any $x\in \uu_{\jj} (\varphi)$, slightly changing a definition from \cite{kn:I3}, set
$$(A_{\jj}(\varphi)\, h)(x) = 
\Lambda_{\varphi,\jj}(x)\,  h (X^{-m}(x, \nabla \varphi_\jj))\;, $$
where
$$\Lambda_{\varphi,\jj}(x) = \Lambda_{\varphi_{(j_1, \ldots,j_m)}}(x)\,
\Lambda_{\varphi_{(j_1, \ldots,j_{m-1})}}(X^{-1}(x, \nabla \varphi_{\jj}))\,
\ldots  \Lambda_{\varphi } (X^{-m}(x, \nabla \varphi_{\jj } )) \in (0,1)\;.$$

The following facts can be derived from \cite{kn:I1}, \cite{kn:I3} 
(see also Proposition 5.1 in \cite{kn:B}).

\begin{prop} For every integer $p \geq 1$ there exists a global constant $\Con_p > 0$ 
such that for any admissible configuration $\jj =  (j_0,j_1,  \ldots, j_m)$ and any
phase function $\varphi$ satisfying the Condition (${\mathcal P}$) on $\Gamma_{j_0}$ on some open set $\uu$, we have
$\|\Lambda_{\varphi,\jj}\|_{p}(x) \leq \Con_p\, \|\nabla \varphi\|_{(p)} (\uu\cap B_0)$
for $x\in  \uu_{\jj}(\varphi)\cap B_0\;.$
\end{prop}

\bs

\def\vxij{\varphi_{\xi,j}}
\def\tR{\widetilde{R}}
\def\lip{\mbox{\footnotesize\rm Lip}}
\def\Lip{\mbox{\rm Lip}}
\def\clip{C^{\lip}}

\section{Ruelle operator and asymptotic solutions}
\renewcommand{\theequation}{\arabic{section}.\arabic{equation}}
\setcounter{equation}{0}

Given $\xi \in \sa$, let
$\ldots, P_{-2}(\xi), P_{-1}(\xi), P_0(\xi), P_1(\xi), P_2(\xi) , \ldots$
be the successive reflection points of the unique billiard trajectory in the exterior of $K$ such that
$P_j(\xi) \in K_{\xi_j}$ for all $j \in \Z$. Set 
$$f(\xi) = \| P_0(\xi) - P_1(\xi)\|\;.$$
Following \cite{kn:I3} (see also Appendix A), one constructs a sequence
$\{ \varphi_{\xi,j}\}_{j=-\infty}^\infty$ of phase functions such that for each $j$,
$\vxij$ is defined and smooth in a neighborhood $U_{\xi,j}$ of the segment
$[P_j(\xi), P_{j+1}(\xi)]$ in $\Omega$ and:

\medskip

(i) $\| \nabla \vxij\| = 1$ on $U_{\xi,j}$ and $\nabla \vxij$ satisfies the part $(i)$ of condition (${\mathcal P}$) on $U_{\xi,j}$;

\medskip

(ii) $\di \nabla \vxij (P_j(\xi)) = \frac{P_{j+1}(\xi) - P_j(\xi)}{\| P_{j+1}(\xi) - P_j(\xi)\|}$ ;

\medskip

(iii) $\vxij = \varphi_{\xi,j+1}$ on $\Gamma_{\xi_{j+1}} \cap U_{\xi,j} \cap U_{\xi,j+1}$ ;

\medskip

(iv) for each $x\in U_{\xi,j}$ the surface  $C_{\xi,j}(x) = \{ y\in U_{\xi,j} : \vxij(y) = \vxij(x)\}\;$
is strictly convex with respect to its normal field $\nabla \vxij$.

More precisely, one can proceed as follows. Given $\xi \in \sa$, let  $\xi^- = (\ldots, \xi_{-2}, \xi_{-1}, \xi_0)$
and let $\psi_{\xi^-}$ be the phase function with $\psi_{\xi^-}(P_0) = 0$ and $\nabla \psi_{\xi^-}(P_0) = (P_1-P_0)/\|P_1-P_0\|$
constructed in Proposition 4(a) in Appendix A. Set $\varphi_{\xi,0} = \psi_{\xi^-}$ and $\varphi_{\xi,j} = (\psi_{\xi^-})_{(\xi_0, \xi_1, \ldots, \xi_j)}$
for any $j > 0$. For $j < 0$, setting $\xi^{(j)} = (\ldots, \xi_{j-2},\xi_{j-1}, \xi_j)$  and using again Proposition 4,
we get a phase function $\psi_{\xi^{(j)}}$ with $\psi_{\xi^{(j)}}(P_j) = 0$ and
$\nabla \psi_{\xi^{(j)}}(P_j) = (P_{j+1}-P_j)/\|P_{j+1}-P_j\|$. By the uniqueness of the phase functions
$\psi_\eta$ (see Proposition 4(c)), it follows that there exists a constant $c_j$ such that 
$\psi_{\xi^-} = (\psi_{\xi^{(j)}} + c_j)_{(\xi_j,\xi_{j+1}, \ldots,\xi_0)}$ (locally near the segment 
$[P_0,P_1]$). Setting $\varphi_{\xi,j} = \psi_{\xi^{(j)}} + c_j$, one obtains a phase function defined on some naturally
determined (see the proof of Proposition 4 (a) in Appendix A) open set $\uu_{\xi^-,j}$ such that
\be
(\varphi_{\xi,j})_{(\xi_j, \xi_{j+1}, \ldots, \xi_{-1},\xi_0)} = \psi_{\xi^-},\: j < 0\;.
\ee
This completes the construction of the phase functions $\varphi_{\xi,j}$. 

Moreover, it follows from Proposition 2 that for any $p \geq 1$ there exists a global constant $\Con_p > 0$ such that
\be
\|\nabla \varphi_{\xi,j}\|_{(p)} \leq \Con_p \quad
\ee
for all $\xi\in \sa$ and $ j \in \Z$.\\

\noindent
{\bf Remark 3.} Notice that the above construction can be carried out for $j < 0$ for any $\xi\in \san$  and any
billiard trajectory $\gamma$  in $\Omega$ with reflection points $\ldots, P_{-2}(\xi), P_{-1}(\xi), P_0(\xi)$  such that
$P_j(\xi) \in K_{\xi_j}$ for all $j \leq 0$. Then one defines a phase function  $\psi_{\xi^-}$ with $\psi_{\xi^-}(P_0) = 0$
as above, and using (3.1) one gets a  sequence $\{\varphi_{\xi,j}\}_{j\leq 0}$ of phase functions 
such that for each $j < 0$, $\vxij$ is defined and smooth in a neighborhood $U_{\xi,j}$ of the segment 
$[P_j(\xi), P_{j+1}(\xi)]$ in $\Omega$ and
satisfies the conditions (i)-(iv). Moreover (3.2) holds for any $p \geq 1$ and any $j \leq 0$.\\

\medskip

For any $y\in U_{\xi,j}$ denote by $G_{\xi,j}(y)$ the {\it Gauss curvature} of $C_{\xi,j}(x)$
at $y$. Now define $g : \sa \longrightarrow \R$ by
$$g(\xi) = \frac{1}{N-1} \, \log \frac{G_{\xi,1}(P_1(\xi))}{G_{\xi,0}(P_0(\xi))}\;.$$
Clearly, $g(\xi) = \log \Lambda_{\varphi_{\xi,1}}(P_1(\xi))$, where $\Lambda_\varphi$ is the function 
introduced in Sect. 2.

Given a function $F: \sa \longrightarrow \C$ and an integer $n \geq 0$, set
$$\var_n F = \sup\{|F(\xi) - F(\eta)|:\:\xi_i = \eta_i\: \:{\rm for}\: |i| < n\},$$
and for $0 < \theta < 1$ we define $\|F\|_{\theta} = \sup_n \frac{\var_n F}{\theta^n},\: \||F\||_{\theta} = \|F\|_{\infty} + \|F\|_{\theta}$ 
and introduce the space $\ff_{\theta}(\sa) = \{ F: \: \||F\||_{\theta} < \infty\}.$ 
Clearly $\ff_\theta(\sa)$ is the space of all Lipschitz functions with respect to the metric
$d_\theta$ on $\sa$ defined by $d_\theta(\xi,\xi) = 0$ and $d_\theta(\xi,\eta) = \theta^n$, where
$n \geq 0$ is the least integer with $\xi_i = \eta_i$ for $|i| < n$. 

It follows from Proposition 1 that $f, g \in {\mathcal F}_{\alpha}(\sa)$. By Sinai's Lemma (see e.g. \cite{kn:PP}), there exist 
$\tf,\tg \in \ff_{\sqrt{\alpha}}(\sa)$ depending on future coordinates only and $\chi_1, \chi_2 \in \ff_{\sqrt{\alpha}}(\sa)$ such that
$$f(\xi) = \tf(\xi) + \chi_1(\xi) - \chi_1(\sigma \xi) \quad , \quad 
g(\xi) = \tg(\xi) + \chi_2(\xi) - \chi_2(\sigma \xi),\: \xi \in \sa\;.$$
As in  the proof of Sinai's Lemma, for any $k = 1, \ldots, \ka_0$ choose and fix an arbitrary sequence 
$\eta^{(k)} = (\ldots, \eta^{(k)}_{-m}, \ldots, \eta^{(k)}_{-1}, \eta^{(k)}_0)\in \Sigma_A^-$  with $\eta^{(k)}_0 \neq k$. Then for any $\xi \in \sa$ (or $\xi \in \saa$) set
$$e(\xi) = (\ldots, \eta^{(\xi_0)}_{-m}, \ldots, \eta^{(\xi_0)}_{-1}, \eta^{(\xi_0)}_0= \xi_0, \xi_1, \ldots, \xi_m, \ldots)\in \Sigma_A \;.$$ 
Then we have
$$\chi_1(\xi) = \sum_{n=0}^\infty [ f(\sigma^n (\xi)) - f(\sigma^n e(\xi))]\;,$$
and the function $\chi_2$ is defined similarly, replacing $f$ by $g$.

Setting $\chi(\xi,s)  = -s \chi_1(\xi) + \chi_2(\xi)$, for the function $R(\xi,s) = -s\, f(\xi) + g(\xi) + \i\, \pi$ we have
$R(\xi,s) = \tR(\xi,s) + \chi(\xi,s) - \chi(\sigma \xi,s)$ for  $\xi \in \sa$,  $s\in \C$, where $\tR(\xi,s) =  -s\, \tf(\xi) + \tg(\xi) + \i\, \pi$
depends on future coordinates of $\xi$ only (so it can be regarded as a function on $\saa\times \C$).
Below we need the {\it Ruelle transfer operator} $L_s : C(\saa) \longrightarrow C(\saa)$ defined by
$$L_s  u (\xi)  =  \sum_{\sigma \eta = \xi} e^{\tR(\eta,s)}\,  u (\eta)\;$$
for any continuous (complex-valued) function $ u$ on $\saa$ and any $\xi\in \saa$. Notice that
$$L_s^n u(\xi) = (-1)^n\,\sum_{\sigma \eta = \xi} e^{-s\tf(\eta) + \tg(\eta)}u(\eta) = (-1)^n\, L^n_{-s\tf + \tg}u(\xi) \quad ,\quad  n \geq 0\;,$$
hence $\|L_s^n\|_{\infty} = \big\| L^n_{-s\tf + \tg} \big\|_{\infty}$. Set $\tilde{L}_s = L_{-s\tf + \tg}$.

Define the map $\Phi : \sa \longrightarrow \mtb = \mt \cap S^*_{\dk}(\Omega)$ by 
$$\Phi(\xi) = (P_0(\xi), (P_1(\xi) - P_0(\xi))/ \|P_1(\xi) - P_0(\xi)\|)\;.$$
Then $\Phi$ is a bijection such that $\Phi\circ \sigma = B \circ \Phi$, where $B : \mtb \longrightarrow \mtb$ is the
{\it billiard ball map}. It is well-known (and relatively easy to see) that there exist global constants $0 < \alpha' < \alpha < 1$, $\Con > 0$ and $\con >0$
($\alpha$ is actually the constant from Proposition 1) such that 
$$\con\, d_{\alpha'}(\xi,\theta) \leq \dist(\Phi(\xi), \Phi(\eta)) \leq \Con\, d_\alpha(\xi,\eta) \quad , \quad \xi,\eta\in \sa\;,$$
where $\dist$ is the Euclidean distance in $S^*(\Omega) \subset \R^N\times \S^{N-1}$. Thus, if
$h : \mtb \longrightarrow \C$ is Lipschitz, then $h\circ \Phi\in \ff_{\alpha}(\sa)$, and if $v\in \ff_{\alpha'}(\sa)$, then 
$v\circ \Phi^{-1}$ is a Lipschitz function on $\mtb$.

Let $\pi : \sa \longrightarrow \saa$ be the natural projection. Notice that
for any function $v : \saa \longrightarrow \C$ the function $v\circ \pi : \sa \longrightarrow \C$
depends on future coordinates only, so $(v\circ \pi) \circ \Phi^{-1} : \mtb \longrightarrow \C$ is constant
on local stable manifolds. Conversely, if $h : \mtb \longrightarrow \C$ is constant
on local stable manifolds, then $v = h\circ \Phi : \sa \longrightarrow \C$ depends on future
coordinates only, so it can be regarded as a function on $\saa$. 
For any $(p,u) \in S^*(\Omega)$ sufficiently close to $\mt$, let $\omega(p,u) \in S_{\dk}^*(\Omega)$
be the backward shift of $(p,u)$ along the flow to the first point at the boundary. That is,
$\omega(p,u) = (q,u)\in S_{\dk}^*(\Omega)$, where $p = q + t\, u$ and $(p,u) = \phi_t(q,u)$ for some 
some $t \geq 0$ and $\la u , \nu(q) \ra  > 0$. Thus, $\omega : V_0 \longrightarrow S_{\dk}^*(\Omega)$ is a 
smooth map defined on an open subset $V_0$ of $S^*(\Omega)$ containing $\mt$.

Denote by $\clip_u(\mtb)$ the {\it space of Lipschitz functions} $h : \mtb \longrightarrow \C$ such that 
$h \circ \omega$ is constant on any local stable manifold $\wsloc(x)$ of the flow $\phi_t$ contained in the
interior of $V_0 \setminus  S_{\dk}^*(\Omega)$. For such $h$ 
let $\Lip(h)$ denote the {\it Lipschitz constant} of $h$, and for $t\in \R$, $|t| \geq 1$, define
$$\|h\|_{\lip,t} = \|h\|_0 + \frac{\Lip(h)}{|t|} \quad , \quad \|h\|_0 = \sup_{x\in \mtb} |h(x)|\;.$$

To estimate the norm of $\tilde{L}_s^n$,  we will apply Dolgopyat type estimates (\cite{kn:D}) 
established in the case of open billiard flows in \cite{kn:St2} 
for $N = 2$ and in \cite{kn:St4} for $N \geq 3$ under certain assumptions (see Appendix C below). 
It follows from these results that there exist constants 
$\sigma_0 < s_0, \: t_0 > 1$ and $0 < \rho  < 1$ so that for $s = \tau + \ii t$ with 
$\tau \geq \sigma_0$, $|t|\geq t_0$ and $n = p[\log|t|] + l,\: p \in \N,\: 0 \leq l \leq [\log|t|] - 1$, for
any function $v \in C(\saa)$ of the form $v = h\circ \Phi$ for some $h\in \clip_u(\mtb)$ we have
\begin{equation} \label{eq:6.6}
\|\tilde{L}_s^n v\|_{\infty} \leq C \rho^{p[\log|t|]}e^{l\Pr(-\tau\tf + \tg)}\, \|h\|_{\lip,t} \;.
\end{equation}
Here ${\rm Pr}\,(F)$ denotes the {\it topological pressure} of $F$ defined by
$${\rm Pr}(F) = \sup_{\mu \in {\mathcal M}_{\sigma}} \bigl[h_\mu(\sigma) + \int_{\saa} F\, d\mu ],$$
where ${\mathcal M_{\sigma}}$ is the set of all probability measures on $\saa$ invariant with respect to $\sigma$ and 
$h_\mu(\sigma)$ is the {\it measure-theoretic entropy} of $\sigma$ with respect to $\mu$.

The abscissa of absolute convergence $s_0$ introduced in Sect. 1 is determined by the equality ${\rm Pr}(-s_0 f + g) = 0.$ Thus,
$$h_\nu(\sigma) -s_0 \int f d\nu + \int g d\nu \leq 0 \quad ,\quad \forall \nu \in {\mathcal M}_{\sigma}\,.$$
Let $\nu_g$ be the equilibrium state of $g$ such that ${\rm Pr}\: (g) = h_{\nu_g}(\sigma)  + \int g \, d \nu_g.$ Then 
${\rm Pr}\:(g) \leq s_0 \int f d\nu_g$. Next, let $\nu_0 \in {\mathcal M}_{\sigma}$ be the equilibrium state of $-s_0 f + g$ with
$$h_{\nu_0}(\sigma) - s_0 \int f \, d \nu_0 + \int g \, d \nu_0 = 0\, .$$
This yields $s_0 \int f \, d \nu_0 = h_{\nu_0}(\sigma) + \int g \, d \nu_0 \leq {\rm Pr}\: (g).$ Consequently,
$$\frac{{\rm Pr}\: (g)}{\int f \, d\nu_g} \leq s_0 \leq \frac{{\rm Pr}\: (g)}{\int f \, d \nu_0}$$
and we deduce that $s_0 < 0$ if only if ${\rm Pr}\:(g) < 0$.

We will deal with oscillatory data on $\Gamma_1$ (which can be replaced by any $\Gamma_j$) of the form
$$u_1 (x, s) = e^{- s\, \varphi(x)}\, h (x) \quad , \quad x\in \Gamma_1\;, s\in \C,\: \sigma_0 \leq \Re (s) \leq 1 \;.$$
Here $\varphi$ is a $C^{\infty}$ phase function defined on some open subset $\uu = \uu(\varphi)$
and satisfying the condition $(\mathcal P)$ on $\Gamma_1$ (see Sect. 2 above)
and $h$ is a $C^{\infty}(\Gamma)$ function with small support on $\Gamma_1$.  In fact, using a $C^\infty$ 
extension, we may assume that $h$ is a $C^\infty$ function on $\R^N$, so in particular $h$ is $C^\infty$
on $\uu$, as well.
For every configuration 
$\jj = (j_0, j_1,\ldots,j_m), \: j_0 =1,\: |\jj| = m + 1$, we can construct a function $u_{\jj}(x, s)$ following a recurrent procedure (see \cite{kn:I5}). 
We construct a sequence of phase functions $\varphi_{\jj}(x)$ and amplitudes $a_{\jj}(x)$ and define 
$$u_{\jj}(x, s) = (-1)^{|\jj|-1}e^{- s \varphi_{\jj}(x)} a_{\jj}(x)\;.$$
For the configurations $\jj$ and $\jj' = (j_0, j_1,\ldots, j_m, j_{m+1})$ we have 
$u_{j_0}(x, s) = u_1(x, s)$ on $\Gamma_1$ and  $u_{\jj}(x, s) + u_{\jj'}(x, s) = 0$ on $\Gamma_{j_{m+1}}$.

The phase functions $\varphi_{\jj}$ and their domains $\uu_{\jj} (\varphi)$ are determined following the procedure in Section 2.
In particular, each $\varphi_{\jj}$ satisfies the condition (${\mathcal P}$) on $\Gamma_{j_m}$, so it follows from the definition of the condition (${\mathcal P}$)
(see (ii) there) that  $\Gamma_{i} \subset \uu_{\jj}(\varphi)$ for every $i = 1, \ldots, \kappa_0$, $i \neq j_m$.
The amplitudes $a_{\jj}(x)$ are determined  on $\uu_{\jj}(\varphi)$ as the solutions of the transport equations
$$2 \la \nabla \varphi_{\jj}, \nabla a_{\jj} \ra + (\Delta \varphi_{\jj})a_{\jj} = 0.$$ 
More precisely, using the notations of Sect. 2 (see also Sect. 4 in  \cite{kn:I3} and Sect. 4.1 in \cite{kn:I5}), we will assume that $a_{\jj}(x)$ has the form
\be
a_{\jj}(x) = (A_{\jj}(\varphi)h) (x) \quad , \quad x\in \uu_{\jj}(\varphi) \;.
\ee

Next, let  $\mu = (\mu_0 = 1, \mu_1, \ldots) \in \saa$. It follows from \cite{kn:I3} that there exists a unique point  $y(\mu) \in \Gamma_1$
such that the ray $\gamma(y, \varphi)$ issued from a point $y(\mu)$ in direction
$\nabla\varphi(y(\mu))$ follows the configuration $\mu$. Let $Q_0(\mu) = y(\mu), Q_1(\mu), \ldots, $ 
be the consecutive reflection points of this ray. Define
$$f_i^+(\mu) = \| Q_i(\mu) - Q_{i+1}(\mu)\|\quad , \quad
g_i^+(\mu) = \frac{1}{N-1}\, \log \frac{ G^\varphi_{\mu,i}(Q_{i+1}(\mu))}{G^\varphi_{\mu,i}(Q_i (\mu))} < 0\;,$$
where  $G^\varphi_{\mu,i}(y)$ denotes the {\it Gauss curvature} of the surface
$$C^\varphi_{\mu,i} (x) = \{ z\in \uu_{(\mu_0, \mu_1, \ldots, \mu_i)}(\varphi) :
\varphi_{(\mu_0, \mu_1, \ldots, \mu_i)}(z)  =  \varphi_{(\mu_0, \mu_1, \ldots, \mu_i)}(x)\}\;$$
at $y$.  As for $g(\xi)$, the function $g_i^+(\mu)$ can be expressed by means of the function $\Lambda_\varphi$
introduced in Sect. 2, namely $g_i^+(\mu) = \log \Lambda_{\varphi_{(\mu_0, \mu_1, \ldots, \mu_i)}}(Q_{i+1}(\mu))$.

Using the points $Q_j(\mu)$ constructed above, define $\tv \in \ff_\theta(\saa)$ by
$$\tv_s(\xi) =  e^{-s\, \varphi(Q_0(\xi))}\, h (Q_0(\xi))$$
if $\xi_0 = 1$ and $\tv_s(\xi) = 0$ otherwise. Here the function $h$ comes from the boundary data $u_1(x, s)$.

Next, for $s\in \C$ and $\xi\in \saa$ with $\xi_0 = 1$, following \cite{kn:I5}, set
\be
\phi^+(\xi,s) = \sum_{n=0}^\infty \left( -s \,  [f(\sigma^n e(\xi)) - f^+_n(\xi)] 
+ [g(\sigma^n e(\xi)) - g^+_n(\xi)]\right) \;.
\ee
Formally, define $\phi^+(\xi,s) = 0$ when $\xi_0 \neq 1$, thus obtaining a function
$\phi^+ : \saa \times \C \longrightarrow \C$.  

Now for any  $s\in \C$ define the operator $\gg_{s} : C (\saa) \longrightarrow C(\saa)$ by 
$$(\gg_{s} v) (\xi) =  \sum_{\sigma \eta = \xi} e^{-\phi^+(\eta,s) - s \tf(\eta) + \tg(\eta)}\, v(\eta)
\quad , \quad v\in C(\saa)\;, \; \xi\in \saa\;.$$
(Although similar, this is different from the corresponding definition in \cite{kn:I5}.) 

{\bf Fix an arbitrary $\ell = 1, \ldots, \ka_0$ and an arbitrary point $x_0 \in \Gamma_\ell$}.
Define the function $\phi^-(x_0;\cdot, \cdot) : \sa \times \C \longrightarrow \C$ (depending on $\ell$ as well)
as follows. First, set $\phi^-(x_0;\eta,s) = 0$ if $\eta_0 \neq \ell$. Next, assume that $\eta\in \sa$ satisfies $\eta_0 = \ell$.
There exists a unique billiard trajectory in $\Omega$ with successive reflection points $ \tP_i (x_0;\eta) \in \dk_{\eta_i}$ 
($-\infty < i \leq 0$) such that $x_0 = \tP_{-1}(x_0;\eta) + t\nabla \psi_{\eta^-}(\tP_{-1}(x_0;\eta))$ for some $t > 0$. 
(See the beginning of this section and Appendix A for the definition of $\psi_{\eta^-}$.)
Notice that in general the segment $[\tP_{-1}(x_0;\eta),x_0]$ may intersect the interior of $K_\ell$. 
Denote $\tP_0(x_0;\eta) = x_0$, and for any $i < 0$ set
$$f^-_i(x_0;\eta) = \| \tP_{i+1}(x_0;\eta) - \tP_{i}(x_0;\eta)\|\quad , \quad
g^-_i(x_0;\eta) = \frac{1}{N-1}\, \log \frac{G_{\eta,i}(\tP_{i+1}(x_0;\eta))}{G_{\eta,i}(\tP_{i}(x_0;\eta))}\;.$$
Then define
$$\phi^-(x_0;\eta,s) = -s\, \sum_{i = -1}^{-\infty} [f(\sigma^i(\eta)) - f^-_i(x_0;\eta)] 
+ \sum_{i = -1}^{-\infty} [g(\sigma^i (\eta)) - g^-_i (x_0;\eta)]\;.$$
We will show later that this series is absolutely convergent.

Next, define the operator $\mm_{n,s}(x_0) : C(\saa) \longrightarrow C(\saa)$ (depending also on $\ell$) by
$$(\mm_{n,s}(x_0) v)\, (\xi) = \sum_{\sigma \eta = \xi} e^{-\phi^-(x_0;\sigma^{n+1}e(\eta),s) -
\chi(\sigma^{n+1} e(\eta),s) - s \tf(\eta) + \tg(\eta)} \, v(\eta) $$
for any $v\in C(\saa)$, any $x_0\in \Gamma$  and any $\xi \in \saa$.

Let $s_0 \in \R$ be the abscissa of absolute convergence of the dynamical zeta function (see Sect. 1)
determined by ${\rm Pr}\, (-s_0 \tf + \tg) = 0.$

The first part  in the following theorem is similar to (4.10) in \cite{kn:I5}:

\begin{thm} There exist global constants $c > 0$, $a > 0$, $\theta \in (0,1)$ and $\Con_p > 0$ for every integer $p \geq 0$  
such that for any choice of $\ell = 1, \ldots, \ka_0$ and $x_0\in \Gamma_\ell$ the following hold:
\ms

(a)  For all integers $n \geq 1$, all $\xi\in \saa$ with $\xi_0 = \ell$ and all $s \in \C$ with $\Re (s) \geq s_0 - a$ we have
\begin{eqnarray}
&        & \left| \left(L^n_s \mm_{n,s}(x_0) \gg_{s} \tv_s\right)(\xi) - \sum_{|\jj | = n+3, j_{n+2} = \ell} 
u_{\jj}(x_0, s)\right|\nonumber \\
& \leq &  \Con_0\,(\theta+ c \, a)^n\, e^{\Con_0 [ \Ref(s)\, (1+\|\varphi\|_{\Gamma,0}) + \|\nabla  \varphi\|_{\Gamma,(1)}]}
\left[\Bigl(|s|+ \|\nabla \varphi\|_{\Gamma, (1)}\Bigr)\,  \|h\|_{\Gamma,0} + 
\|h\|_{\Gamma, (1)}\right]\; .
\end{eqnarray}

(b) For  all $n \geq 1$, all $\xi\in \saa$ with $\xi_0 = \ell$ and  all $s \in \C$ with $\Re (s) \geq s_0 - a$  we have
\begin{eqnarray}
&        & \left\| \left(L^n_s \mm_{n,s}(\cdot) \gg_{s} \tv_s\right)(\xi) - \sum_{|\jj | = n+3, j_{n+2} = \ell } 
u_{\jj}(\cdot,\: s)\right\|_{\Gamma,p}\nonumber \\
& \leq &  \Con_p\,(\theta+ c \, a)^n\, e^{\Con_p [ |\Ref(s)|\, (1+\|\varphi\|_{\Gamma,0}) + \|\nabla  \varphi\|_{\Gamma,(1)}]}
\sum_{i = 0}^{p} \Bigl(|s| \|\nabla \varphi\|_{\Gamma, i} + \|\nabla
\varphi\|_{\Gamma, i+1} \Bigr)^{i + 1}  \|h\|_{\Gamma, p- i}\, .
\end{eqnarray}

\end{thm}

In this section we deal with part (a). The proof of part (b) is given in Section 4 below.\\

{\it Proof of Theorem 3(a).}
Fix $\ell$, $x_0 \in \Gamma_\ell$ and $\xi\in \saa$ with $\xi_0 = \ell$. Then for any $s\in \C$ and $n \geq 1$,
using Sect. 4.1 in \cite{kn:I5}, setting $\jj = (1,j_1,j_2,\ldots, j_{n+1},  \ell)$, we get
\be
u_{(1,j_1,j_2,\ldots, j_{n+1},  \ell)}(x_0, s) = 
(-1)^{n+2}\, e^{- s\, [\varphi (Q_0(\jj)) + f^+_0(x_0;\jj) + \ldots + f^+_{n+1}(x_0;\jj)]}\, a_{\jj} (x_0)\;,
\ee
where $f^+_i(x_0; \jj) = \| Q_i(x_0;\jj) - Q_{i+1}(x_0;\jj)\|$ ($i = 0,1, \ldots, n+1$), $Q_i(x_0;\jj)$ being the
reflection points of the billiard trajectory issued from a point $y \in \Gamma_1$ in direction $\nabla \varphi(y)$
which follows the configuration $\jj$ for its first $n+1$ reflections and is such that  $Q_{n+2}(x_0;\jj) = x_0$. 
Notice that the segment $[Q_{n+1}(x_0;\jj), x_0]$ may intersect\footnote{In fact one can define the functions $f^+_i(x_0; \jj)$ ($i = 0,1, \ldots, n+1$) and therefore
$u_{\jj}(x_0, s)$ for any $x_0\in \uu_{\jj}(\varphi)$ in a similar way. Just consider the (unique) billiard trajectory issued from
a point $y = Q_0(x_0; \jj) \in \Gamma_1$ in direction $\nabla \varphi(y)$ following the configuration 
$\jj$ for its first $n+1$ reflections and  such that if $v$ is the reflected direction
of the trajectory at $Q_{n+1}(x_0;\jj)$, then $x_0 = Q_{n+1}(x_0,\jj) + t\, v$ for some $t \geq 0$.} the interior of $K_\ell$. 
Then there is exactly one such trajectory. Given a function  $F(\xi): \Sigma_A^+ \longrightarrow \C$, introduce the notation
$$F_{n}(\xi) = F(\xi) + F(\sigma(\xi)) +...+F(\sigma^{n-1}(\xi)).$$

We have
\begin{eqnarray}
\left(L^n_s \mm_{n,s} (x_0)\gg_{s} \tv_s\right)(\xi) 
& = & (-1)^n\, \sum_{\sigma^n\eta = \xi} 
e^{-s\tf_n(\eta) + \tg_n(\eta)}\, \left(\mm_{n,s} (x_0) \gg_{s} \tv_s\right)(\eta)\nonumber\\
& = &  (-1)^n\, \sum_{\sigma^n\eta = \xi} e^{-s\tf_n(\eta)+ \tg_n(\eta)}\, \sum_{\sigma\zeta = \eta} 
e^{-\phi^-(x_0; \sigma^{n+1} e(\zeta),s)- \chi(\sigma^{n+1} e(\zeta),s) - s \tf(\zeta) +\tg(\zeta)}\nonumber\\ 
&    & \quad \times \sum_{\sigma\mu = \zeta} e^{-\phi^+(\mu,s) +  \chi(e(\mu),s) -s \tf(\mu)+ \tg(\mu)}\,  \tv_s(\mu)\nonumber\\
& = & (-1)^n\, \sum_{\sigma^{n+2}\mu = \xi, \mu_0 = 1}  e^{ -s \tf_{n+2}(\mu) + \tg_{n+2}(\mu)}\,  \Wn(x_0;\mu,s)\;, 
\end{eqnarray}
where the function
$$\Wn(x_0;\cdot, \cdot) = \Wn_{1,\ell}(x_0;\cdot, \cdot) :  \saa\times \C \longrightarrow \C$$
is defined by $\Wn(x_0;\mu,s) = 0$ when $\mu_0 \neq 1$ or $\mu_{n+2} \neq \ell$ and
\begin{eqnarray}
\Wn(x_0; \mu,s)
=  e^{ - \phi^-(x_0; \sigma^{n+1} e(\sigma \mu),s) -
\chi(\sigma^{n+1} e(\sigma \mu),s)  - \phi^+(\mu,s)
+ \chi(e(\mu),s)  -s\, \varphi(Q_0(\mu))}\, \,  h (Q_0(\mu))\;
\end{eqnarray}
whenever $\mu_0 = 1$ and $\mu_{n+2} = \ell$.  It follows from (3.9)  that
\be
\left[ L^n_{s} \mm_{n,s}(x_0) \gg_{s} \tv_s\right] (\xi) 
= (-1)^n\, \left[ L_{-s \tf + \tg}^{n+2} \left( \Wn(x_0;\cdot ,s)\right) \right](\xi)\;.
\ee

Clearly, in (3.9) the summation is over sequences
\be
\mu = (1,j_1,j_2, \ldots, j_{n+1}, \ell, \xi_1, \xi_2, \ldots) = (\jj ,\xi)\;,
\ee
with $\mu_{n+2} = \ell$, where $\jj  = (1,j_1,j_2,\ldots,j_{n+1}, \ell)$. 

Write for convenience
\be
\Wn(x_0;\mu,s) = e^{z (x_0;\mu,s)}\, e^{-s\, \varphi(Q_0(\mu))}\,  h (Q_0(\mu))\;,
\ee
where
\be
z (x_0;\mu,s) = - \phi^-(x_0;\sigma^{n+1} e(\sigma \mu),s) - \chi(\sigma^{n+1} e(\sigma \mu),s)  - \phi^+(\mu,s)
+ \chi(e(\mu),s) \;.
\ee

It follows from Propositions 1 and 2 that there exist global constants $C > 0$ and $\alpha \in (0,1)$ such that
$$|f(\sigma^n e(\xi)) - f^+_n(\xi)| \leq C \, \alpha^n \quad, \quad
|g(\sigma^n e(\xi)) - g^+_n(\xi)| \leq C\, \|\nabla \varphi \|_{\Gamma,(1)}\, \alpha^n\;$$
for all $\xi \in \sa$ and all integers $n \geq 1$, so by (3.5),
\begin{eqnarray*}
\phi^+(\mu,s) 
& = & (|s|  + \|\nabla \varphi\|_{\Gamma,(1)})\, O(\alpha^{n}) +
\sum_{i = 0}^{n+1} \left(-s\, [f(\sigma^i e(\mu)) - f^+_i(\mu)] + [g(\sigma^i e(\mu)) - g^+_i (\mu)] \right)\;.
\end{eqnarray*}
Thus, using the definitions of $\tf$, $\tg$ and $\chi$ and the fact that
$\chi(\sigma^{n+2} e(\mu),s) = \chi(\sigma^{n+1} e(\sigma \mu),s) + |s|\, O(\alpha^n)$, we get
\begin{eqnarray*}
&    & -s [f^+_0(\mu) + f^+_1(\mu) + \ldots +  f^+_{n+1}(\mu)]
+ [g^+_0(\mu) + g^+_1(\mu) + \ldots +  g^+_{n+1}(\mu)] \\
& = &  (s + \|\nabla \varphi \|_{\Gamma,(1)})\, O(\alpha^{n}) -\phi^+(\mu,s) -s[f(e(\mu))+f(\sigma e(\mu)) +\ldots +  
f(\sigma^{n+1} e(\mu)]\\
&    & \quad + [g(e(\mu))+g(\sigma e(\mu)) +\ldots +  g(\sigma^{n+1} e(\mu)]\\
& = &  (|s| + \|\nabla \varphi \|_{\Gamma,(1)})\, O(\alpha^{n}) -\phi^+(\mu,s) -s \tf_{n+2}(\mu) + 
\tg_{n+2}(\mu)+ \chi(e(\mu),s) - \chi (\sigma^{n+1} \, e(\sigma \mu),s)\;.
\end{eqnarray*}

Now, {\bf fix for a moment $n \geq 1$ and $\mu$ as in (3.12)}, and set $\eta = \sigma^{n+1} e(\sigma(\mu))$. Then  we have
\be
\eta = \sigma^{n+1} e(\sigma(\mu)) = 
( \ldots, *, *, \mu_1, \mu_2, \ldots, \mu_{n+1}; \mu_{n+2} = \ell , \mu_{n+3}, \ldots)\;,
\ee
and as for $\phi^+$ one gets 
$$\phi^-(x_0;\eta,s) =  ( |s| + \|\nabla \varphi \|_{\Gamma,(1)})\, O(\alpha^{n}) - s\; 
\sum_{i = -1}^{-n-1} [ f(\sigma^{i} \eta) - f^-_i (x_0;\eta)] + \sum_{i = -1}^{-n-1} [ g(\sigma^{i} \eta) - g^-_i (x_0;\eta)]\;.$$
From these estimates  and (3.14) one derives that
\begin{eqnarray}
z (x_0;\mu,s)  
& = & s\tf_{n+2}(\mu) - \tg_{n+2}(\mu) -\phi^-(x_0;\eta,s)  -s\,  \sum_{i = 0}^{n+1} f^+_i(\mu) +   \sum_{i=0}^{n+1}
g^+_i(\mu)  +  ( |s| + \|\nabla \varphi\|_{\Gamma,(1)})\, O(\alpha^{n}) \nonumber \\
& = &  s\tf_{n+2}(\mu) - \tg_{n+2}(\mu) -s \, c (x_0;\mu) + d (x_0;\mu) +   (|s| + \|\nabla \varphi\|_{\Gamma,(1)})\, O(\alpha^{n})\;,
\end{eqnarray} 
where 
$$c (x_0;\mu) = - \sum_{i =0}^{n+1} [ f(\sigma^{i} \eta) - f^-_i (x_0;\eta)]
+  \sum_{i = 0}^{n+1} f^+_i (\mu) \quad , \quad
d (x_0;\mu) =   -\sum_{i = -1}^{-n-1} [g(\sigma^{i}\eta ) - g^-_i (x_0;\eta)] +\sum_{i = 0}^{n+1} g^+_i (\mu) \;.$$

We will show that
\be
\left| c (x_0;\mu) -  \sum_{i = 0}^{n+1} f^+_i (x_0;\jj )\right|  \leq \Con\, \alpha^{n}\;, 
\ee
and
\be
\left| e^{d (x_0;\mu)}\, h (Q_0(\mu)) -  (A_{\jj }(\varphi) h) (x_0) \right| \leq \Con \, 
\left( \|\nabla \varphi \|_{\Gamma,(1)}\, \|h\|_{\Gamma,0} +  \|h\|_{\Gamma,(1)}\right)\, \theta^{n} 
\ee
for some global constant $\Con > 0$, where
$$\theta = \sqrt{\alpha} \in (0,1)\;.$$

There exists a unique ray $\gamma(y,\varphi)$ issued from a point $y = y_n(x_0;\mu) \in \Gamma_1$ in direction
$\nabla\varphi(y)$, following the configuration $\mu$ for its first $n+1$ reflections and such that if 
$\tQ_i (x_0;\mu)$ ($1 \leq i \leq n+1$) are its first $n+1$ reflection points and $v$ is the reflected direction
of the trajectory at $Q_{n+1}(x_0;\jj)$, then $x_0 = Q_{n+1}(x_0,\jj) + t\, v$ for some $t \geq 0$.
Set $\tQ_{n+2}(x_0;\mu) = x_0$. Notice that as before the segment $[\tQ_{n+1}(x_0;\mu) , x_0]$ may intersect the
interior of $K_{\ell}$ (or be tangent to $\Gamma_\ell$ at $x_0$).

Before we continue, let us make a few simple (however essential) remarks concerning the  sequences of points 
\be
Q_0(\mu) \in \Gamma_1 = \Gamma_{\mu_0}, Q_1(\mu) \in \Gamma_{\mu_1},
\ldots, Q_{n+1}(\mu) \in \Gamma_{\mu_{n+1}}, Q_{n+2}(\mu) \in \Gamma_{\mu_{n+2}} = \Gamma_\ell , \ldots\;,
\ee
\be
\tQ_0(x_0;\mu) \in \Gamma_1 = \Gamma_{\mu_0}, \tQ_1(x_0;\mu) \in \Gamma_{\mu_1},
\ldots, \tQ_{n+1}(x_0;\mu) \in \Gamma_{\mu_{n+1}}, \tQ_{n+2}(x_0;\mu) \in \Gamma_\ell\;,
\ee
\be
\ldots, P_{\eta_{-n-1}}(\eta) \in \Gamma_{\eta_{-n-1}}=\Gamma_{\mu_1},
\ldots, P_{-1}(\eta) \in \Gamma_{\eta_{-1}} = \Gamma_{\mu_{n+1}}, 
P_{0}(\eta) \in \Gamma_{\eta_0} = \Gamma_{\mu_{n+2}} = \Gamma_\ell, \ldots\;,
\ee
\be
\ldots, \tP_{\eta_{-n-1}}(x_0;\eta) \in \Gamma_{\eta_{-n-1}}=\Gamma_{\mu_1},
\ldots, \tP_{-1}(x_0;\mu) \in \Gamma_{\eta_{-1}} = \Gamma_{\mu_{n+1}}, 
\tP_{0}(x_0;\eta) \in \Gamma_{\eta_0} = \Gamma_{\mu_{n+2}}  = \Gamma_\ell\;.
\ee
It is clear that the sequences (3.19) and (3.20) 'start' from the same convex level surface $\varphi = \const$, therefore by Proposition 1 
there exist constants $\Con > 0$ and $\alpha \in (0,1)$ such that
\be
\| Q_i(\mu) - \tQ_i (x_0;\mu)\| \leq \Con\, \alpha^{n+2-i}\quad, \quad 0 \leq i \leq n+2\;.
\ee
(Notice that $\tQ_{n+2}(x_0;\mu) = x_0\in \Gamma_\ell$, so $\|Q_{n+2}(\mu)-\tQ_{n+2}(x_0;\mu)\|\leq \diam (K) \leq \Con$.)
Similarly, the right ends of sequences (3.21) and (3.22) determine points on the same unstable manifold of the billiard flow $\phi_t$,
so by Proposition 1 these sequences `converge backwards', i.e.
\be
\| P_{i}(\eta) - \tP_i (x_0;\eta)\| \leq \Con\, \alpha^{|i|}\quad, \quad i \leq 0\;.
\ee
On the other hand, notice that the sequences (3.19) and (3.21) continue indefinitely
to the right following the same patterns. Thus, these sequences `converge forwards'. More
precisely, using  Proposition 1 again, we have
\be
\| Q_i (\mu) - P_{i -n-2}(\eta)\| \leq \Con\, \alpha^{i}\quad, \quad 1 \leq i \;.
\ee
Similarly, the sequences (3.20) and (3.22) `converge forwards' to 
$\tQ_{n+2}(x_0;\mu) = \tP_0(x_0;\eta) = x_0$, namely
\be
\| \tQ_i (x_0;\mu) - \tP_{i -n-2}(x_0;\eta)\| \leq \Con\, \alpha^{i}\quad, \quad 1 \leq i  \leq n+2\;.
\ee

It now follows from (3.2) and (3.24)  that, 
\be
|g(\sigma^i (\eta)) - g^-_i (x_0;\eta)| = \left| \frac{1}{N-1} \, \log \frac{G_{\eta,i}(P_{i+1}(\eta))}{G_{\eta,i}(P_i (\eta))}
- \frac{1}{N-1} \, \log \frac{G_{\eta,i}(\tP_{i+1}(x_0;\eta))}{G_{\eta,i}(\tP_i (x_0;\eta))}\right| \leq  \Con\, \alpha^{|i|}\;
\ee
for all $i \leq 0$. In particular, the second series in (3.5) is absolutely
convergent, and by (3.27) and Proposition 3,
$|d (x_0;\mu) | \leq \Con$ for some global constant $\Con > 0$.

Next, setting
\be
\ta_i (x_0;\mu) = 
\frac{1}{N-1}\, \log\left( \frac{G^\varphi_{\mu,i }(\tQ_{i+1}(x_0;\mu))}{G^\varphi_{\mu,i}(\tQ_{i}(x_0;\mu))} \right) \;,
\ee
and using (3.23) and Proposition 2, one gets
\begin{eqnarray}
| \ta_i (x_0;\mu) - g^+_i (\mu) |
& =     & \frac{1}{N-1}\, 
\left| \log \frac{G^\varphi_{\mu,i}(\tQ_{i+1}(x_0;\mu))}{G^\varphi_{\mu,i}(\tQ_{i}(x_0;\mu))}
- \log \frac{ G^\varphi_{\mu,i}(Q_{i+1}(\mu))}{G^\varphi_{\mu,i}(Q_i (\mu))} \right|\nonumber\\
& \leq & C\, \|\nabla \varphi\|_{\Gamma,(1)}\, (\|\tQ_{i}(x_0;\mu) - Q_i (\mu)\| + 
\|\tQ_{i+1}(x_0;\mu) - Q_{i +1}(\mu)\|)\nonumber\\
& \leq & \Con \, \|\nabla \varphi\|_{\Gamma,(1)}\, \alpha^{n + 2 - i}\;.
\end{eqnarray}
for all $i = 0, 1, \ldots, n+2$.

Next, notice that by construction 
$\varphi_{\eta,i } = (\varphi_{\eta,-n-2})_{(\mu_1, \ldots,\mu_{n+2+i})} + \const$ for $ -n-1\leq i \leq -1\;.$
Thus, by (2.2), (3.2)  and (3.25), for all $-n-1 \leq i \leq -1$ we have
\begin{eqnarray}
| g^+_{n+2+i}(\mu) - g(\sigma^i \eta) |
& =     & \frac{1}{N-1}\, \left| 
\log \frac{ G^\varphi_{\mu,n+2+i}(Q_{n+2+i+1}(\mu))}{G^\varphi_{\mu,n+2+i}(Q_{n+2+i}(\mu))} -
\log \frac{G_{\eta,i}(P_{i+1}(\eta))}{G_{\eta,i}(P_i(\eta))}\right|\nonumber\\
& \leq & \Con\, (\|\nabla \varphi_{(\mu_1, \ldots,\mu_{n+2+i})} - 
\nabla (\varphi_{\eta,-n-2})_{(\mu_1, \ldots,\mu_{n+2+i})}\|_{\Gamma, (1)}\nonumber\\
&        & \quad + \|Q_{n+2+i+1}(\mu)) - P_{i+1}(\eta)\| + \|Q_{n+2+i}(\mu)) - P_{i}(\eta)\|)\nonumber\\
& \leq &  \Con\,  \|\nabla \varphi - \nabla (\varphi_{\eta,-n-2})\|_{\Gamma, (1)}\, \alpha^{n+2+i} + \Con\, \alpha^{n+2+i} \nonumber \\
& \leq &  \Con\,  \|\nabla \varphi\|_{\Gamma, (1)} \, \alpha^{n+2+i}  \;.
\end{eqnarray}

In a similar way (3.26) implies
\be
| \ta_{n+2+i}(x_0;\mu) - g^-_i (x_0;\eta)|
\leq   \Con \,\|\nabla \varphi\|_{\Gamma, (1)}\, \alpha^{n+2+i }\quad, \quad -n-1 \leq i \leq -1\;.
\ee

To prove (3.18), notice that 
$(A_{\jj}(\varphi) h) (x_0) = \Lambda_{\varphi,\jj}(x_0)\, h (\tQ_0(x_0;\mu))\;.$
The definition of $\Lambda_{\varphi,\jj}$ and $\tQ_{n+2}(x_0;\mu) = x_0$ give
\begin{eqnarray}
\log \Lambda_{\varphi,\jj}(x_0)  = \log \Lambda_{\varphi,\jj}(\tQ_{n+2}(x_0;\mu)) 
=  \sum_{i=0}^{n+1} \ta_i (x_0;\mu)\;.
\end{eqnarray}

Next, assume for simplicity that $n$ is odd (the other case is similar), and set $m = (n+1)/2$. 
Using (3.27) -- (3.31), we  get
\begin{eqnarray}
\log \Lambda_{\varphi,\jj}(x_0) - d (x_0;\mu)
& = &  \sum_{i=0}^{n+1} \ta_i(x_0;\mu)  + \sum_{i = -1}^{-n-1} [g(\sigma^{i}\eta) - g^-_i(x_0;\eta)] - \sum_{i=0}^{n+1} g^+_i(\mu)\nonumber\\
& = & \sum_{i = -m-1}^{-n-1} [g(\sigma^{i}\eta) - g^-_i(x_0;\eta)] + \sum_{i = 0}^m \left[\ta_i (x_0;\mu) -  g^+_i (\mu) \right] \nonumber\\
&    & +  \sum_{i = m+1}^{n+1} \left[\ta_i (x_0;\mu) - g^-_{i - n-2}(x_0;\eta)\right] + 
\sum_{i = -1}^{-m} [g(\sigma^{i}\eta) - g^+_{n+2+i}(\mu)]\nonumber\\
& = & O(\alpha^m)\, \|\nabla \varphi\|_{\Gamma, (1)} = O(\theta^n)\, \|\nabla \varphi\|_{\Gamma, (1)}\;.
\end{eqnarray}
Since by (3.23),
\be
|h (\tQ_0(x_0;\mu)) - h (Q_0(\mu))| = \|h\|_{\Gamma,1}\; O (\alpha^{n})\, ,
\ee
the above gives
\begin{eqnarray*}
\left| e^{d (x_0;\mu)}\, h (Q_0(\mu)) -  (A_{\jj}(\varphi) h) (x_0) \right| 
& \leq & \left| e^{d (x_0;\mu)} - e^{\log \Lambda_{\varphi,\jj}(x_0)} \right|\, \| h (Q_0(\mu))\| \\
&        & + \Lambda_{\varphi,\jj}(x_0)\, \left\| h (Q_0(\mu)) - h(\tQ_0(x_0;\mu)  \right\| \\
& \leq & e^{\max\{ d (x_0;\mu), \log \Lambda_{\varphi,\jj}(x_0)\}}\, \\
&        & \quad |d (x_0;\mu) -  \log \Lambda_{\varphi,\jj}(x_0)| \, \|h\|_{\Gamma,0} + \|h\|_{\Gamma,(1)}\, O(\alpha^n)\\
& \leq & \Con \,  \left( \|\nabla \varphi \|_{\Gamma,(1)}\, \|h\|_{\Gamma,0} +  \|h\|_{\Gamma,(1)}\right)\, \theta^{n}\;,
\end{eqnarray*}
which proves (3.18).

Similarly to (3.27) one gets $|f(\sigma^i (\eta)) - f^-_i (x_0;\eta)|  \leq  \Con\, \alpha^{|i|}\;,$ and also
$$|f^+_i(\mu) - f^+_i (x_0;\jj)| =  \left| \| Q_i (\mu) - Q_{i+1}(\mu)\| - \| Q_i (x_0;\jj) - Q_{i+1}(x_0;\jj)\|\right|
\leq \Con\, \alpha^{n+2-i}\;.$$
Combining these two estimates yields (3.17).

Next, using the notation from the beginning of this proof, notice that for any
$\mu$ as in (3.12) we have $Q_i(x_0;\jj) = \tQ_i(x_0;\mu)$ for all $i = 0,1 \ldots, n+2$ and
therefore $f^+_i (x_0;\jj) =  \| \tQ_i (x_0;\mu) - \tQ_{i+1}(x_0;\mu)\|$
for all $i = 0,1, \ldots, n+1$. (This has been used already in the proof of (3.17).) 

Define the function
$$\tWn (x_0;\cdot, \cdot) = \tWn_{1,\ell} (x_0;\cdot, \cdot):  \saa\times \C \longrightarrow \C$$
by $\tWn(x_0;\mu,s) = 0$ when $\mu_0 \neq 1$ or $\mu_{n+2} \neq \ell$ and
\begin{eqnarray}
\tWn(x_0;\mu,s)
& =  & e^{ s\tf_{n+2}(\mu) - \tg_{n+2}(\mu) -s\, \varphi(\tQ_0(x_0;\mu)) - s\, \sum_{i =0}^{n+1}  \| \tQ_i (x_0;\mu) - \tQ_{i+1}(x_0;\mu)\|}\nonumber \\
&    & \quad \times \Lambda_{\varphi,\jj}(x_0)\, h (\tQ_0(x_0;\mu))\;,
\end{eqnarray}
whenever $\mu_0 = 1$ and $\mu_{n+2} = \ell$, where $\jj = \jj^{(n+2)}(\mu)$ is defined by (3.12).

Using  (3.8), we can now write
\begin{eqnarray*}
&    &\sum_{|\jj| = n+3, j_0=1, j_{n+2} = \ell} u_{\jj}(x_0,-\i\, s) \\
& = & (-1)^{n}\, \sum_{\sigma^{n+2}\mu = \xi, \mu_0=1} e^{ -s\, \varphi(\tQ_0(x_0;\mu)) -s \sum_{i = 0}^{n+1}  
\| \tQ_i(x_0;\mu) - \tQ_{i+1}(x_0;\mu)\|}\, \Lambda_{\varphi,\jj}(x_0)\, h (\tQ_0(x_0;\mu))\\
& = & (-1)^{n}\, \sum_{\sigma^{n+2}\mu = \xi}  e^{ -s\tf_{n+2}(\mu) + \tg_{n+2}(\mu)}\,\tWn(x_0;\mu,s)
=  (-1)^{n}\, \left[ L_{-s\tf + \tg}^{n+2} \left(\tWn(x_0; \cdot, s)\right)\right] (\xi)\;.
\end{eqnarray*}
This and (3.11) imply
\begin{eqnarray}
&    & \Big| \left(L^n_s \mm_{n,s}(x_0) \gg_{s} \tv_s\right)(\xi) - \sum_{{|\jj | = n+3}\atop{j_{n+2} = \ell}}  u_{\jj}(x_0, s) \Big|\nonumber\\
& = &  \left| L^{n+2}_{-s\tf+\tg} \left[  \left(\Wn(x_0;\cdot,s) - \tWn(x_0;\cdot,s)\right)\right] (\xi) \right|.
\end{eqnarray} 

Standard estimates for Ruelle transfer operators  yield that there exists a global constant  $\Con > 0$  such that
\be
\left\| L^p_{-s\tf + \tg} H \right\|_\infty \leq \Con\,
e^{\Conf |\Ref(s)|}\, e^{p\, \Prf(-\Ref(s)\,\tf + \tg)}\, \|H\|_\infty \quad, \quad p \geq 0\:,\:\:
s\in \C\:,
\ee
for any continuous function $H : \saa \longrightarrow \C$. 

\bs

\noindent
{\bf Remark 4.} The above estimate can be derived e.g. from \cite{kn:St3} -- see the proof of Theorem 2.2, Case 1, there
which uses arguments from \cite{kn:Bo2} (see also the proof of Theorem 2.2 in \cite{kn:PP}).
More precisely, since $f,g\in \ff_\alpha(\sa)$, where $\alpha > 0$ is as in Proposition 1, we have $\tf,\tg \in \ff_\theta(\saa)$,
where $\theta = \sqrt{\alpha} \in (0,1)$. Setting $ u = -\Re(s)\, \tf + \tg$, $v = -\Im(s)\, \tf$, $\lambda = e^{\Prf(-\Ref(s)\, \tf + \tg)}$,
we have $-s\, \tf + \tg = u + \i\, v$, and $\lambda > 0$ is the maximal eigenvalue of the operator $L_u$ on $\ff_\theta(\saa)$.
Let $h \in \ff_\theta(\saa)$ be a positive corresponding eigenfunction, i.e. $L_u h = \lambda\, h$. It is then
easy to check (see e.g. (2.2) in \cite{kn:St3}) that 
$\|L^p_{-s\tf+\tg}H\|_\infty \leq \frac{\|h\|_\infty}{\min h}\, \lambda^p\, \|H\|_\infty$ for any $p \geq 0$ and any
continuous functions $H$ on $\saa$. To estimate $\frac{\|h\|_\infty}{\min h}$ one can use e.g. (3.6) in \cite{kn:St3} --
it follows from there that
$$\frac{\|h\|_\infty}{\min h} \leq K = e^{2\theta\, b/(1-\theta)} \lambda^M\, e^{M\, \|u\|_\infty}\;,$$
where $M \geq 1$ is a constant (one can take $M = 2$ in the situation considered here) and $b = \max\{ 1, \|u\|_\theta\}$.
Clearly, $\|u\|_\theta \leq |\Re(s)|\, \|\tf\|_\theta + \|\tg\|_\theta \leq \Con\, (|\Re(s)| +1)$ and similarly,
$\|u\|_\infty \leq \Con\, (|\Re(s)| +1)$, so (3.37) follows.

\bs

To use (3.37), we need to estimate
$$\sup_{\xi\in \saa} \left|  \left(\Wn(x_0;\cdot,s) - \tWn(x_0;\cdot,s)\right) (\xi)\right|\;.$$

Fix for a moment  $s\in \C$.  
According to the definitions
of $\Wn$ and $\tWn$, it is enough to consider $\mu \in \saa$ with $\mu_0 = 1$ and $\mu_{n+2} = \ell$. 
For such $\mu$, using (3.13), (3.16), (3.32), (3.33) and (3.35),  we have
\begin{eqnarray}
&        & |\Wn(x_0;\mu,s) - \tWn(x_0;\mu,s)|\nonumber\\
& =     & \left| e^{s\tf_{n+2}(\mu) - \tg_{n+2}(\mu) -s\, \varphi(\tQ_0(x_0;\mu)) 
-s\, \sum_{i=0}^{n+1}f^+_i(x_0;\jj) + \sum_{i=0}^{n+1} \ta_i (x_0;\mu)}\right|\nonumber \\
&        &  \times | e^{(s+\|\nabla \varphi\|_{\Gamma,(1)})O(\theta^n)  - s [c (x_0;\mu) - 
\sum_{i=0}^{n+1}f^+_i(x_0;\jj)]-s [\varphi(Q_0(\mu)) -\varphi(\tQ_0(x_0;\mu)]}\, 
h (Q_0(\mu))\nonumber \\
&        &   -  h (\tQ_0(x_0;\mu)) |\;.
\end{eqnarray}

To estimate (3.38), first notice that by (3.15) and Proposition 1,
$$|f(\sigma^i \mu) - f(\sigma^{i -(n+2)}\eta)| \leq \Con \, \alpha^i \quad, \quad 0 \leq i \leq n+2\;.$$
Using this, (3.24), (3.26) and Proposition 1 again, one gets
\begin{eqnarray}
\left| \tf_{n+2}(\mu) -  \sum_{i = 0}^{n+1}f^+_i (x_0;\jj)\right| 
& \leq & \Con + \left| f_{n+2}(\mu) -  \sum_{i=0}^{n+1}f^+_i(x_0;\jj) \right| \nonumber\\
& \leq & \Con + \sum_{i=0}^{n+1} \left| f(\sigma^i \mu) - f^+_i (x_0;\jj) \right| 
\leq  \Con\;,
\end{eqnarray}
for some global constant $\Con > 0$. Similarly,  it follows from (3.15), (3.29) and (3.30)  that
\be
\left|    \tg_{n+2}(\mu)  -  \sum_{i =0}^{n+1} \ta_i (x_0;\mu)\right| \leq \Con\, \|\varphi\|_{\Gamma,(1)}\;.
\ee

Next, notice that 
$$|e^{(s +\|\nabla \varphi\|_{\Gamma,(1)})O(\theta^n)}  - 1| \leq \Con \, e^{\Conf\, ( |\Ref(s)| +\|\nabla \varphi\|_{\Gamma,(1)})}\, 
(|s| + \|\nabla \varphi\|_{\Gamma,(1)})\; \theta^n \;.$$
Using the latter, (3.17), (3.18), (3.39) and (3.40) in (3.38) yields
\begin{eqnarray*}
&        & |\Wn(x_0;\mu,s) - \tWn(x_0;\mu,s)|\nonumber\\
& \leq & \Con\, e^{\Conf [ |\Ref(s)|\, (1+\|\varphi\|_{\Gamma,0}) + \|\nabla \varphi\|_{\Gamma,(1)}]}
    \; \left| e^{(s+\|\nabla \varphi\|_{\Gamma,(1)})O(\theta^n)} \,
      h (Q_0(\mu)) -  h (\tQ_0(x_0;\mu)) \right|\\
& \leq & \Con\, e^{\Conf [ |\Ref(s)|\, (1+\|\varphi\|_{\Gamma,0}) + \|\nabla \varphi\|_{\Gamma,(1)}]}
\left| e^{(s +\|\nabla \varphi\|_{\Gamma,(1)})O(\theta^n)} -1\right|\, | h (Q_0(\mu))| \\
&        & + \Con\, e^{\Conf [ |\Ref(s)|\, (1+\|\varphi\|_{\Gamma,0}) + \|\nabla \varphi\|_{\Gamma,(1)}]}         
\left| h (Q_0(\mu)) - h (\tQ_0(x_0;\mu)) \right|\\
& \leq & \Con\, e^{\Conf [ |\Ref(s)|\, (1+\|\varphi\|_{\Gamma,0}) + \|\nabla \varphi\|_{\Gamma,(1)}]}\,
\left[ (|s|+\|\nabla \varphi\|_{\Gamma,(1)})\,  \|h\|_{\Gamma,0} +  \|h\|_{\Gamma,(1)}\right]\, \theta^n\;.
\end{eqnarray*}
Thus, choosing the global constant $\Con > 0$ sufficiently large, combining the above with (3.37) gives
\begin{eqnarray}
&        & \left| L^{n+2}_{-s\tf+\tg} \left[\left( \Wn(x_0;\cdot,s) 
- \tWn(x_0;\cdot,s)\right)\right] (\xi)\right| \nonumber \\
& \leq & \Con\, e^{\Con [ |\Ref(s)|\, (1+\|\varphi\|_{\Gamma,0}) + \|\nabla \varphi\|_{\Gamma,(1)}]}\,
\left[ (|s|+\|\nabla \varphi\|_{\Gamma,(1)})\,  \|h\|_{\Gamma,0} +  \|h\|_{\Gamma,(1)}\right] \nonumber\\
&        & \times \Bigl(e^{\Prf(-\Ref (s)\,\tf + \tg)}\, \theta\Bigr)^{n+2}\;.
\end{eqnarray}
Next we have (see for example Ch. 4 in \cite{kn:PP}) 
$$\frac{d}{ds}{\rm Pr}(-s\tf + \tg)\Big\vert_{s = s_0} = - \int_{\Sigma_A^{+}} \tf d\nu = -\int_{\Sigma_A^{+}} f d\nu = -c_0 < 0,$$
where $\nu$ is the equilibrium state of $(-s_0 \tf + \tg)$. Recall that $\Pr(-s_0 \tf + \tg) = 0$, so
$e^{\Pr( - \Re(s)\, \tf + \tg)} < 1$ for $\Re(s) > s_0$. 
Now assume $s_0-a \leq \Re (s)$ with some small constant $a > 0$. Then
$$e^{\Prf(- \Ref (s) \tf + \tg)} = 1 + c_0(s_0 - \Re (s) ) + O((\Re (s) - s_0)^2) \leq 1 + c_1 a$$
for some constant $c_1 > 0$. Thus,
$$e^{\Prf(- \Ref (s) \tf + \tg)}\theta \leq \theta + c\, a$$
for some global constant $c = c_1\theta > 0$. Combining this with (3.41), completes the proof of (3.6).

\section{Estimates for the derivatives}
\renewcommand{\theequation}{\arabic{section}.\arabic{equation}}
\setcounter{equation}{0}

In this section we prove Theorem 3(b).  Throughout we assume that $p \geq 1$.

For any $x\in \Gamma_\ell$ close to $x_0$ and any $\eta\in \sa$
with $\eta_0 = \ell$ define the points $\tP_j(x;\eta)$ and the functions $f_i^-(x;\eta)$, $g_i^-(x;\eta)$,
$\phi^-(x;\eta,s)$, etc.,  as in the beginning of Section 3 replacing
the point $x_0$ by $x$. 
We will assume that the segment $[\tP_{-1}(x_0;\eta),x_0]$ has no
common points with the interior of $K_\ell$ and $x$ is close enough to
$x_0$ so that the same holds with $x_0$ replaced by $x$.

By Proposition 4 in Appendix A below there exists a unique phase function $\psi_\eta$ (also depending on $x_0$) defined in a 
neighborhood $U$ of $x_0$ in $\Gamma_\ell$, such that $\psi_\eta(x_0)  = 0$ and the backward trajectory 
$\gamma_-(x,\nabla \psi_\eta (x))$
of any point $x\in U$ with $\psi_\eta(x) = 0$ has an itinerary $(\ldots, \eta_{-\ell}, \ldots, \eta_{-1}, \eta_0)$, that is
$$\nabla \psi_\eta (x) =  \frac{\tP_{0}(x;\eta) - \tP_{-1}(x;\eta)}{\|\tP_{0}(x;\eta) - \tP_{-1}(x;\eta)\|}$$
for any $x\in \cc_{\psi_\eta}\cap U$. (Notice that in general $\psi_\eta$ is different from the functions $\varphi_{\eta,j}$
defined in the beginning of Sect. 3.) 
For any $i < 0$, denoting $J = (\eta_i, \eta_{i+1}, \ldots, \eta_{-1},\eta_0)$, 
we can write $\psi_\eta = (\psi_{\eta,i})_J$
for some phase function $\psi_{\eta,i}$ (defined on some naturally defined open subset $V_{\eta,i}$of $\R^N$)
satisfying Ikawa's condition (${\mathcal P}$) on $\Gamma_{\eta_i}$. We then have
$\tP_i (x;\eta) = X^{-i}(x, \nabla (\psi_{\eta,i})_J) \;.$
As in the beginning of Sect. 3 (see (3.2) there) one derives that there exists
a global constant $\Con_p > 0$ such that
$\|\psi_{\eta,i}\|_{(p)}(V_{\eta,i} \cap B_0) \leq \Con_p\;$
for all $\eta$ and $i < 0$.
Using (2.4) in Proposition 2 with $\varphi = \psi_{\eta,m}$ for some $m \geq i$ and replacing $C_p$ with a larger
global constant if necessary, we get
\be
\|\tP_i (\cdot; \eta)\|_{\Gamma,p} (x) \leq \Con_p\, \alpha^{|i|} \quad, \quad i < 0\;.
\ee
Similarly, for any $\mu \in \saa$ with $\mu_0 = 0$ and
$\mu_{n+2} = k$ we have
\be
\|\tQ_i(\cdot;\eta)\|_{\Gamma,p} (x) \leq \Con_p\, \alpha^{n+2-i} \quad, \quad 0\leq i \leq n+2\;,
\ee
and 
\be
\|\tQ_i (\cdot;\mu) - \tP_{i -n-2}(\cdot;\eta)\|_{\Gamma, p} (x) \leq \Con_p\, \alpha^{i} \quad, \quad 0\leq i \leq n+2\;.
\ee

Next, recall the function $\Lambda_\varphi$ from the beginning of this section. By Proposition 2,
\be
\|\nabla \varphi_J\|_{\Gamma, p}  \leq \Con_p\, \|\nabla \varphi\|_{\Gamma,(p)}
\ee
for any finite admissible configuration $J$.

Since for any $i < 0$ we have $g_i^-(x;\eta) = \log \Lambda_{\psi_{\eta,i}}(\tP_{i+1}(x;\eta))\;,$
it follows from (4.1)--(4.3) and Proposition 3 that for any $p \geq 1$ there exists a global constant $\Con_p > 0$  such that 
\be
\|g_i^-(\cdot; \eta)\|_{\Gamma,p} (x) \leq \Con_p \, \alpha^{|i|}\quad, \quad  i < 0\;.
\ee
Similarly, according to (3.28) and Proposition 2, 
\be
\| \ta_i(\cdot;\mu)\|_p (x) \leq \Con_p \,\|\nabla \varphi\|_{\Gamma,(p)}\, \alpha^{n+2-i} \quad, \quad 0\leq i \leq n+2\;,
\ee
and as in the proof of (3.31) one gets,
\be
\| \ta_i (\cdot; \mu) - g^-_{i-n-2} (\cdot; \eta) \|_p (x) \leq \Con_p \, 
\|\nabla \varphi\|_{\Gamma,(p+1)}\, \alpha^{i} \quad, \quad 0\leq i \leq n+2\;.
\ee

Next, given $x$ as above, $\mu$ and $n$ with $\mu_{n+2} = \ell$, define 
$\Wn(x;\mu,s)$ by (3.10), $\eta$ by (3.15) and $\tWn(x;\mu,s)$ by (3.35) replacing $x_0$ by $x$. 
We will estimate the derivatives of
$$\Wn(x;\mu,s) - \tWn(x;\mu,s)$$
with respect to $x$. 

First look at the first derivatives $D_v[\Wn(\cdot;\mu,s) - \tWn(\cdot;\mu,s)] (x)$, where $v \in S_x \Gamma$.
Writing
$\phi^-(x;\eta, s) = -s \, \phi^-_1(x;\eta) + \phi^-_2(x;\eta)\;,$ where
$$\phi_1^-(x ;\eta) =  \sum_{i = -1}^{-\infty} [f(\sigma^i(\eta)) - f^-_i(x;\eta)] \quad , \quad
\phi_2^-(x,\eta) =  \sum_{i = -1}^{-\infty} [g(\sigma^i (\eta)) - g^-_i (x;\eta)]\;,$$
notice that for any $x, x'\in \Gamma_\ell$ (close to $x_0$) we have
$$\phi_1^-(x;\eta) - \phi_1^-(x';\eta) = -\psi_\eta(x) + \psi_\eta(x')\;,$$
so $D_v (\phi_1^-(\cdot;\eta))(x) = D_v (\psi_\eta(x))$. Therefore by (3.14)
\be
D_v z (\cdot;\mu,s)(x) = - s\, D_v\psi_\eta (x) + \sum_{i = -1}^{-\infty} D_v (g^-_i (\cdot;\eta))(x) \;.
\ee

Next,  using the notation
$\jj = (\mu_0, \mu_1, \mu_2, \ldots, \mu_{n+2})$
and
\begin{eqnarray*}
\tz(x;\mu,s) 
=   s\tf_{n+2}(\mu) - \tg_{n+2}(\mu) -s\, (\varphi_{\mu_0})_{\jj}(x) \;,
\end{eqnarray*}
it follows from (3.38) that
\begin{eqnarray}
&     & \Wn(\cdot;\mu,s) - \tWn(\cdot;\mu,s)(x) \nonumber\\
&  = & e^{z (x;\mu,s) -s\, \varphi(Q_0(\mu))}\,  h (Q_0(\mu)) - e^{\tz (x;\mu,s)}\,
\Lambda_{\varphi,\jj}(\tQ_{n+2}(x;\mu))\, h (\tQ_0(x;\mu)) \nonumber\\
& = & (I)(x) + (II)(x)\;,
\end{eqnarray}
where
$$(I)(x) = [ e^{z (x;\mu,s) - s\, \varphi(Q_0(\mu))}  - e^{\tz (x;\mu,s) + \log \Lambda_{\varphi,\jj}(\tQ_{n+2}(x;\mu))} ]\,  h (Q_0(\mu)) \;,$$
and
$$(II)(x) = e^{\tz (x;\mu,s)}\, \Lambda_{\varphi,\jj}(\tQ_{n+2}(x;\mu))\, [h(Q_0(\mu)) - h (\tQ_0(x;\mu))]\;.$$

Let $\oo$ be a small compact connected neighborhood of $x$ in $\Gamma$. 
Fix temporarily $\mu$, $s$, $n$ and $\eta$ with (3.15), and set
$$A(y) = z (y;\mu,s) - s\, \varphi(Q_0(\mu)) \quad, \quad
B(y) = \tz (x;\mu,s)+ \log \Lambda_{\varphi,\jj}(\tQ_{n+2}(x;\mu))\quad , \quad y\in \oo\;.$$

To estimate $(I)$ first notice that by the estimates in Sect. 3,
$$\| A \|_0 (\oo) = O (|s| + |s|\, \|\varphi\|_{\Gamma,0} +  \|\nabla \varphi\|_{\Gamma,(1)})\;,$$
and
\be
|e^A|_{\Gamma,0}(\oo) \leq \Con\, e^{\Conf\,[ |\Ref(s)|\, (1+ \|\varphi\|_{\Gamma,0}) + \| \nabla \varphi\|_{\Gamma,(1)}] }\;.
\ee

It follows from (4.6) and (3.40) that $|\tg_{n+2}(\mu)| \leq \Con\,   \|\nabla \varphi\|_{\Gamma,(1)}$.
Combining this with the definition of $\tz (x;\mu,s)$ and (3.39) implies
$$\| \tz (\cdot;\mu,s) \|_0 (\oo) = O (|s| + |s|\, \|\varphi\|_{\Gamma,0} +  \|\nabla \varphi\|_{\Gamma,(1)})\quad, \quad
\| B \|_0 (\oo) = O (|s| + |s|\, \|\varphi\|_{\Gamma,0} +  \|\nabla \varphi\|_{\Gamma,(1)}) \;.$$

Next, we will estimate the derivatives of $A$ and $B$.
For any $q\geq 1$ and any $y\in \oo$, using (4.8), (2.1) and (4.5), we get
\begin{eqnarray}
\|A\|_{\Gamma,q}(y)
& =     & \| s \, \phi_1^-(\cdot;  \eta) - \phi_2^-(\cdot;\eta)\|_{\Gamma,q}(y)
\leq  |s|\, \|\nabla \psi_\eta\|_{\Gamma,q}(y) + \sum_{i=-1}^{-\infty} \|g_i^-(\cdot;\eta)\|_{\Gamma,q}(y)\nonumber\\
& \leq & |s|\, \Con_q + \Con_q\, \sum_{i =-1}^{-\infty} \alpha^{|i|} \leq \Con_q\, (|s|+1)\;.
\end{eqnarray}
Thus, for any $q\geq 0$,
\begin{eqnarray*}
\|e^A\|_{\Gamma,q}(\oo)
& \leq & \Con_q \| e^A \|_{\Gamma,0}(\oo) \, (\max_{1\leq i \leq q} \|A \|_{\Gamma,i}(\oo))^q \nonumber\\
& \leq & \Con_q\, e^{ \Conf \,[ |\Ref(s)| \, (1+  |\varphi|_{\Gamma,0}) + \| \nabla \varphi\|_{\Gamma,(1)}]}\, (|s|+1)^q\;.
\end{eqnarray*}

Similarly, (4.4) gives
$$\|\tz (\cdot;\mu,s)\|_{\Gamma,q}(y) = \|s\, (\varphi_{\mu_0})_{\jj}\|_{\Gamma,q}(y) \leq \Con_q\, |s|\, \|\nabla \varphi\|_{\Gamma,(q)}\;,$$
while (3.31) and (4.6) imply
$$\|\log \Lambda_{\varphi,\jj}(\cdot)\|_{\Gamma,q}(y)
\leq   \sum_{i=0}^{n+1} \| \ta_i (\cdot;\mu)\|_{\Gamma,q}(y) \leq \Con_q\, \|\nabla \varphi\|_{\Gamma,(q)}\;$$
for any $q\geq 0$, so
$$\|B\|_{\Gamma,q}(y) \leq \Con_q\, (|s|+1)\, \|\nabla \varphi\|_{\Gamma,(q)}  \quad, \quad y\in \oo\;.$$

The next step is to estimate the derivatives of $A - B$. First notice that by Proposition 2 and (2.1) we have
$$\| \nabla \psi_\eta  -  \nabla (\varphi_{\mu_0})_J \|_{\Gamma,q} (\oo) \leq 
\Con_q \, \alpha^{n}\, \|\nabla \psi_\eta - \nabla \varphi_{\mu_0}\|_{\Gamma,(q)}
\leq \Con_q \, \alpha^{n}\, \|\nabla \varphi \|_{\Gamma,(q)}\;.$$

Set again $m = \frac{n+1}{2}$, assuming for simplicity that $n$ is odd, and $\theta = \sqrt{\alpha} \in (0,1)$.
As in the proof of (3.18) above, for any $y\in \oo$ and any $q \geq 1$, using (4.5), (4.6) and (4.7), we have
\begin{eqnarray*}
\|A - B\|_{\Gamma,q}(y)
& \leq &  
\left\| - s\, \psi_\eta + \sum_{i=-1}^{-\infty} g_i^-(\cdot;\eta) + s\, (\varphi_{\mu_0})_J - \sum_{i=0}^{n+1} \ta_i(\cdot;\mu) \right\|_{\Gamma,q}(y)\\
& \leq & |s|\, \| \psi_\eta  -  (\varphi_{\mu_0})_J \|_{\Gamma,q}(y) + 
\sum_{i = -m -1}^{-\infty} \|g_i^-(\cdot;\eta)\|_{\Gamma,q}(y)\\
&        & + \sum_{i= 0}^{m} \|\ta_i (\cdot;\mu)\|_{\Gamma,q}(y) + 
\sum_{i = m+1}^{n+1} \|\ta_i (\cdot;\mu) - g_{i -n-2}^-(\cdot;\eta)\|_{\Gamma,q}(y)\\
& \leq & \Con_q (|s|\, \|\nabla \varphi \|_{\Gamma,(q)} + \|\nabla \varphi \|_{\Gamma,(q+1)})\, \theta^n\;. 
\end{eqnarray*}
From Sect. 3, a similar estimate holds for $q = 0$. Consequently,
\begin{eqnarray*}
\|e^{B-A}\|_{\Gamma,q}(\oo)
& \leq & \Con_q \| e^{B - A} \|_{0}(\oo) \, (\max_{1\leq i \leq q} \|B - A\|_{\Gamma,i}(\oo))^q \nonumber\\
& \leq & \Con_q\, e^{ \Conf \, (|\Ref(s)| + \|\nabla \varphi \|_{\Gamma,(1)})}\, 
(|s|\, \|\nabla \varphi \|_{\Gamma,(q)} + \|\nabla \varphi \|_{\Gamma,(q+1)})^q \theta^{nq}\;.
\end{eqnarray*}

Finally, as in the estimate just after (3.40), it follows that
$$\|e^{B - A} - 1\|_{0} (\oo) \leq  \Con\, e^{\Conf \, (|\Ref(s)| + \|\nabla \varphi \|_{\Gamma,(1)})}\,(|s| + \|\nabla \varphi\|_{\Gamma,(1)})\,  \theta^n\;.$$

The above, (4.10) and (4.11) imply that for any $q\geq 1$,
\begin{eqnarray*}
\|(I)\|_{\Gamma,q}(\oo)
\leq  \Con_q\,  \|h\|_{0}(\Gamma)\,\, e^{ \Conf \, [ |\Ref(s)|\, (1+ \|\varphi\|_{\Gamma,0}) + \|\nabla \varphi \|_{\Gamma,(1)} ]}\, 
(|s|\, \|\nabla \varphi \|_{\Gamma,(q)} + \|\nabla \varphi \|_{\Gamma,(q+1)})^q\, \theta^n\; .
\end{eqnarray*}

Using similar estimates, for any $q \geq 1$ one gets
\begin{eqnarray*}
|(II)|_{\Gamma,q}(\oo)
& \leq & \Con_q\, \alpha^n\,  e^{ \Conf \, [ |\Ref(s)|\, (1+
\|\varphi\|_{\Gamma,0} ) + \| \nabla \varphi\|_{\Gamma,(1)}]}\,
   \sum_{r=0}^{q-1} 
(|s|+1)^{r+1} \, (\|\nabla \varphi\|_{\Gamma,(r)})^{r+1}\, \|h\|_{\Gamma,q-r}(\oo)\,\;.
\end{eqnarray*}

It now follows from (4.9) and the estimates for $(I)$ and $(II)$ found above that for any $p \geq 1$ we have
\begin{eqnarray*}
&        & \| \Wn(\cdot;\mu,s) - \tWn(\cdot;\mu,s)\|_{\Gamma,(p)}(\oo)\\
& \leq & \Con_p\,\theta^n\,  e^{ \Conf \, [ |\Ref(s)|\, (1+ \|\varphi\|_{\Gamma,0}) + \|\nabla \varphi \|_{\Gamma,(1)} ]}\, 
\times  \sum_{r=0}^{q}  \left[ |s|\, \|\nabla \varphi \|_{\Gamma,(r)} + 
\|\nabla \varphi  \|_{\Gamma,(r+1)}\right]^{r+1}\, \|h\|_{\Gamma,q-r}(\oo)\;.
\end{eqnarray*}
Combining this with (3.6), (3.36) and the argument from the end of Sect. 3 completes the proof of Theorem 3.
\endofproof


\section{Estimates for $w_{0, j}(x, s)$}
\renewcommand{\theequation}{\arabic{section}.\arabic{equation}}
\setcounter{equation}{0}

\def\mn{{\mathcal M}_{n,s}}
\def\mc{{\mathcal M}}
\def\rc{{\mathcal R}_s}
\def\ii{{\bf \i}}
\def\ta{\tilde{a}}
\def\Ep{E_p(s, \varphi, h)}
\def\sj{\sum_{|i| = n+2, i_{n+2} = j}}

Our purpose in this section is to prove that the series 
$$w_{0,j}(x, s) = \sum_{n= n_j}^{\infty} \sum_{|\jj | = n+3,\: j_{n+2} = j} u_{\jj}(x,s),\: x \in \Gamma_j$$
is convergent and that $\wj(x, s)$ is an analytic function for $s \in \d1$ with values in $C^{\infty}(\Gamma_j).$ Since we deal with initial data 
$m(x,s) = u_1(x,s)$ on $\Gamma_1$ we set $n_1 = -2$ and $n_j = -1,\: j=2,...,\kappa_0$. 
By Theorem 3, it is clear that the problem is reduced to the convergence of the series
$$\sum_{n=0}^{\infty} (L_s^n {\mathcal M}_{n,s}(x){\mathcal G}_s\tilde{v}_s)(\xi),\: x \in \Gamma_j.$$

Throughout this and the following sections we will use the notation
$$\Ep = \begin{cases} e^{C_p [ |\Ref(s)|\, (1+\|\varphi\|_{\Gamma,0}) + \|\nabla \varphi\|_{\Gamma,(1)}]}
\sum_{j = 0}^p \Bigl(|s| \|\nabla \varphi\|_{\Gamma, j} + \|\nabla \varphi\|_{\Gamma,j+1} \Bigr)^{j + 1}  \|h\|_{\Gamma, p - j}\: {\text if} \: p \geq 1,\\
C_0e^{C_p [ |\Ref(s)|\, (1+\|\varphi\|_{\Gamma,0}) + \|\nabla \varphi\|_{\Gamma,(1)}]}\left[\Bigl(|s|+ \|\nabla \varphi\|_{\Gamma, (1)}\Bigr)\,  \|h\|_{\Gamma,0} + 
\|h\|_{\Gamma, (1)}\right]\:{\text if}\: p = 0,\end{cases}$$
where as before by $C_p$ we denote  positive global constants depending on $p$ which may change from line to line.

First we will establish for $\sigma_0 \leq \Re (s) \leq 1$ the inequality
\begin{equation} \label{eq:5.1}
\|L_s^n \mn(.) - L_s^{n-1} \mc_{n-1, s}(.) L_s\|_{\Gamma, p} \leq C_p \Ep \theta^{n},
\end{equation}
where $L_s = -L_{-s\tilde{f} + \tilde{g}}$  and $\sigma_0 < s_0$. The precise choice of $\sigma_0$ depends on the estimates (3.3) and will be discussed below. For this purpose we write
$$\Bigl(L_s^n \mn - L_s^{n-1} \mc_{n-1, s} L_s\Bigr)w(\xi) = -L^{n+1}_s\Bigl[Y^{(n)}(x; s, \mu) - \widetilde{Y}^{(n)}(x; s, \mu)\Bigr](\xi),$$
where
$$Y^{(n)}(x; s, \mu) = \exp\Bigl(-\phi^-(x; \sigma^{n+1} e(\mu), s) - \chi(\sigma^{n+1} e(\mu), s)\Bigr) w(\mu),$$
$$\widetilde{Y}^{(n)}(x; s, \mu) = \exp \Bigl(-\phi^-(x; \sigma^n e(\sigma \mu), s) - \chi(\sigma^n e (\sigma \mu), s)\Bigr) w(\mu).$$

The inequality (\ref{eq:5.1}) follows from the estimates
\begin{eqnarray} \label{eq:5.2}
& &\Bigl\|\phi^{-}(x; \sigma^{n+1}e(\xi), s) - \phi^{-}(x; \sigma^{n}e(\sigma(\xi)), s)\Bigr\|_{\Gamma, p} \leq C_p \Ep \theta^n, \\
& &|\chi(\sigma^{n+1} e(\xi), s) - \chi(\sigma^n e(\sigma(\xi)), s)| \leq C (1 + |s|) \theta^n \label{eq:5.3}
\end{eqnarray}
and the form of the operators $\mc_{n,s}(x)$. The estimate (\ref{eq:5.3}) is a consequence of the choice of $\chi_1, \chi_2$ and the fact that 
$f, g \in {\mathcal F}_{\theta}(\Sigma_A)$. To prove (\ref{eq:5.2}), notice that 
$$\Bigl|\sum_{i = -1}^{-\infty}[f(\sigma^{n+1+i}e(\xi)) - f(\sigma^{n+i}e(\sigma(\xi)))]\Bigr| \leq C\theta^n,$$
and similar estimates hold for the function $g$. The terms involving $f$ and $g$ are independent of $x$ and they are not important for the estimates 
of the derivatives. To deal with the terms depending on $x$, recall that
$$\phi^-(x; \eta) = - s \phi_1^-(x; \eta) + \phi_2^-(x; \eta)$$
with $D_v(\phi_1^-(.;\eta)(x) = D_v(\psi_{\eta}(x)).$ Here and below we use the notation of the previous section. On the other hand,
\begin{equation} \label{eq:5.4}
\|\nabla\psi_{\sigma^{n+1}e(\mu)}(x) - \nabla\psi_{\sigma^{n}e(\sigma(\mu))}(x)\|_{\Gamma, p} \leq C_p \alpha^n.
\end{equation}
In fact, the backward trajectories $\gamma_{-}(x, \nabla \psi_{\sigma^{n+1}e(\mu)}(x))$ and $\gamma_{-}(x, \nabla \psi_{\sigma^{n}e(\sigma(\mu))}(x))$ follow an itinerary
$(\mu_{n+1}, \mu_n,\ldots,\mu_1)$ and we can apply Proposition 2. Now we repeat the argument used in the previous section for the estimate of $\|A - B\|_{\Gamma, p}$. 
Set $m = \frac{n + 1}{2}$ and assume for simplicity  that $n$ is odd. For fixed $n$ we set $\eta = \sigma^{n+1} e(\mu),\: \tilde{\eta} =\sigma^n e(\sigma(\mu))$. The estimate of
$$\|\phi_1^{-}(x; \eta) - \phi_1^{-}(x; \tilde{\eta})\|_{\Gamma, p}$$
follows from (\ref{eq:5.4}). Next we write
$$\sum_{i = -1}^{-\infty} \Bigl(g_i^-(x;\eta) - g_i^-(x; \tilde{\eta}) \Bigr)= \sum_{i = -m-1}^{-\infty}\Bigl(g_i^-(x; \eta) - g_i^-(x; \tilde{\eta})\Bigr)$$
$$+ \sum_{i = m+1}^{n+1} (g_{i -n-2}^-(x;\eta) - \ta_i(x; \mu)) - \sum_{i = m+1}^{n+1} (g_{i-n-2}^-(x, \tilde{\eta})- \ta_i(x; \mu)).$$ 
The $\|.\|_{\Gamma, p}$ norms of the sums from $i = m +1$ to $n +1$ can be estimated as in Section 4 by using (4.7), since 
$$\eta = \sigma^{n+1}e(\mu) = (\ldots,*,*,\mu_0, \mu_1,\ldots, \mu_{n+1}= \ell ,\mu_{n+2},\ldots ),$$
$$\tilde{\eta} = \sigma^n e(\sigma(\mu)) = (\ldots,*,*,\mu_1,\ldots, \mu_{n+1} = \ell , \mu_{n+2},\ldots),$$
and 
$$\sum_{i = m+1}^{n+1} \|g_{i-n-2}^-(x;\eta) - \ta_i(x; \mu))\|_{\Gamma, p} \leq \sum_{i = m+1}^{n+1} \alpha^i\,,$$
$$\sum_{i = m+1}^{n+1} \|g_{i -n-2}^-(x;\tilde{\eta}) - \ta_i(x; \mu))\|_{\Gamma, p}\leq \sum_{i = m+1}^{n+1} \alpha^i\; .$$
To estimate the sums from $i = -m -1$ to $-\infty$, we apply (4.5) and this completes the proof of (\ref{eq:5.1}).\\

From  the representation
$$L_s^n \mn = \sum_{k = 1}^n \Bigl(L_s^k \mc_{k,s} - L_s^{k-1} \mc_{k-1, s} L_s\Bigr) L_s^{n-k} + \mc_{0,s}L_s^n$$
we get
$$\sum_{n=1}^{\infty} L_s^n \mn w = \sum_{n=1}^{\infty} \Bigl[\sum_{k=1}^n \Bigl(L_s^k \mc_{k,s} - L_s^{k-1} \mc_{k-1,s} L_s\Bigr) L_s^{n-k}w + \mc_{0,s} L_s^n w\Bigr].$$
Since $s_0 \in \R$ is the abscissa of absolute convergence, for $\Re
(s) > s_0$ we have $\Pr(-\Re (s) \tf + \tg) < 0$ and $\|L_s^n\|_{\infty} \leq 1,\:\forall n$. Consequently, the double sum in the right hand side is absolutely convergent 
for $\Re (s) > s_0$ and we can change the order of summation. Applying Fubini's theorem, we are going to examine
\begin{equation} \label{eq:5.5}
\sum_{n=0}^{\infty} L_s^n \mn {\mathcal G}_s \tilde{v}_s = \Bigl(\mc_{0,s} + {\mathcal Q}_s \Bigr)\sum_{n=0}^{\infty} L_s^n {\mathcal G}_s \tilde{v}_s,
\end{equation}
where
$${\mathcal Q}_s = \sum_{k=1}^{\infty}\Bigl(L_s^k \mc_{k,s} - L_s^{k-1} \mc_{k-1, s} L_s\Bigr).$$
According to (\ref{eq:5.1}), the series defining ${\mathcal Q}_s$ is absolutely convergent for $\sigma_0 \leq \Re (s) \leq 1$ and 
$$\|{\mathcal Q}_s\|_{\Gamma,p} \leq C_p \Ep.$$
Consequently, the problem of the analytic continuation of the left hand side of (\ref{eq:5.5}) for $\Re (s) < s_0$ is reduced to that of the series $\sum_{n=0}^{\infty} L_s^n w_s, w_s = {\mathcal G}_s \tilde{v}_s.$\\

The analysis of $\sum_{n=0}^{\infty} L_s^n w_s$ is based on Dolgopyat type estimates (3.3) and we must show that
$w_s = h_s \circ \Phi$ with some $h_s \in \clip_u(\mtb)$ (see Appendix C for the definition of the map $\Phi$ and the space $\clip_u(\mtb)$). 
This assertion is proved in Appendix C, where we show that for $|\Re (s)| \leq a$ we have $\|h_s\|_{\lip, t} \leq C_0$ with $C_0$ independent on $s$.
Thus for $s = \tau + \ii t,\: \sigma_0 \leq \tau \leq 1,\: |t| \geq t_0 > 1$, we get
$$\sum_{n = 0}^{\infty}\|\tilde{L}_s^n w_s\|_{\infty} \leq \sum_{p=0}^{\infty} \sum_{l=0}^{[\log|t|] -1}C\rho^{p[\log|t|]}e^{l\Pr(-\tau\tf + \tg)} \|h_s\|_{\lip, t }$$
$$\leq \frac{C C_0}{1 - \rho^{[\log|t|]}} \sum_{l=0}^{[\log|t|]- 1} e^{l\Pr(-\tau\tf + \tg)}$$
$$\leq C_1 \max\{\log|t|,\: |t|^{\Pr(-\tau\tf + \tg)}\}.$$
On the other hand, for $\sigma_0$ sufficiently close to $s_0$ we have
$\Pr(-\sigma_0 \tf + \tg) = \tilde{\beta}_0 < 1.$ Combining this with the estimate for ${\mathcal Q}_s$, we conclude that for 
$\sigma_0 \leq \Re (s) $ and $|t| \geq t_0$ we have
$$\Bigl\|\sum_{n=0}^{\infty} L_s^n \mn {\mathcal G_s}\tilde{v}_s\Bigr\|_{\Gamma,0} \leq C_2 |t|^{1 + \tilde{\beta}_0}.$$

The analysis in Sect. 5 of \cite{kn:I1} implies that the series defining $w_{0, j}(x, s)$ is absolutely convergent for 
$x \in \Gamma_j,\:\Re (s) \geq s_0 + d > s_0$ and we have
\begin{equation} \label{eq:5.6}
\|w_{0, j}(x, s)\|_{\Gamma_j,0} \leq C_{j, d},\: \Re (s) \geq s_0 + d.
\end{equation}
On the other hand, the analytic continuation of the series
$\sum_{n=0}^{\infty}L_s^n \mn {\mathcal G}_s\tilde{v}_s$ established above and Theorem 3(a) with a sufficiently small $\epsilon = s_0 - \Re (s) > 0$ guarantee an 
analytic continuation of $w_{0,j}(x, s)$ for $x \in \Gamma_j,\: \Re (s) \geq \sigma_0, \: |\Im (s)| \geq t_0$ with $\sigma_0 = s_0 - \epsilon.$ Applying Theorem 3(a) 
once more for $s = \sigma_0 + \ii t$, we get the estimate
$$\|w_{0, j}(x,\sigma_0 + \ii t)\|_{\Gamma_j,0} \leq D_{j} |t|^{1 + \tilde{\beta}_0}.$$
The same argument works for all $\ell = 1,...,\kappa_0$ and we get the same estimate for
$$w_{0,\ell}(x,  s) = \sum_{n=n_{\ell}}^{\infty} \sum_{|\jj | = n+3, j_0 = 1,\: j_{n+2} = \ell} u_{\jj}(x, s),\: x \in \Gamma_{\ell}\, .$$
Clearly, we can choose $0 < \tilde{\beta}_0 < 1$ independent of $\ell = 1,\ldots, \kappa_0.$\\

Now we will obtain $C^p(\Gamma_j)$ estimates for $w_{0, j}(x, s).$ To examine the regularity of the functions $w_{0,j}(x, s)$ on $\Gamma_j$, set 
$$U_{n+2,j}(x, s) = \sum_{|\jj| = n+3, j_{n+2} = j} u_{\jj}(x, s).$$
We start with an estimate of the $C^p(\Gamma_j)$ norms of $U_{n+2,j}(x, s)\bv_{\Gamma_j}$. For this purpose, applying Theorem 3(b) with $p \geq 1$, 
we must estimate the norms $\|L^s \mc_{n,s}(.)w_s\|_{\Gamma_j, p}$, where $w_s = {\mathcal G}_s \tilde{v}_s$ and $L_s^n$ are independent of $x \in \Gamma.$ We write
$$L_s^n \mc_{n,s} w_s= \mc_{0,s} L_s^n w_s+ \sum_{k = 1}^m \Bigl(L_s^k \mc_{k,s} - L_s^{k-1} \mc_{k-1, s} L_s \Bigr) L_s^{n-k}w_s$$
$$+ \sum_{k = m+1}^n \Bigl(L_s^k \mc_{k,s} - L_s^{k-1} \mc_{k-1, s} L_s \Bigr) L_s^{n-k}w_s = B_0 + B_1 + B_2,$$
where $m = [n/2]$. We apply the estimate (3.3) combined with $\|h_s\|_{\lip, t} \leq C_0,\: t = \Im (s)$ and we obtain
$$\|L_s^n w_s\|_0 \leq C \rho^n e^{\log |t| \Bigl[\Pr (-s\tf + \tg) - \log\rho \Bigr]} \leq C \rho^n |t|^{\beta_0},\: \forall n \in \N$$
with $0 < \rho < 1,\:\beta_0 = \Pr(-\sigma_0 \tf + \tg) - \log \rho > 0.$ Increasing $\rho$, we can arrange $\beta_0 < 1$ but this is not important for our argument (see also Remark 7 in Appendix C).

For the term $B_0$  we get
$$\|B_0\|_{\Gamma_j, p} \leq C_{p}|\Im (s)|^{\beta_0}E_p(s, \varphi, h) \rho^n.$$
In the same way for the term $B_1$ we have
$$\|B_1\|_{\Gamma_j, p} \leq C_p' |\Im (s)|^{\beta_0}\Ep \sum_{k = 1}^m \theta^k \rho^m \leq C''_p |\Im (s)|^{\beta_0}\Ep (\sqrt{\rho})^n.$$
Finally, for $B_2$ we obtain
$$\|B_2\|_{\Gamma_j, p} \leq D_p |\Im (s)|^{\beta_0}\Ep \sum_{k = m +1}^n \theta^k \leq D'_p |\Im (s)|^{\beta_0}\Ep \theta^{m + 1}.$$
So, changing $\theta$ by another global constant $0 < \tilde{\theta} < 1, \: \tilde{\theta} \geq \max{\{\sqrt{\rho}, \sqrt{\theta}}\}$,
we arrange an estimate
$$\|L_s^n \mc_{n,s}w_s\|_{\Gamma_j, p} \leq B_p |\Im (s)|^{\beta_0}\Ep \tilde{\theta}^n.$$
Thus, with global constants $C_p,\: D_p$ we deduce
\begin{equation} \label{eq:5.7}
\|U_{n+2,j}(x, s)\|_{\Gamma_j, p} \leq C_p|\Im (s)|^{\beta_0}\Ep (\theta^n + \tilde{\theta}^n) \leq D_p |\Im (s)|^{\beta_0}\Ep \tilde{\theta}^n, \forall n \in \N.
\end{equation}
Consequently, the series $w_{0, j}(x, s)$ is convergent in $C^p(\Gamma_j)$ norm and for $\sigma_0 \leq \tau \leq s_0 + 1$ we have the estimates
\begin{equation} \label{eq:5.8}
\|w_{0, j}(x, \tau  + \ii t)\|_{\Gamma_j, p} \leq B_p |t|^{\beta_0}\Ep, \: p \geq 1,
\end{equation}
where the constants $B_p$ are independent of $j$. Summing over $\ell = 1,\ldots,\kappa_0$, we obtain the same estimate for 
$\|w_0(x, \tau + \ii t)\|_{\Gamma, p}$ and for $\Re (s) \geq \sigma_0$ the trace $w_0(x, s)$ is an analytic function in $s$ with values in $C^{\infty}(\Gamma).$\

It is interesting to observe that contracting the domain $\sigma_0 \leq \Re (s) \leq s_0 + 1$ we may obtain  better bounds for the $C^p(\Gamma)$ norms. 
For example, we treat below the case $p = 0$ and the same argument works for $p \geq 1.$
In the domain $\sigma_0 \leq \Re( s) \leq s_0 + d,\: d > 0,\:\Im (s) \geq t_0,$ we apply the Phragmen-Lindel\"of theorem (see 5.65 in \cite{kn:T}). 
Notice that when we decrease $d > 0$ the constant $C_{j, d}$ in (\ref{eq:5.6}) change but we  always have the bound (\ref{eq:5.6}). Consequently, for 
$\sigma_0 \leq \tau \leq s_0 + d$ we deduce
$$\|w_{0,j}(x, \tau + \ii t)\|_{\Gamma_j,0} \leq B|t|^{\kappa(\tau)}, \: t \geq t_0,$$
where $\kappa(x)$ is a linear function such that 
$$\kappa(\sigma_0) = 1 + \tilde{\beta}_0\; ,\: \kappa(s_0 + d) = 0.$$
It is clear that choosing $d > 0$ small enough, there exist $\sigma_0'$ with  $\sigma_0 < \sigma_0' < s_0$ and $0 < \beta < 1$ so that for $\tau \geq \sigma_0'$ we have
$$\|w_{0,j}(x, \tau + \ii t)\|_{\Gamma_j,0} \leq A_j\, |t|^{\beta},\: t \geq t_0\;,$$
and similarly we treat the case $t \leq -t_0.$ Finally, for $\tau \geq \sigma_0',\:|t| \geq t_0$ we have
\begin{equation} \label{eq:5.9}
\|w_{0,j}(x, \tau+ \ii t)\|_{\Gamma_j,0} \leq A_j |t|^{\beta}.
\end{equation}
Here the constants $A_j$ depend on the norms of $\nabla \varphi$ and $h.$ 

\ms

\noindent
{\bf Remark 5.} In the following we will not use the estimate (\ref{eq:5.9}), however a similar argument based on Phragmen-Lindel\"of theorem 
will be crucial in Section 7, where we need to control the behavior of the remainder ${\mathcal Q}_M(x, s; k)$ and its bounds when $|\Im (s)| \to \infty.$ 
On the other hand, (\ref{eq:5.9}) is related to  the assumption (\ref{eq:1.6}) of Ikawa mentioned in the Introduction. The estimate (\ref{eq:1.6}) can be 
established choosing $\sigma_0' < s_0$ close to $s_0$ and applying (3.3).  This is not necessary for our exposition and we leave the details to the reader.


\section {The leading term $V^{(0)}(x, s; k)$}
\renewcommand{\theequation}{\arabic{section}.\arabic{equation}}
\setcounter{equation}{0}

\def\do{{\mathcal D}_0}
\def\d1{{\mathcal D}_0}

Our purpose here is to apply the construction in Section 3 with boundary data 
$$m(x,s; k) = e^{\ii k \psi(x)}b(x, s; k),\: x \in \Gamma_j,$$
where $k \geq 1$ and $s \in \d1 = \{ s \in \C: \sigma_0 \leq \Re (s) \leq 1,\: |\Im \:s| \geq J > 0 \}$, with some constant $J$ which will be chosen below. 
We suppose that there exists a phase function $\varphi(x)$ satisfying the condition $(\mathcal P)$ in $\Gamma_j$ such that 
$\varphi(x)\vert_{\Gamma_j} = \psi(x)$ for $x \in {\rm supp}_x\: b(x, s; k)$. The amplitude $b(x, s; k)$ is analytic with respect to $s \in \do$ and 
$\bigcup_{s, k}{\rm supp}_x \: b \subset \Gamma_j,$
$$\|b(x, s; k)\|_{\Gamma_j, p} \leq C_p,\: \forall k \geq 1, \: s \in \do, \: \forall p \in \N.$$ 
In the following we will use the notation $\la z \ra = (1 + |z|)$. For our construction it is convenient to write the oscillatory data $m(x, s; k)$ with phase $e^{-s\psi(x)}$ and we set 
$$m(x, s; k) = e^{-s \psi(x)} e^{(s + \ii k)\psi(x)}b(x, s; k) = e^{-s \psi(x)} b_1(x, s; k).$$
Then
$$\|b_1(x, s; k)\|_{\Gamma_j, p} \leq C_p'\la s + \ii k \ra^p, \forall p \in \N.$$
Notice that our data depends on {\bf two parameters} $s \in \do$ and $k \geq 1.$ The complex parameter $s$ will be related to the convergence of the series $\wj(x,s; k)$ 
constructed in Sect. 5 starting with initial data $m(x, s; k)$,  while the real parameter $k$ is connected with the oscillatory data 
$G(x)e^{\ii k \la x, \eta\ra}\vert_{y \in \Gamma_j},\: |\eta| \leq 1- \delta_1/2 < 1,$ coming from a Fourier transform (see Sect. 8). It is important to note that up to the end of 
Sect. 7 the parameters $s$ and $k$ will not be related and the estimates obtained depend on expressions of the form $\la s + \ii k \ra ^M$. After the application of 
Phragmen-Lindel\"of argument at the end of Sect. 7, we take $|s + \ii k| \leq$ Const in order to get bounds by powers of $k$. We consider amplitudes $b(x, s; k)$ 
depending on $s$ and $k$ to cover higher order approximations in Sect. 7. Starting with boundary data $e^{-s \psi} b_1$ and following the procedure in Sects. 3-5, 
we can justify the convergence of the series $w_{0, j}(x, s; k)$ which are analytic for $s \in \do$.

\def\tc{\tilde{C}}
\def\fs{\widetilde{WF}}

Now we will discuss the domain where the parameter $s$ is running. For $\Im (z) < 0$ we define the resolvent
$(-\Delta_K -z^2)^{-1}$ of the Dirichlet Laplacian $-\Delta_K$ related to $K$ by the spectral calculus and we get
$$\|(-\Delta_K - z^2)^{-1}\|_{L^2(\Omega) \to L^2(\Omega)} \leq \frac{C}{|z| |\Im (z)|},\: \Im (z)< 0.$$
The cut-off resolvent $\psi (-\Delta_K - z^2)^{-1}\psi, \: \psi \in C_0^{\infty}(\Omega),$ has a meromorphic continuation in $\C$ for $N$ odd and in 
$\C \setminus \ii \R^+$ for $N$ even. This resolvent is called {\bf outgoing}.
Setting $z = -\ii s$, we obtain an outgoing resolvent $(\Delta_K - s^2)^{-1}$ which is a bounded operator in $L^2(\Omega)$ for $\Re (s) > 0$ and the 
analytic singularities of $\psi(\Delta_K - s^2)^{-1} \psi$ are included in $\Re (s) < 0$.  Set $\Omega_j = \R^N \setminus K_j$ and suppose that 
$K \subset \{x \in \R^N: \: |x| < \rho_0\}$. Since the real parameter $k \geq 1$ is positive, we assume in this and in the following sections that 
$\Im (s) < 0$. To treat the case $\Im (s) > 0$, we must take $k \leq -1$ and repeat the argument. For our analysis it is more {\it convenient} to consider 
the outgoing resolvent ${\mathcal R}(s)$ acting on functions $f \in H^2(\Gamma)$ defined for $s$ outside the set of resonances (and also for $s \notin \ii \R^{+}$ for $N$ even). 
More precisely, given $f \in H^2(\Gamma)$ we define ${\mathcal R}(s)f = v(x,s)$, where $v(x, s)$ is the unique outgoing solution of the problem
$$\begin{cases}  (\Delta - s^2) v = 0,\: x \in \ovo,\cr v\vert_{\Gamma} = f.\end{cases}$$
Here outgoing means that
$$v(r\theta) = r^{-\frac{N-1}{2}} e^{- s r}(w(\theta) + o(1)),\: \partial_r v + s v = o(1)v,\: r \to +\infty$$
uniformly with respect to $\theta  \in S^{N-1}$ with some $w \in C^{\infty}(\S^{N-1})$. This condition is equivalent to
\begin{equation} \label{eq:6.1}
v\vert_{|x| \geq \rho_1} = (S_0(s) u)\bv_{|x| \geq \rho_1}
\end{equation}
for some $\rho_1 >> \rho_0$ and a compactly supported (in a compact set independent of $s$) function $u$, where
$$S_0(s) = (\Delta - s^2)^{-1}: L^2_{\rm comp}(\R^N) \longrightarrow H^2_{\rm loc}(\R^N)$$
is the outgoing resolvent of the Laplacian in $\R^N.$
If we replace above $K$ by the strictly convex obstacle $K_j$, we can choose $J \geq 2$ so that the outgoing resolvents
$${\mathcal R}_j(s): H^{p + 2}(\Gamma_j) \to H^{p + 1}(\Omega_j \cap \{|x| \leq R\}),\: p \in \N$$
are analytic  (see \cite{kn:V}, \cite{kn:G}) for  
$$s \in \d1 = \{s \in \C: \: \sigma_0 \leq \Re (s) \leq 1,\: |\Im (s)| \geq J\}.$$ 
and $w_j = {\mathcal R}_j(s) f$ is outgoing solution of the problem
$$\begin{cases} (\Delta - s^2) w_j = 0,\: x \in \Omega_j,\cr w_j\bv_{\Gamma_j} = f.\end{cases}$$
Moreover, for $s \in \d1$ and $R \geq \rho_0 + 1$ we have the estimate
\begin{equation} \label{eq:6.2}
\|{\mathcal R}_j(s) f \|_{H^{p + 1}(\Omega_j \cap \{|x| \leq R\})} \leq C_{R, p}\langle s \rangle^{p +2}\|f\|_{H^{p +2}(\Gamma_j)},\: j = 1,...,\kappa_0\;,
\end{equation}
with some constant $C_{R, p} > 0$. The above estimate was established for $p = 0$ in Proposition A. II. 2 in \cite{kn:G}. 
For the sake of completeness we give the argument for $p \geq 1.$ Let $\chi \in C_0^{\infty}(\R^N)$ be a cut-off function such that $\chi(x) = 1$ for $|x| \leq R$ and $\chi(x) = 0$ for 
$|x| \geq R + 1.$ Set $w_j = {\mathcal R}_j(s) f$ and observe that
$$\Delta(\chi w_j) = 2<\nabla \chi, \nabla w_j> +s^2 \chi w_j + \Delta(\chi) w_j = F_j.$$
The function $\chi w_j$ is a solution of the Dirichlet problem in $\omega_R = (|x| \leq R + 1) \cap \Omega_j$ and the standard estimates for boundary problems imply
$$\|\chi w_j\|_{H^{2}(\omega_R)} \leq C_{R, 2} \Bigl( \|F_j\|_{L^2(\omega_R)} + \|f\|_{H^{3/2}(\Gamma_j)}\Bigr).$$
To estimate $\|\chi w_j\|_{L^2(\omega_R)}$, write $w_j = e(f) - (\Delta_{K_j} - s^2)^{-1}(\Delta - s^2)e(f)$, where $e(f)$ is extension operator from $H^2(\Gamma_j)$ to 
$H^{5/2}_{\rm comp}(\omega_{R-1}).$ This implies $\|\chi w_j\|_{L^2(\omega_R)} \leq B_R \langle s \rangle \|f\|_{H^2(\Gamma_j)}$, since for strictly convex obstacles we have (see for instance, Chapter X in \cite{kn:V})
$$\|\chi (\Delta_{K_j} - s^2)^{-1}\chi\|_{L^2 \to L^2} \leq C\langle s \rangle ^{-1}.$$
In the same way one estimates $\|\Delta(\chi) w_j\|_{L^2(\omega_R)}$ by using another cut-off, and applying  (\ref{eq:6.2}) for $p = 0$ we obtain this estimate for $p = 1$. 
The general case can be considered by using an inductive argument. More precise estimates than (\ref{eq:6.2}) can be obtained following a construction of outgoing parametrix for the Dirichlet problem outside $K_j$ (see Appendix II in \cite{kn:G}). 

Finally, notice that for $v$ with supp $v \subset \{|x| \leq R\}$ we have the estimates  (see \cite{kn:V})
\begin{equation} \label{eq:6.3}
\|S_0(s) v\|_{H^{p+1}(|x| \leq R)}\leq C_{R, p} \|v\|_{H^p(|x| \leq R)},\:p \in \N,\; s \in \d1.
\end{equation}

For our construction we need to introduce some pseudodifferential operators depending on the parameter $s \in \d1$. 
For this purpose we will use the notation and the results in Appendices A.I, II in {\cite{kn:G}} (see also \cite{kn:SV}, Appendix). Given a set $X \in \R^{N-1}$, 
we denote by $\tc^{\infty}(X)$ the space of the functions $u(x, s),\: s \in \d1,$ such that $u(. , s) \in C^{\infty}(X)$ and $p(u(., s)) = {\mathcal O}(\la s \ra^{-\infty})$ 
for all seminorms $p$ in $C^{\infty}(X)$. In a similar way we define distributions $\tilde{D}'(X)$. Next,
given two open sets $X \subset \R^{N-1},  Y \subset \R^{N-1}$, consider the spaces of symbols $a(x, y, \eta, s) \in S_{\rho, \delta}^{m, l}(X \times Y)$ such that
for every compact $U \subset X \times Y$, all multiindices $\alpha, \beta, \gamma$ and $s \in \d1$ we have
$$\sup_{(x, y) \in U} |\partial_x^{\alpha} \partial_y^{\beta} \partial_{\eta}^{\gamma} a(x,y,\eta,s)| \leq C_{\alpha, \beta, \gamma, U} |s|^{l + \rho|\gamma| + 
\delta|\alpha + \beta|}(1 +|\eta|)^{m - |\gamma|}.$$
Consider the pseudodifferential operator $Op(a) \in L^{m, l}_{\rho, \delta}(X)$ defined by
$$(Op(a)u)(x, s) = \Bigl(\frac{s}{2\pi}\Bigr)^{N-1} \int e^{-s \langle x - y, \eta \rangle} a(x, y, \eta, s) u(y, s) dy d\eta,$$
where the support of $a(x, y, \eta, s)\in S_{\rho, \delta}^{m, l}(X \times Y)$ with respect to $(y, \eta)$ is uniformly bounded for $s\in \d1$ and $a(x, y, \eta, s)$ is analytic 
for $s \in \d1$. The operator $Op(a)$ maps $\tc_0^{\infty} (Y)$ into $\tc^{\infty} (X).$ Below we will take $Y = \Gamma_j$ and the symbols $a(x,y,\eta,s)$ 
will have compact supports with respect to $(y,\eta)$. Moreover, we will work with symbols in $S^{m,l}_{0,0}.$ We say that $Op(a)$ is properly supported if the kernel $K(x,y,s)$ 
of $Op(a)$ is properly supported uniformly with respect to $s$. Recall that $K(x,y,s)$ is properly supported if both projections from the support of $K(x,y,s)$ to $X$ and $Y$ are 
proper maps (see Definition 18.1.21 in \cite{kn:H}).  We refer to Appendix A. I in \cite{kn:G} for the properties of pseudodifferential operators depending on $s$. Notice that a 
properly supported pseudodifferential operator $Op(a)$ can be defined also by a symbol $a(x, \eta, s)$. A properly supported  pseudodifferential operator $Op(a)$ is called 
{\it elliptic} at $(x_0, \eta_0) \in T^*(X)$ if $a(x,\eta, s)$ satisfies the estimate
$$|a(x,\eta, s)| \geq C\la s \ra^p,\:p \geq 0,\: (x,\eta) \in {\mathcal V},\: s \in \d1,$$
${\mathcal V}$ being a neighborhood of $(x_0, \eta_0)$  independent of $s$.\\

Next, consider Fourier integral operators with real phase function $\varphi(x, \eta)$ and complex parameter $s \in \d1$ having the form
$$I(u)(x, s) = \Bigl(\frac{s}{2\pi}\Bigr)^{N-1} \int e^{-s (\varphi(x,\eta) -\la y, \eta \ra)} a(x, y, \eta, s) u(y, s) dy d\eta,$$
where as above the support of $a(x, y, \eta, s)\in S_{\rho, \delta}^{m, l}(X \times Y)$ with respect to $(y, \eta)$ is uniformly bounded for $s\in \d1$ and $a(x, y, \eta, s)$ is analytic 
for $s \in \d1$. For example, the local parametrix constructed in the hyperbolic region defined below is a Fourier integral operator in this form.

\def\qc{{\mathcal Q}}
\def\uc{{\mathcal U}}

To examine the asymptotic behavior with respect to the parameter $s$  we will use the frequency set $\fs(u)$ introduced in \cite{kn:G} (see also \cite{kn:GS}, \cite{kn:SV}) 
(The notation $\fs(u)$ is used to avoid the confusion with the wave front set $WF(u)$ of a distribution).  We recall the definition of $\fs(u)$ only for the so called {\it finite points} 
$(x, \eta) \in T^*(X)$, since this is sufficient for our argument. Let  $u(x, s) \in \tilde{{\mathcal D}}'(X)$ be a distribution depending on the parameter $s$ so that for every compact 
$X' \subset X$ there exists M such that
$u(x, s)\vert_{X'} \in H^{-M}(X')$ and $\|u(.,s)\bv_{X'}\|_{H^{-M}} \leq C_M \la s \ra^{-M}.$
We say that $(x_0, \eta_0) \in T^*(X)$ is not in $\fs(u)$ if there exists $Op(a) \in L^{0,0}_{\rho, \delta}(X), \: \rho + \delta < 1,$ 
properly supported and elliptic at $(x_0, \eta_0)$ such that for every compact $U \subset X$ we have
$$\|(Op(a)u)(x, s)\|_{C^j(U)} \leq C_{U, M, j}\la s \ra^{-M},\:\forall j\in \N, \: \forall M \in \N,\: s \in \d1.$$
If $\uc$ is a neighborhood of $K$ and if the distribution kernel $Q(x, y, s)$ of an operator $\qc(s): C^{\infty}(\Gamma) \longrightarrow C^{\infty} ({\mathcal U} \setminus K)$ 
belongs to $\tc^{\infty}(\uc \setminus K \times \Gamma)$,  we will say briefly that $\qc(s)u$ is a {\it negligible term}. The terms having behavior ${\mathcal O}(\la s \ra^{-M})$ 
with large $M$ will also be called negligible. It is important to note that a series of negligible terms in general is not negligible, and one needs to have uniform estimates with 
respect to $s$ of the terms of the series to conclude that such a series is negligible. 

\subsection{\underline{ Construction of the operators $P_h,P_g,P_e$}}

In the analysis below we fix $j \in \{1,...,\kappa_0\}.$ Consider the {\it hyperbolic}, {\it glancing} and {\it elliptic sets} on $T^*(\Gamma_j)$ defined respectively by 
$${\mathcal H} = \{(y, \eta) \in T^*(\Gamma_j):\: |\eta| < 1\},\: {\mathcal G} = \{(y, \eta) \in T^*(\Gamma_j): \: |\eta| = 1\},$$
$${\mathcal E} = \{(y, \eta) \in T^*(\Gamma_j):\: |\eta| > 1\},$$
where $(y, \eta)$ are local coordinates in $T^*(\Gamma_j)$.
Let  $\chi_0 \in C_0^{\infty}(T^*(\Gamma_j))$ be a function such that $0 \leq \chi_0 \leq 1$ and $\chi_0(y, \eta) = 0$ in a small neighborhood $G_0$ of 
${\mathcal G} \cup {\mathcal E}$, while $\chi_0(y, \eta) = 1$ for $(y, \theta)  \in G_1,\: G_1 \subset T^*(\Gamma_j) \setminus G_0 \subset {\mathcal H}.$ 
Choosing a finite covering of $\Gamma_j$, we may suppose that in local coordinates $(y, \eta)$ we have $\chi_0(y, \eta) = 1$ for 
$y \in \Gamma_j,\: |\eta| \leq 1 -\delta_1,$ where $\sqrt{1 - \delta_0^2} < 1 - \delta_1 < 1$ and $\delta_0 \in (0,1)$ is a global constant chosen 
as in Lemma 1 (see Sect. 2). Thus if a ray $\gamma_{in}$ issued from $\cup_{\ell \neq j} K_{\ell}$ meets $\Gamma_j$ at $y \in \Gamma_j$ 
with direction $\xi \in \S^{N-1}$ so that $\chi_0(y, \xi\vert_{T_y(\Gamma_j)}) \not= 1,$  then the reflected or diffractive outgoing ray 
$\gamma_{out}$ issued from $(y, \xi - 2 \la \xi, \nu(y) \ra \nu(y))$ does not meet a neighborhood of $\cup_{\nu \neq j} K_{\nu}$ depending only on $\delta_0.$

\def\tc{\tilde{C}^{\infty}}
\def\uc{{\mathcal U}_j}
\def\qc{{\mathcal Q}}

Consider a finite partition of unity of the set ${\rm supp} (\chi_0) \subset {\mathcal H}$ and, as in \cite{kn:G},  a finite partition of unity of psedodifferential 
operators to localize the construction. Let $(y_0, \eta_0) \in {\rm supp} (\chi_0) \subset {\mathcal H}$ and let $\chi(y, \eta) \in C_0^{\infty}(T^*(\Gamma_j)),\: 0 \leq \chi(y, \eta) \leq 1$,
be a function such that $\chi = 1$ in a neighborhood of $(y_0, \eta_0)$. Let $\widetilde{{\mathcal U}_j}$ be a small neighborhood of $K_j$ and let $\uc = \widetilde{{\mathcal U}_j} \setminus K_j$. 
Let $\Gamma_{\chi} \subset \Gamma_j$ be the projection of supp $\chi(x, \eta)$ on $\gj$.

We will omit again the dependence on $k$ in the notation  if the context is clear. Given boundary data $u(y,s)$, in the hyperbolic region we construct an outgoing parametrix
$H_{h, \chi}: \tc(\Gamma_{\chi}) \longrightarrow \tc(\uc)$ of the form

$$(H_{h, \chi}u) (x, s) = \Bigl(\frac{s}{2 \pi}\Bigr)^{N-1} \int e^{-s\Bigl(\psi(x, \eta) - \la y, \eta \ra \Bigr)} \sum_{\nu = 0}^M a_{\nu}(x, y, \eta)s^{-\nu} u(y, s) dy d\eta.$$
We have
$$\begin{cases} (\Delta_x - s^2) (H_{h,\chi}u)(x, s) = s^{-M}A_M(s)u,\:\:  x \in \uc,\\
(H_{h,\chi} u)(x,s)\bv_{\Gamma_j} = Op(\chi)u, \end{cases}$$
where 

$$A_M (s)u= \Bigl(\frac{s}{2 \pi}\Bigr)^{N-1}\int e^{-s(\psi(x, \eta)- \la y, \eta \ra )} (\Delta_x - s^2)\Bigl(a_M(x, y, \eta) \Bigr)u(y, s)dy d\eta.$$
The construction of $H_{h, \chi}$ is given in Appendix A. II. 2 in \cite{kn:G}. Here the phase $\psi(x, \eta)$ satisfies the equation
$$|\nabla_x\psi|^2 = 1,\: \psi\vert_{\Gamma_j} = \la x, \eta \ra,\: (x,\eta)\: \mbox{\text close to}\: (y_0, \eta_0).$$
The amplitudes $a_{\nu}(x,y, \eta)$ are determined from the transport equations with initial data
$$a_0\vert_{x \in \Gamma_j} = \chi(y, \eta),\: a_{\nu}\vert_{x \in \Gamma_j} = 0, \: \nu \geq 1.$$ 
Notice  that $a_{\nu}$ depend only on $\chi(y, \eta)$ and the integration in $H_{h, \chi} u$ is over a compact domain with respect to $y$ and $\eta$, so for $s \in \d1$ the integral is well defined.
 Applying a finite partition of unity, we construct an outgoing parametrix $H_h: \tc(\Gamma_j) \longrightarrow  \tc(\uc)$ such that
$$\begin{cases} (\Delta_x - s^2) (H_h u)(x, s) = s^{-M}B_M(s)u,\:\:  x \in \uc,\\
(H_h u)(x, s)\bv_{\Gamma_j} = Op(\chi_0)u, \end{cases}$$
where the operator $B_M(s)$ is analytic with respect to $s$ and satisfies the estimates 
$$\|B_M(s)u\|_{H^p(\uc)} \leq C_{p}|s|^{p+ 2}\|u\|_{0, \Gamma_j},\: \forall p \in \N$$
 with some global constants.
Let $\Psi(x) \in C_0^{\infty}(\uc)$ be a cut-off function such that $\Psi(x) = 1$ in a small neighborhood of $K_j.$ Then we obtain
$$(\Delta_x - s^2) [\Psi H_h u] = s^{-M} \Psi B_M(s)u + [\Delta, \Psi] H_h u, \: x \in \uc$$
and we define the outgoing parametrix
$$(P_h u) (x, s)= \Psi H_h u - S_0(s)\Bigl(s^{-M} \Psi B_M(s)u + [\Delta, \Psi] H_h u\Bigr),\: x \in \Omega_j.$$
 Thus we get
$$\begin{cases} (\Delta_x - s^2)(P_h u)(x, s) = 0,\: x \in\: \Omega_j,\: s \in \d1,\\
(P_h u)(.,s) \in L^2(\Omega_j) \:{\rm if}\: \Re (s) > 0,\\
(P_h u)(x, s)\bv_{\gj} = Op(\chi_0) u + {\mathcal Q}_h(s) u,
\end{cases}$$
where for large $M$ we obtain a negligible operator ${\mathcal Q}_h(s)$ coming from the trace of the action of $S_0(s).$ Here we use the fact that the frequency set of $S_0(s) w$ is given by the outgoing rays issued from $\fs (w)$ and the outgoing rays issued from $[\Delta, \Psi] H_hu$ do not meet $\gj$. Notice that the operator $P_h$ depends analytically on $s$.\\

Next, let $\chi_1(x, \eta) + \chi_2(x, \eta) = 1 - \chi_0(x, \eta)$, where $\chi_1(x, \eta) \in C_0^{\infty}(T^*(\Gamma_j))$ is a function with support in 
$\{(x, \eta):\:1 - \delta_1 \leq 1 - 2 \epsilon_0 \leq |\eta| \leq 1 + 2 \epsilon_0\}$, while $\chi_2(x, \eta) \in C^{\infty}(T^*(\Gamma_j))$ has support in 
$\{(x, \eta):\:|\eta| \geq 1 + \epsilon_0\},\: \epsilon_0 > 0$ being small enough.
In the glancing region following the construction in Appendix A. II. 3 in \cite{kn:G} and Appendix A. 3 in \cite{kn:SV}), we construct an outgoing parametrix $H_g$ such that
$$\begin{cases} (\Delta_x - s^2) (H_g u) = s^{-M}B_g(s) u,\:\:  x \in \uc,\\
(H_g u)(., s) \in L^2(\Omega_j) \:{\rm if}\: \Re (s) > 0,\\
H_g u\vert_{\Gamma_j} = Op(\chi_1)u + s^{-M}B_g'(s) u \;, \end{cases}$$
where $B_g(s)$ and $B_g'(s)$ are Fourier-Airy operators with complex parameter. The only difference with the construction in \cite{kn:G} is that we have $s^{-M}B_g(s)$ and $s^{-M}B_g'(s)$ instead of operators with kernel in $\tc(\uc \times \Gamma_j)$ and $\tc(\Gamma_j \times \Gamma_j)$, respectively. For this purpose, as in the hyperbolic case, we use a finite sum of amplitudes instead of an asymptotic infinite sum of symbols. The advantage is that our parametrix $H_g$, as well as $B_g(s)$ and $B_g'(s)$, depend analytically on $s$. 
Now define
$$(P_g u) (x, s)= \Psi H_g u - S_0(s)\Bigl(s^{-M}\Psi B_g(s) u + [\Delta, \Psi] H_g u\Bigr),\: x \in \Omega_j.$$

In the elliptic region the construction of a parametrix in Appendix A. II. 4, \cite{kn:G} is given by a Fourier integral operator with big parameter $\lambda$ and complex phase function. 
When $\lambda$ is complex, there are some difficulties to justify this construction (see Appendix A.4 in \cite{kn:SV}). For this reason in the elliptic region we 
introduce $P_e u= {\mathcal R}_j(s)\Bigl(Op(\chi_2)u\Bigr)$ keeping the analytic dependence on $s$.\\
\def\sj{{\mathcal S}_j(s)}
\def\rj{{\mathcal R}_j(s)}
Thus setting ${\mathcal S}_j(s) = P_h + P_g + P_e$, we have
$$\begin{cases} (\Delta_x - s^2)(\sj u)(x, s) = 0,\: x \in\:\Omega_j,\: s \in \d1,\cr
(\sj u)(.,s) \in L^2(\Omega_j) \:{\rm if}\: \Re (s) > 0,\cr
(\sj u)(x, s)\bv_{\gj} = u + {\mathcal Q}_j(s) u
\end{cases}$$
where for large $M$ the operator ${\mathcal Q}_j(s)$ is negligible.\\

Our strategy is to apply the above construction to the function 

$$w_{0,j}(x, s) = \sum_{n = n_j}^{\infty} U_{n+2, j}(x, s)\bv_{\Gamma_j},$$
where
$$U_{n+2, j}(x, s) = \sum_{|\jj| = n + 3, j_{n+2} = j} u_{\jj}(x, s)$$
and $u_{\jj}(x, s)$ are defined in Sect. 3 starting with initial data $e^{-s\varphi}b_1(x, s; .)$. 
Recall that in the previous section we obtained estimates for the $C^p(\Gamma_j)$ norms of $U_{n+2, j}(x, s)$ for $s \in \d1.$ Thus applying $P_h, P_g$ and $P_e$ to $\wj(x, s)$ 
we obtain convergent series. 
Consequently, the function $(\sj\wj)(x, s)$ is analytic for $s \in \d1$ with values in $C^{\infty}(\overline{\Omega_j})$ and here we use the fact that $\wj(x, s) \in C^{\infty}(\Gamma_j)$. It is convenient to introduce the following 

\begin{deff} Let $\omega \subset \R^N$ be an open set and let
${\mathcal D}$ be a domain in $\C$. We say that the function $U(x, s; k)$ satisfies the condition {\rm (S)} in $(\omega, {\mathcal D})$ if 
the following hold:\\

(i) for $k \geq 1$, $U(., s; k)$ is a $C^{\infty}(\overline{\omega})$-valued analytic function in ${\mathcal D}$,\\

(ii) $U(., s; k) \in L^2(\omega) \: {\text for}\: \Re\: s > 0,$\\

(iii) $(\Delta_x - s^2)U(x, s; k) = 0\;\;\; {\text in}\;\;\; \omega$ for every $s \in {\mathcal D}.$
\end{deff}

It is clear that  $(\sj(s) \wj)(x, s)$ satisfies the condition (S) in $(\Omega_j, \d1).$ Taking the sum over $j = 1,\ldots,\kappa_0$, we conclude that the function
$$V^{(0)}(x, s) = \sum_{j=1}^{\kappa_0} (\sj\wj)(x, s)$$
satisfies the condition (S) in $(\ovo, \d1)$.

\subsection{\underline{Traces of $\sj \wj$ on $\Gamma_{\ell}$}}

The analysis of the traces $\Bigl(\sj \wj\Bigr)(x, s)\bv_{\Gamma_{\ell}},\: \ell \neq j,$ is more difficult.
The main contributions come from $(P_h\wj)\bv_{\Gamma_{\ell}},\: \ell \neq j.$ 
Our goal is to find the leading term of  $P_h\Bigl(U_{n + 2, j}(x, s)\bv_{\Gamma_j}\Bigr)\bv_{\Gamma_{\ell}},\: \ell \neq j$. Let $\jj$ be a configuration such that $|\jj| = n + 3,\: j_{n+2} = j$ 
and let $e^{-s \varphi_{\jj}(x)} a_{\jj}(x, s)$ be a term in $U_{n + 2, j}(x, s).$ 
For $x \in \Gamma_j$ consider
$$Op(\chi_0) \Bigl(e^{-s \varphi_{\jj}(x)} a_{\jj}(x, s)\vert_{\Gamma_j}\Bigr) = \int e^{-s(\la x -  y, \eta \ra +  \varphi_{\jj}(y))} \chi_0(y,\eta) a_{\jj}(y, s) dyd\eta$$
$$= \sum_{\mu = 1}^T \int e^{-s(\la x -  y, \eta \ra +  \varphi_{\jj}(y))} \chi_0(y,\eta) a_{\jj}(y, s) \beta_{\mu}(y, \eta)dyd\eta = \sum_{\mu = 1}^T I_{\mu}(x, s),$$
where $\beta_{\mu} \in C_0^{\infty}(T^*(\Gamma_j))$ are cut-off functions such that $\sum_{\mu = 1}^T\beta_{\mu}(y, \eta) = 1$ for $(y, \eta) \in {\rm supp} \chi_0(y, \eta).$

For $I_{\mu}(x, s)$ we will apply the stationary phase argument with big complex parameter $s \in \d1$ (see for instance, Lemma 2.3 in \cite{kn:G}). 
The critical points of $I_{\mu}(x, s)$ satisfy the equations $x = y,\: \eta = \nabla_y \varphi(y),$  the matrix
$$G_{\jj}(y) = \Bigl( \begin{matrix} \varphi_{\jj,y, y}& & - I \\ -I & & 0\end{matrix}\Bigl)$$
is invertible and  we have
$$(G_{\jj}(y))^{-1} = \Bigl(\begin{matrix} 0 & & -I\\ -I & & -\varphi_{\jj, y, y}\end{matrix}\Bigr).$$ 
An application of the stationary phase argument yields
\begin {eqnarray}\label{eq:6.5}
Op(\chi_0) \Bigl(e^{-s \varphi_{\jj}(x)} a_{\jj}(x, s)\vert_{\Gamma_j}\Bigr)  = e^{-s \varphi_{\jj}(x)} \Bigl[ \chi_0(x, \nabla_y \varphi_{\jj}(x)) a_{\jj}(x, s)\nonumber\\
+ \sum_{q = 1}^{M - 1} L_{q, \jj}(y, D_{y}, D_{\eta})(\chi_0 a_{\jj})(x, \nabla_y \varphi_{\jj}(x)) s^{-q} + A_{M,\jj}(x, s) s^{-M}\Bigr], \: x \in \Gamma_j.
\end{eqnarray}
Here $L_{q, \jj}(y, D_y, D_{\eta})$ are operators of order $2q$ and the form of $(G_{\jj}(y))^{-1}$ shows that all terms in
$L_{q, \jj}$ contain derivatives with respect to one of the variables $\eta_i,\: i = 1,...,N-1.$ Thus, the terms in (\ref{eq:6.5}) with coefficients $s^{-q},\: 1 \leq q \leq M-1,$ vanish if  $ |\nabla_y \varphi_{\jj}(x)| \leq 1 - \delta_1.$\\

For $s \in \d1$ we have
$$P_h\Bigl[\sum_{j = 1, \: j \neq \ell}^{\kappa_0}U_{n+2, j}\bv_{\Gamma_j}\Bigr] =
\rj\Bigl[\Bigl(Op(\chi_0) + \qc_h(s)\Bigr)\Bigl(\sum_{j = 1, \: j \neq \ell}^{\kappa_0}U_{n+2, j}\bv_{\Gamma_j}\Bigr)\Bigr]$$
and for large $M$, the operator $\qc_{h, j, \ell}u = (\rj\qc_h(s)u)\bv_{\Gamma_{\ell}},\: j \neq \ell,$ is negligible.

The leading contribution in the traces on $\Gamma_{\ell}$ comes from the trace of the terms 
$${\mathcal R}_j(s) \Bigl(e^{-s \varphi_{\jj}(x)} \chi_0(x, \nabla_y\varphi_{\jj}(x)) a_{\jj}(x, s)\bv_{\Gamma_j}\Bigr),$$
that is from the action of $\rj$ on the leading term in (\ref{eq:6.5}). To examine this contribution 
we construct, as in Sect. 4 in \cite{kn:I3}, an asymptotic outgoing {\it global} solution
$$v_{\jj, M}(x, s)  = e^{-s\psi_{\jj}(x)} \sum_{\mu = 1}^{M} c_{\jj, \mu}(x, s) s^{-\mu}$$ 
of the problem
$$\begin{cases} (\Delta_x - s^2)v_{\jj, M}(x, s) = s^{-M}r_{\jj, M}(x, s),\: x \in \Omega_j,\\
v_{\jj, M}(x, s)\bv_{\Gamma_j} = e^{-s \varphi_{\jj}(x)} \chi_0(x, \nabla_y\varphi_{\jj}(x)) a_{\jj}(x, s) \bv_{\Gamma_j}. \end{cases}$$
We have $\psi_{\jj}(x) = \varphi_{\jj}(x)$ on $\Gamma_j$ and the phase $\psi_{\jj}(x)$ is defined following the procedure in Sect. 2. Moreover, $\psi_{\jj}(x)$ satisfies the 
condition $({\mathcal P})$ on $\Gamma_j$. Next, the amplitudes  $c_{\jj, \mu}(x, s)$ are determined globally by the transport equations. 
  It is easy to see that 
$$c_{\jj, 0}(x, s) \vert_{\Gamma_{\ell}} = -a_{(\jj, \ell)}(x, s)\vert_{\Gamma_{\ell}},\:\ell \neq j,$$
where $(\jj, \ell)$ is the configuration $(j_0,j_1,...,j_{n+2} = j, \ell)$. This follows from the definition of $a_{(\jj, \ell)}(x, s)$ in Sect. 3 and from the transport equation for the leading term $c_{\jj, 0}$ (see Section 4 in \cite{kn:I3}) combined with the 
fact that if $c_{\jj, 0}(x, s)\vert_{\Gamma_{\ell}} \neq 0$, then $x$ must lie on a ray issued from $(y, \nabla_y \varphi_{\jj}(y))$ with 
$\chi_0(y, \nabla_y \varphi_{\jj}(y)) = 1.$ The sign (-) appears since for the configurations $(\jj, \ell)$ we have to include the factor $(-1)^{n+4}.$ 
Next, we choose a function $\Phi \in C_0^{\infty}(|x| \leq \rho_0 + 1)$ which is equal to 1 in a neighborhood of $K$ and introduce
$$V_{\jj, M}(x, s) = \Phi v_{\jj,M}(x, s) - S_0(s) \Bigl(s^{-M} r_{\jj, M}(x,s) + [\Delta, \Phi] v_{\jj, M}(x, s)\Bigr).$$
We have $(\Delta_x - s^2) V_{\jj, M}(x, s) = 0$ in $\Omega_j$ and for $M$ large  the traces  
$$V_{\jj, M}(x, s)\bv_{\Gamma_{\ell}} -\Bigl({\mathcal R}_j(s) \Bigl[\Bigl(e^{-s \varphi_{\jj}(x)} \chi_0(x, \nabla_y\varphi_{\jj}(x)) 
a_{\jj}(x, s)\Bigr)\bv_{\Gamma_j}\Bigr]\Bigr)\bv_{\Gamma_{\ell}},\: \ell = 1,...,\kappa_0$$
are negligible terms coming from the action of $S_0(s)$. We obtain this first for the trace on $\gj$ and then use the estimates for 
the resolvent ${\mathcal R}_j(s).$  On the other hand, for large $M$ we get $V_{\jj, M}(x, s)\bv_{\Gamma_{\ell}} = v_{\jj, M}(x, s)\bv_{\Gamma_{\ell}}$ 
modulo negligible terms related to the action of $S_0(s).$ Thus the leading 
term of the trace on $\Gamma_{\ell}$ is $e^{-s \varphi_{\jj}(x)} c_{\jj, 0}(x,s)\bv_{\Gamma_{\ell}}.$ 

Next, consider $e^{-s \varphi_{\jj}(x)} b_{\jj}(x, s)\bv_{\Gamma_j}$ with $b_{\jj}(x, s)\vert_{\Gamma_j} = 0$ for 
$|\nabla_y \varphi_{\jj}(x)| \leq 1 - \delta_1$. Moreover, assume that if $b_{\jj}(x, s) \neq 0$ for $x \in \Gamma_j$, then 
$x$ is lying on a segment issued from some obstacle $K_{\ell},\: \ell \neq j.$ From (\ref{eq:6.5}) we see that the terms with 
coefficients $s^{-q},\: 1 \leq q \leq M-1,$ have these properties. According to Theorem A. II. 12 in \cite{kn:G}, the frequency set of 
${\mathcal R}_j(s) \Bigl(e^{-s \varphi_{\jj}(x)} b_{\jj}(x, s) \vert_{\Gamma_j}\Bigr)$ is included in the set determined by the outgoing rays issued from 
$\fs \Bigl(e^{-s \varphi_{\jj}(x)} b_{\jj}(x, s) \bv_{\Gamma_j}\Bigr)$. According to Lemma 1, our choice of $\delta_1$ shows that these rays do not meet a 
neighborhood of $\cup_{\ell \neq j} K_{\ell}.$ Consequently, the traces of ${\mathcal R}_j(s) \Bigl(e^{-s \varphi_{\jj}(x)} b_{\jj}(x, s) \vert_{\Gamma_j}\Bigr)$ on 
$\Gamma_{\ell},\: \ell \neq j,$ are negligible. It is clear also that all terms with factors $s^{-q}$ will produce traces with this factor.\\

 For fixed $n$ and fixed $j, \ell \neq j$ we take 
the finite sum over the configurations $|\jj| = n+3$ of all terms having coefficient $s^{-q},\:1 \leq  q \leq M,$ in the trace 
$\rj \Bigl(Op(\chi_0)U_{n+2, j}\bv_{\gj}\Bigr)\bv_{\Gamma_{\ell}}$ and we denote this sum by $s^{-1} R_{h, n, j, \ell}(x,s)$. 
Since we cannot estimate directly the series with the contributions $s^{-q},$ we are going to include in $s^{-1} R_{h, n, j, \ell}(x, s)$ 
all terms mentioned above as negligible and appearing with coefficients $s^{-q},\: 1 \leq q \leq M.$ 

Thus for fixed $n$, summing over $j = 1,...,\kappa_0, \: j \not= \ell$ and $\jj$, we obtain all configurations $\jj$ with $|\jj| = n + 4,\: j_{n+3} = \ell$ and we conclude that
\begin{eqnarray} \label{eq:6.6}
\Bigl(P_h\sum_{j = 1, \: j \neq \ell}^{\kappa_0}U_{n+2, j}\bv_{\Gamma_j}\Bigr)\bv_{\Gamma_{\ell}} 
= -\sum_{|\jj| = n + 4, \j_{n+3} = \ell} e^{-s \varphi_{\jj}(x)} a_{\jj}(x,s)\bv_{\Gamma_{\ell}} + s^{-1} R_{h, n, j,\ell}(x, s)\\ \nonumber
+{\mathcal Q}_{h, j, \ell}\Bigl(\sum_{j = 1, \: j \neq \ell}^{\kappa_0}U_{n+2, j}\bv_{\Gamma_j}\Bigr).
\end{eqnarray}

To treat $(P_g\wj) \vert_{\Gamma_{\ell}},\: \ell \neq j$, we apply the same argument. Observe that according to the results in Appendix II in \cite{kn:G}, the frequency set of 
$\rj \Bigl(Op(\chi_1)U_{n + 2, j}(x,s)\bv_{\Gamma_j}\Bigr)$ is related to the outgoing rays issued from the frequency set of 
$$Op(\chi_1) \Bigl(\sum_{|\jj| = n + 3,\: j_{n +2} = j}e^{-s \varphi_{\jj}(x)} a_{\jj}(y, s)\vert_{\Gamma_j}\Bigr).$$
For every $\jj$ the frequency set of $Op(\chi_1) \Bigl(e^{-s\varphi_{\jj}(y)} a_{\jj}(y,.)\bv_{\gj}\Bigr)$ is given by 
$(y, \nabla_y \varphi_{\jj}(y))$ such that 
$$y \in {\rm supp}\:\: a_{\jj}(y, .)\bv_{\gj},\: |\nabla_y\varphi_{\jj}(y)| \geq 1 - \delta_1.$$
If $y \in \gj$ has this property and $a_{\jj}(y, .)\bv_{\gj} \neq 0$ for some configuration $\jj$, then $y$ is lying on a segment issued from some $\Gamma_{\mu}, \: \mu \neq j.$ Our choice of $\delta_1$ guarantees that the  outgoing rays mentioned above pass outside a neighborhood of $\cup_{\ell \neq j} K_j$. 
Thus, we deduce
\begin{equation} \label{eq:6.6g}
\Bigl(P_g\sum_{j = 1, \: j \neq \ell}^{\kappa_0}U_{n+2, j}\bv_{\Gamma_j}\Bigr)\bv_{\Gamma_{\ell}} 
= s^{-M} R_{g, n, j,\ell}(x, s).
\end{equation}
Here the series $\sum_{n = 0}^{\infty} R_{g, n, j, \ell}$ is convergent but we cannot show that
$s^{-M} \sum_{n =0}^{\infty} R_{g, n, j, \ell}$ is negligible. In fact, the results of Theorem 3 cannot be applied to this series and for this reason we take $M =1$ in (\ref{eq:6.6g}) and consider $R_{g, n, j, \ell}$ together with the terms $R_{h, n, j, \ell}.$ A similar analysis can be applied to 
$\rj \Bigl(Op(\chi_2)U_{n+2, j}\bv_{\gj}\Bigr)\bv_{\Gamma_{\ell}}$ since there are no outgoing rays issued from the elliptic region, and we obtain
$$\Bigl(\rj \Bigl( Op(\chi_2)U_{n+2, j}\Bigr)\bv_{\gj}\Bigr)\bv_{\Gamma_{\ell}}  = {\mathcal Q}_{e, j, \ell} \Bigl(U_{n+2, j}\bv_{\gj}\Bigr),$$
where the operator ${\mathcal Q}_{e, j, \ell}$ has kernel in $\tc(\Gamma_{\ell} \times \gj).$\\

Summing over $n$ and $j = 1,...,\kappa_0$, we conclude that for $x \in \Gamma$ we have
\begin{equation} \label{eq:6.7}
V^{(0)}(x, s; k) = m(x, s; k) + s^{-1} R_1(x, s; k) + s^{-M}{\mathcal Q}_{M, 0}(x, s; k),
\end{equation}
where in the notations the dependence on $k$ is involved. The cancellation of the leading terms follows from the equality
$$\Bigl(a_{(\jj, \ell)}(x, s) + a_{\jj}(x, s)\Bigr)\bv_{x \in \Gamma_{\ell}} = 0,\: \ell \neq j\;,$$
 and the  representation (\ref{eq:6.6}).  The negligible terms coming 
from the action  of ${\mathcal Q}_{h, j, \ell},\: {\mathcal Q}_{e, j. \ell},\:j, \ell = 1,...,\kappa_0$ to $\wj$ are included in $s^{-M}{\mathcal Q}_{M, 0}(x, s; k),$ while $R_1(x, s; k)$ is the
 sum over $n$, $j$ and ${\ell}$ of the contributions $R_{h, n,j,\ell}(x, s; k)$ and $ R_{g, n, j, \ell}(x, s; k)$ coming from (\ref{eq:6.6g}) with $M = 1.$
 Applying the estimates for $U_{n+2, j}\bv_{\Gamma_j}$ and the analyticity of $P_h, P_g$ and $P_e$,  we deduce that $Q_{M, 0}(x, s; k)$ and $V^{(0)}(x, s; k)\bv_{\Gamma}$ are analytic for $s \in \d1.$  Thus we conclude that $R_1(x, s; k)$ is analytic for $s \in \d1.$ We can prove directly that $R_1(x, s; k)$ is analytic examining the series
$$\sum_{n = n_j}^{\infty} P_{h, n, j, \ell}(x, s; k),\: \sum_{n = n_j}^{\infty}P_{g, n, j, \ell}(x, s; k).$$
In fact, it suffices to obtain estimates $|P_{h, n, j, \ell}| \leq B_{h, j, \ell} \tilde{\theta}^n, \forall n \in \N,$ and we treat this question in the next subsection. Thus the analyticity of $R_1(x, s; k)$ is not related to the analyticity of $V^{(0)}$ and $Q_M$ and we may work with a parametrix $P_e$ which is not analytic in $s$ (see Appendix A. 4 in \cite{kn:SV} and Sect. 8). This could simplify a little bit our argument, but we arrange $V^{(0)}$ to be analytic in order to have similarity with the construction in \cite{kn:I3}. On the other hand, to obtain estimates for the outgoing resolvent better than (\ref{eq:6.2}) we must use an approximation by a parametrix.

\subsection{\underline{Estimates of $R_1(x, s; k)$}}

To estimate $R_1(x, s; k)$ we need to estimate $R_{h, n, j, \ell}$ and $R_{g, n, j, \ell}$. To deal with $R_{h, n, j, \ell}$, we use the equality (\ref{eq:6.6}). Notice that the trace
$$\Bigl(P_h\sum_{j = 1, \: j \neq \ell}^{\kappa_0}U_{n+2, j}\bv_{\Gamma_j}\Bigr)\bv_{\Gamma_{\ell}} $$
 is given by the trace on $\Gamma_{\ell}$ of
$$S_0(s) \Bigl[\Bigl (s^{-M} B_M(s) + [\Delta, \Psi]H_h\Bigl) \sum_{j = 1, \: j \neq \ell}^{\kappa_0}U_{n+2, j}\bv_{\Gamma_j}\Bigr].$$
The term involving $s^{-M}$ is easy to be handled, and we treat the term with $[\Delta, \Psi]$.  Applying the estimates (\ref{eq:5.7}) with $p = 0$ and applying the $L^2$ estimates for the action of the Fourier integral operator $H_h$, we get
$$\|[\Delta, \Psi] H_h\Bigl(\sum_{j = 1, \: j \neq \ell}^{\kappa_0}U_{n+2, j}\Bigr)\bv_{\Gamma_j}\|_{0} \leq C_{j, \ell}|s|^{2 + \beta_0}\la s + \ii k \ra  \tilde{\theta}^n,$$
where $\beta_0$ and $ 0 < \tilde{\theta} < 1$ were introduced in Sect. 5 and $\la s + \ii k\ra$ comes from (\ref{eq:5.7}).
Next for $g\in C^0(\R^N)$ with compact support we write $S_0(s)g = E_s * g$, where $E_s(x)$ is the kernel of $S_0(s)$. This kernel has the form
$$E_s(x) = \frac{\ii}{4} \Bigl(\frac{s}{2 \pi |x|}\Bigr)^{\gamma} H_{\gamma}^{(1)}(s |x|),\: \gamma = (N - 2)/2,$$
where $H_{\gamma}^{(1)}(z)$ is the Hankel function of first type.
Since $\Gamma_{\ell} \cap {\rm supp}\: \Psi = \emptyset$, we can estimate the $C^p$ norms of $\Bigl(S_0(s)[\Delta, \Psi]w\Bigr)\bv_{\Gamma_{\ell}}$ exploiting the estimates for the derivatives of $H_{\gamma}^{(1)}(z)$. Thus setting $\beta_N = (N - 3)/2 + \beta_0,$ we deduce
\begin{equation} \label{eq:6.s0}
\|S_0(s) [\Delta, \Psi] H_h \sum_{j = 1, \: j \neq \ell}^{\kappa_0}U_{n+2, j}\bv_{\Gamma_j}\|_{\Gamma_{\ell}, p} \leq B_{j, \ell, p} \la s + \ii k \ra |s| ^{2 + p+ \beta_N} \tilde{\theta}^n.
\end{equation}

Next for the sum 
$\sum_{|\jj| = n + 4, \j_{n+3} = \ell} e^{-s \varphi_{\jj}(x)} a_{\jj}(x,s)\bv_{\Gamma_{\ell}}$ in (\ref{eq:6.6})
we apply Theorem 3(b). Consequently, summing over $n$, we obtain  estimates for $s^{-1}\sum_{n = n_j}^{\infty} P_{h, n, j, \ell}$ with the same order as in (\ref{eq:6.s0}).\\

The analysis of $R_{g, n, j, \ell}$ is very similar. To estimate $[\Delta, \Psi] H_g (\sum_{j = 1, \: j \neq \ell}^{\kappa_0}U_{n+2, j}\bv_{\Gamma_j})$, we observe that outside a small neighborhood of $K_j$ the parametrix $H_g$
in the glancing domain can be written as a Fourier integral operator with real phase and we may estimate
$\Bigl(S_0(s)|\Delta, \Psi] H_g w\Bigr)\bv_{\Gamma_{\ell}}$ as in the hyperbolic case discussed above. For the remainder ${\mathcal Q}_{0, M}(x, s; k)$ we have
\begin{equation} \label{eq:6.8}
\|{\mathcal Q}_{M, 0}(x, s; k)\|_{\Gamma, p} \leq D_p \la s + \ii k \ra^{p + 2} |s|^{p + 2+ \beta_0}, \: p \in \N,
\end{equation}
 where $\la s + \ii k \ra^{p + 2}$ comes form the estimates of the amplitude $b_1(x, s; k)$.
Finally,  we get the following crude estimates 
\begin{equation} \label{eq:6.9}
\|R_1(x, s; k)\|_{\Gamma, p} \leq C_p \la s + \ii k \ra ^{p + 2} |s|^{p + 3 + \beta_N},\:s \in \d1,\: \forall p \in \N
\end{equation}
and the term $s^{-1} \|R_1(x, s; k)\|_{\Gamma, 0}$ has no order ${\mathcal O}(|s|^{-1})$ for all $s \in \d1.$ 

It is important to note that in the domain of absolute convergence $\Re (s) > s_0 + d > s_0$ we have better estimates for $R_1(x, s; k)$. First, in this domain 
for all $\gamma$ and $|x| \leq R$ the series
\begin{equation} \label{eq:6.10}
D_x^{\gamma}\Bigl(\sum_{n = 1}^{\infty} \sum_{|\jj| = n} e^{-s \varphi_{\jj}(x)} a_{\jj}(x, s)\Bigr)
\end{equation}
are absolutely convergent (see \cite{kn:I3}). Next Proposition 2 shows that the phases $\varphi_{\jj}(x)$ and their derivatives are uniformly bounded with respect to $\jj$ and by 
recurrence we obtain the absolute convergence of the series 
$$\sum_{n = 1}^{\infty} \sum_{|\jj| = n} e^{- s \varphi_{\jj}(x)} L_{q, \jj} (x, D_x) a_{\jj}(x, s),$$
$L_{q, \jj}(x, D_x)$ being  partial differential operators of order $q$ independent on $\jj$ and $n$ with coefficients uniformly bounded with respect to $\jj$. Now in the equality 
(\ref{eq:6.5}) we can sum over the configurations $\jj$ and after the action of ${\mathcal R}_j(s)$ the sum of all terms with coefficients $s^{-q},\: 1 \leq q \leq M-1$, and the remainder yield contributions which can be included in ${\mathcal Q}_{M, 0}.$ To deal with the traces of
$$\sum_{n = 0}^{\infty}\sum_{|\jj| = n + 3, j_{n+2} = j} {\mathcal R}_j(s)\Bigl(\chi_0(x, \nabla_y \varphi_{\jj}(x)) a_{\jj}(x, s)e^{-s \varphi_{\jj}(x)}\bv_{\Gamma_j}\Bigl),$$
we can exploit the estimates in Sects. 4, 5 in \cite{kn:I3} for the amplitudes $c_{\jj, \mu}(x, s)$ of the asymptotic solutions $v_{\jj, M}(x, s)$. In the same way, 
we can estimate and sum the negligible contributions $s^{-M}R_{g, n, j, \ell}$ coming from the glancing region and show that they yield a negligible term. 
Thus, for $\Re (s)> s_0 + d > s_0$  we deduce
\begin{equation} \label{eq:6.12}
\|R_1(x, s; k)\|_{\Gamma, p} \leq C_{p, d}\la s + \ii k \ra^{p + 2} |s|^p, \: p \in \N,
\end{equation} 
while for $|s + \ii k| \leq a + 1$  we obtain
\begin{equation} \label{eq:6.13}
\|R_1(x, s; k)\|_{\Gamma, p} \leq C_{p, d}'k^p, \: p \in \N.
\end{equation}

\section{Higher order terms of the asymptotic solution}
\renewcommand{\theequation}{\arabic{section}.\arabic{equation}}
\setcounter{equation}{0}

\def\tj{\tilde{\jj}}
\def\jl{(j, \ell)}
\def\do{{\mathcal D}_0}
\def\d2{{\mathcal D_2}}
Our purpose is to improve (\ref{eq:6.7}) by higher order approximations $V^{(j)}(x, s; k),\: j = 1,..,M-1$, where $M$ is an integer such that $M > (N-1)/2.$ In particular, for 
$N = 2$ we can take $M = 1$ and the construction in Sect. 6 is sufficient. Recall that the term $R_1(x, s; k)$ in the previous section has the form  
$$\sum_{n = n_j}^{\infty}\sum_{j, \ell = 1}^{\kappa_0} \Bigl(R_{h, n, j, \ell}(x, s; k) + R_{g, n, j, \ell}(x, s; k)\Bigr)$$
with $n_1 = -2$ and $n_j = -1$ for $j \neq 1.$
Fix $j$ and $\ell$ and set
$$e^{-s \varphi_n(x)}m_{1,n}^{\jl}(x,s; k) = R_{h, n, j, \ell}(x, s; k) + R_{g, n, j, \ell}(x, s; k),\: x \in \Gamma_{\ell},$$
where $\varphi_n(x)$ is one of the phases $\varphi_{\jj}(x)$ in $U_{n + 2, j}(x, s; k).$ The choice of $\varphi_n$ is not important and we omit in the notation the 
dependence of $(j, \ell)$. The analysis in the previous section shows that we have the estimates
\begin{equation} \label{eq:7.1}
\|m_{1, n}^{\jl}(x, s; k)\|_{\Gamma_{\ell}, p}  \leq D_p \la s + \ii k \ra^{p + 2}  |s|^{p + 3 + \beta_N}\tilde{\theta}^n, \forall n \in \N,
\end{equation}
where $0 < \tilde{\theta} < 1$ is the same as in Sect. 5. Here and below we denote by $F^{(j,\ell)}$ some terms depending on the traces on $K_j$ 
and $K_l,\: j, \ell = 1,\dots,\kappa_0,$ while $\jj$, $\jj'$ denote configurations. Now for fixed $n$ we apply the construction of Sects. 3 and 6 to the oscillatory data 
$e^{-s \varphi_n(x)}m_{1, n}^{\jl}(x, s; k)$ and we obtain a series $\sum_{m = -1}^{\infty} U_{1, n, m}^{\jl}(x, s; k)$ with 
$$U_{1, n, m}^{\jl}(x, s; k) = \sum_{|\jj'| = m +3, j'_{m+2} = l} (-1)^{m + 2} e^{- s \varphi_{1, n, \jj'}(x)} a_{1, n, \jj'}^{(j, l)}(x, s; k),$$
where the phase functions $\varphi_{1, n, \jj'}(x)$ depend on the configurations $\jj'$. Taking the summation over $n$, we are going to study the double series
\begin{equation} \label{eq:7.2}
w_{1, j, \ell}(x, s; k) = \sum_{n = n_j}^{\infty} \sum_{m = -1}^{\infty} U_{1,n, m}^{(j, \ell)}(x, s; k)\bv_{\Gamma_{\ell}}, \: x \in \Gamma_{\ell}.
\end{equation}
We repeat the argument of the Sect. 5 for $\sigma_0 \leq \Re(s) \leq 1$ and applying (\ref{eq:7.1}) and Theorem 3(b), we get the estimates
\begin{equation} \label{eq:7.3}
\|U_{1, n, m}^{\jl}(x, s; k)\|_{\Gamma_{\ell}, p} \leq D_p' \la s + \ii k \ra^{p+3} |s|^{p + 4 + \beta_N + \beta_0} \tilde{\theta}^{n +m},\: \forall n \in \N, \forall m \in \N\;,
\end{equation}
with constants $D_p'$ independent of $n, m \in \N$. Thus, the double
series defining $w_{1, j, \ell}(x, s; k)$ is convergent.  Applying $S_{\ell}(s)$ to $w_{1, j, \ell}(x,s; k)$ and exploiting (\ref{eq:7.3}), 
we justify the convergence of the corresponding series and for $s \in \do$ we obtain analytic terms. The function
$$V^{(1)}(x, s; k) = -s^{-1}\sum_{j, \ell = 1}^{\kappa_0}{\mathcal S}_{\ell}(s) \Bigl(w_{1, j, \ell}(x, s; k)\Bigr)$$
satisfies the condition (S) in $(\ovo, \do)$  and for $s \in \do$ and $x \in \Gamma$  we get
\begin{equation} \label{eq:7.4}
V^{(0)}(x, s; k) + V^{(1)}(x, s; k) = m(x,s; k) + s^{-2} R_2(x, s; k) + s^{-M} {\mathcal Q}_{M, 1}(x, s; k).
\end{equation}
Here  $R_2(x, s; k)$ and ${\mathcal Q}_{M, 1}(x, s; k)$ are analytic for $s \in \do$, ${\mathcal Q}_{M, 1}$ satisfies the same estimates as in (\ref{eq:6.8}), while for $R_2(x, s; k)$ we have 
\begin{equation} \label{eq:7.5}
\|R_2(x, s; k)\|_{\Gamma, p} \leq C_p \la s + \ii k \ra^{p + 3} |s|^{p + 6 + 2\beta_N},\: \forall p \in \N.
\end{equation}
For $\Re (s) > s_0 +d > s_0$ we obtain again better estimates, since we can choose $\varphi_n(x) = \varphi_{\jj}(x)$ and $m_{1,n}^{\jl} (x, s; k) = c_{\jj, 1}(x, s; k)\bv_{\Gamma_{\ell}},$ 
where $c_{\jj, 1}(x, s; k)$ is the coefficient in front of $s^{-1}$ in the asymptotic solution $v_{\jj, M}(x, s; k)$ introduced in Sect. 6.  Exploiting the convergence of the series (\ref{eq:6.10}), 
we deduce that in this domain the growth in the right hand side of (\ref{eq:7.5}) is $\la s + \ii k \ra^{p + 3}|s|^p.$\\

Repeating this procedure, we construct $V^{(j)}(x, s; k),\:0 \leq j \leq M -1$, which are analytic functions for $s \in \do$ with values in $C^{\infty}(\Omega)$.
They satisfy the condition (S) in $(\ovo, \do)$ and we have
\begin{equation} \label{eq:7.6}
\sum_{j = 0}^{M-1} V^{(j)}(x, s; k) =  m(x, s; k) + s^{-M} {\mathcal Q}_M(x, s; k),\: x \in \Gamma\;,
\end{equation}
with polynomial estimates
\begin{equation} \label{eq:7.7}
\|{\mathcal Q}_M(x, s; k)\|_{\Gamma, 0} \leq C_{M}\la s + \ii k \ra^{L(M)} |s|^{N(M)}, s \in \do.
\end{equation}
Here ${\mathcal Q}_M(x, s; k)$ is analytic for $s \in \do$ and $C_{M}$ depend on the norms of the derivatives of $\psi(x)$ and $b(x, s; k)$ involved in the boundary data 
$m(x, s; k)$ introduced in the beginning of Sect. 6.  Thus, we establish {\bf crude} estimates with orders $N(M), \: L(M)$ depending on $M$ and it seems quite difficult to 
obtain more precise estimates for $s \in \do$. Of course, we have $N(M) > M$, however we will apply the above estimates for fixed $M$ and the precise value of $N(M)$ 
is not important for our argument. For
$\Re (s) \geq s_0 + d > s_0,\: \Im (s) \leq -J$ the absolutely convergence of (\ref{eq:6.10}) implies
\begin{equation} \label{eq:7.8}
\|{\mathcal Q}_M(x, s; k)\|_{\Gamma, 0} \leq C_{M, d} \la s + \ii k \ra^{L(M)}.
\end{equation}
The constant $C_{M, d}$ depends on $d$ but $L(M)$ is independent of $d$.
Now we fix an integer $M \in \N$ so that $M > \frac{N -1}{2}\:,N(M)$ and $L(M)$ are fixed. Next, we fix $d > 0$ small enough so that
$$d \frac{N(M)}{s_0 + d - \sigma_0} <  M- \frac{N-1}{2}.$$

In the domain $\{s \in \C: \: \sigma_0 \leq \Re (s) \leq s_0 + d < 0,\: \Im (s) \leq -J\}$ consider the analytic with respect to $s$ function
$$F(x, s; k) = \frac{{\mathcal Q}_M(x, s; k)}{(s + \ii k)^{L(M)}}.$$
The estimates (\ref{eq:7.7}) and (\ref{eq:7.8})  combined with the Phragmen-Lindel\"of theorem (see \cite{kn:T}) show that for 
$s \in \{ s  \in \C:\: \Re (s) = t,\: \sigma_0 \leq t \leq s_0 + d,\: \Im (s) \leq -J\},$
we have 
$$\|F(x, s; k)\|_{\Gamma, 0} \leq A_M |s|^{\kappa(t)},$$ 
$\kappa(t)$ being  the linear function such that $\kappa(\sigma_0) = N(M),\: \kappa(s_0 + d) = 0.$ We can choose $\sigma_1 < s_0$ so that $0 \leq \kappa(t) \leq \alpha$ 
for $\sigma_1 \leq t \leq s_0 + d$  with some $0 < \alpha < M - \frac{N - 1}{2}.$
Thus, for $\sigma_1 \leq \Re (s) \leq s_0 + d,\: \Im (s) \leq -J, |s + \ii k| \leq |\sigma_0| + 1$ we get
\begin{equation} \label{eq:7.9}
\|{\mathcal Q}_M(x, s; k)\|_{\Gamma, 0} \leq A_M |s + \ii k|^{L(M)}|s|^{\alpha}  \leq B_M k^{\alpha}, \: k \geq 1.
\end{equation}
Moreover, the constant $B_M$ depends on the derivatives of $\nabla \psi$ and $b(x, s; k)$
involved in the boundary data $m(x, s; k)$ as well as on some global constants depending only on $K$.
The restriction $\sigma_1 \leq \Re (s) \leq s_0 + d$ with $s_0 + d < 0$ was used only to guarantee that the factor $(s + \ii k)^{L(M)} \neq 0$ in this domain. For $\Re (s) > s_0 + d$ we can apply the estimate (\ref{eq:7.8}) to obtain (\ref{eq:7.9}) with another constant $A_M$ and $\alpha = 0.$ Consequently, for some fixed $c$ such that $s_0 + c \geq 1$ the estimates (\ref{eq:7.9}) hold for 
$$s \in \cd1 = \{ s \in \C: \: \sigma_1 \leq \Re (s) \leq s_0 + c,\: \Im (s) \leq -J, \:|s + \ii k| \leq |\sigma_0| +c\}.$$

\section{Integral equation on the boundary}
\renewcommand{\theequation}{\arabic{section}.\arabic{equation}}
\setcounter{equation}{0}

\def\cd1{{\mathcal D}_1}
\def\d1{{\mathcal D}_1}
\def\gj{\Gamma_j}

The purpose in this section is to define for $ s \in \cd1$ an operator $R(s, k): \: L^2(\Gamma) \to C^{\infty}(\ovo),$  where  $k >  J + |\sigma_0| + c$ will be taken sufficiently large and $\cd1$ is the domain introduced in the previous section. The operator $R(s, k)$ satisfies
\begin{equation} \label{eq:8.1}
\begin{cases} (\Delta_x - s^2) R(s, k)f = 0,\: x \in \ovo,\cr
R(s, k)f \in L^2(\Omega),\: {\rm for}\: \Re \: (s) > 0,\cr
R(s, k)f \bv_{\Gamma} = f\end{cases}
\end{equation}
and to arrange the boundary condition we will solve an integral equation on $\Gamma$. After the construction of a solution $\sum_{j = 0}^{M-1} V^{(j)}(x, s; k)$ with the properties in 
Sect. 7, it was mentioned in Proposition 2.4 in \cite{kn:I3} that the existence of $R(s; k)$ can be obtained by the argument in \cite{kn:I2}. On the other hand, \cite{kn:I2} deals with
the case of two strictly convex obstacles and in that case the geometry of the trapping rays is rather different from that in \cite{kn:I3} and our paper. For the sake of completeness 
we will discuss briefly how we can construct $R(s, k)$ by using the construction in Sects. 6-7 in the hyperbolic region and those in \cite{kn:I1}, \cite{kn:I3}, \cite{kn:SV} in the glancing 
and elliptic regions.\\

Fix  $M > (N-1)/2$ and  $0 < \alpha < M - \frac{N -1}{2}$ as in the previous section and  $j \in \{1,...,\kappa_0\}$. Let $Y \subset \Gamma_j$ and let  
$F \in L^2(\gj)$ with ${\rm supp}\: F \subset Y.$ As in Sect. 6, choose  local coordinates $(y, \eta)$ in $T^*(Y)$ with  $y = (y_1,\ldots, y_{N -1}) \in W \subset \R^{N-1}\;,$ 
and write
$$F(y) = (2\pi)^{-N+1} \int e^{\ii<y, \eta>} \hat{F}(\eta) d\eta
= \Bigl(\frac{k}{2\pi}\Bigr)^{N-1} G(y) \int e^{\ii k <y, \eta>}\hat{F}(k \eta) d\eta,$$
where $G(y) \in C_0^{\infty}(\R^{N-1}),\,G(y) = 1$ on $\supp \:F(y)$ and
$$\hat{F}(\eta) = \int e^{-\ii<y, \eta>} F(y) dy.$$

Consider a partition of unity $\chi_0(\eta) + \chi_1(\eta) + \chi_2(\eta) = 1$ with $C^{\infty}$ functions $\chi_i(\eta),\: 0 \leq \chi_i(\eta) \leq 1, \: i =0,1,2,$ such that
$${\rm supp}\: \chi_0(\eta) \subset \{\eta:\: |\eta| \leq 1 - \delta_1/2\},\: {\rm supp}\: \chi_1(\eta) \subset \{\eta:\: 1 - \frac{2}{3}\delta_1 \leq |\eta| \leq 1 + \frac{2}{3}\delta_1\},$$
$${\rm supp}\: \chi_2(\eta) \subset \{\eta: \: |\eta| \geq 1 + \delta_1/2\},$$
$0 < \delta_1 <1$ being the constant in Sect. 6. Set
$$F_i(y) = \Bigl(\frac{k}{2\pi}\Bigr)^{N-1} G(y) \int e^{\ii k <y, \eta>}\chi_i(\eta)\hat{F}(k \eta) d\eta,\: i = 0,1,2.$$
To treat $F_0$ we will apply the results of Sects. 3-7. Consider the function
$$\psi(y; \eta) = <y, \eta>,\:y \in W,\: |\eta| < 1 - \delta_1/2.$$
We can construct a phase function $\varphi = \varphi(x; \eta)$ defined in ${\mathcal V}_j$ such that

$$(i)\:\: \varphi\bv_{{\rm supp}\: G} = \psi(y; \eta),\: y \in W,$$
$$(ii)\:\: \frac{\partial \varphi}{\partial \nu}(x; \eta)\bv_{{\mathcal V}_j \cap \gj} \geq \delta_2 > 0, \:y \in W,$$

$\hspace{3.1cm} (iii)$ the phase $\varphi(x; \eta)$ satisfies the condition $({\mathcal P})$ on $\gj$.\\

The local existence of $\varphi(x; \eta)$ satisfying the conditions (i)-(ii) has been discussed in \cite{kn:I2}, \cite{kn:I3}. To arrange (iii), we use a suitable continuation and we 
treat this problem in  Appendix B below.  Starting with the oscillatory data $m_0(y; \eta) = (2\pi)^{-N + 1} G(y) e^{\ii k \la y, \eta \ra},\: |\eta| \leq 1 - \delta_1/2$ and applying the 
argument of Sects. 6-7, we construct an approximative solution $V_0(x, s; k , \eta)$ which satisfies the condition $(S)$ in $(\ovo, \cd1)$ and such that
$$V_0(y,s;k, \eta) = m_0(y; \eta) + s^{-M} {\mathcal Q}_M(y, s; k, \eta),\: x \in \Gamma.$$
Moreover, for ${\mathcal Q}_{M}(x,s; k, \eta)$ we have the estimate (\ref{eq:7.9}) and it is clear that the constants $B_M$ and $\alpha$ in (\ref{eq:7.9}) can be chosen uniformly 
with respect to  $\eta, \: |\eta| \leq 1 - \delta_1/2.$  Define the operator 
$$U_0(s; k)F = \int V_0(x, s; k, \eta)\chi_0(\eta)\hat{F}(k\eta) k^{N-1}d\eta$$ 
with values in $C^{\infty}(\Omega)$ so that $U_0(s;k)F$  satisfies the condition $(S)$ in $(\ovo, \cd1)$ and
$$U_0(s; k)F\vert_{\Gamma} =  F_0+ s^{-M} \int {\mathcal Q}_{M}(x, s; k, \eta)\chi_0(\eta) \hat{F}(k\eta) k^{N-1} d\eta$$
$$= F_0 + L_0(s; k) F.$$
Therefore
$$\|L_0(s; k) F\|_{L^2(\Gamma)}^2 \leq C_0 \Bigl(\int_{|\eta| \leq 1 - \delta_1/2} 
k^{-M + (N-1)/2 + \alpha} |\hat{F}(k\eta)| k^{(N-1)/2} d\eta\Bigr)^2$$
$$\leq C_0 k^{-2M + N-1 + 2\alpha} \int_{|\eta| \leq 1 -\delta_1/2} d\eta \int_{\R^{N-1}} 
|\hat{F}(k \eta)|^2 k^{N-1} d\eta \leq C_1 k^{-2M + N-1 + 2\alpha} \|F\|_{L^2(\Gamma)}^2$$
with a constant $C_1 > 0$ depending only on $K.$ Moreover, for $s \in \cd1$ we obtain the estimate
\begin{equation} \label{eq:8.2}
\|U_0(s; k)F\|_{L^2(\Omega \cap \{|x| \leq R\})} \leq C_{0, R}k^{m_0} \|F\|_{L^2}.
\end{equation}
To prove this, it is sufficient to show that
\begin{equation} \label{eq:8.3}
\| V_0(x,s; k, \eta)\|_{L^2(\Omega \cap \{|x| \leq R\})} \leq C_{0, R}'k^{p_0},\: s \in \cd1\;,
\end{equation}
uniformly with respect to $|\eta| \leq 1- \delta_1/2.$ On the other hand,
$$V_0(x, s;k, \eta) = V^{(0)}(x, s; k, \eta) - \sum_{m = 1}^{M-1} V^{(m)}(x,s; k, \eta)s^{-m}$$
and
$$V^{(m)}(x, s; k, \eta) = \sum_{j_1, j_2,...,j_{m} = 1}^{\kappa_0}{\mathcal S}_{j_{m}}(s) w_{j_1, j_2,...,j_{m}}(x, s; k, \eta).$$
Here $w_{j_1, j_2,...,j_{m}}(x, s; k, \eta),\: x \in \Gamma_{\ell},$ are infinite series and the estimates of $\| V^{(m)}\|_{L^2(\Omega \cap \{|x| \leq R\})}$ follow from the 
estimates for the operators $H_h, H_g, S_0(s), P_e$ and the estimates for $\|w_{j_1, j_2,...,j_{m}}\|_{H^2(\Gamma_{m})}.$ According to the recurrence procedure in Sect. 7,
we deduce that
$$\|w_{j_1, j_2,...,j_{m}}\|_{H^2(\Gamma_{m})} \leq D_{\ell}|s|^{q(m)},\:s \in \cd1, m = 0,...M-1\;,$$
for some integers $q(m)$, and we get (\ref{eq:8.3}) with $p_0 = \sup_{m} q(m).$

\def\emu{{\mathcal E}_{+, \nu}}
\def\uj{{\mathcal U}_j}

To deal with $F_1(y)$, introduce $\xi(y, \eta)\in \S^{N-1}$ such that 
$$\xi(y, \eta) - \la \nu(y),\xi(y, \eta)\ra = \eta, \:(y, \eta) \in \Xi = {\rm supp}\: G \times \{\eta:\: -\frac{2}{3} \delta_1 \leq |\eta| -1\leq \frac{2}{3}\delta_1\}$$
and consider 
$\zeta(y, \eta) = \xi(y, \eta) - 2 \la \nu(y), \xi(y, \eta) \ra \nu(y) \in \S^{N-1}.$ Our choice of $\delta_1$ in Sect. 6 and Lemma 1 show that at least one of the rays
$\{y + t\xi(y, \eta):\: t \geq 0\},\: \{y + t\zeta(y, \eta):\:t \leq 0\}$ does not meet a $d_0$-neighborhood of $\bigcup_{\ell \not= j} K_{\ell}.$ For every fixed $(y_0, \eta_0) \in \Xi$ we 
have the above property for at least one of the rays related to $\xi(y_0,\eta_0)$ and $\zeta(y_0, \eta_0)$ and the same is true for $(y, \eta)$ sufficiently close to $(y_0, \eta_0).$ 
Consider a microlocal partition of unity  $\sum_{\mu = 1}^{M_1} \psi_{\mu}(y)\Xi_{\mu}(\eta) = 1$ on $\Xi$ so that 
${\supp}\: \Xi_{\mu} \subset \{\eta:\: -\delta_1 \leq |\eta| - 1 \leq \delta_1\},$
while for $(y, \eta) \in {\rm supp}\: \psi_{\mu}\Xi_{\mu}$  we have the property of the rays mentioned above. We fix $\mu$ and assume first that the outgoing rays 
$\{y + t\xi(y, \eta):\:t \geq 0\}$,  $(y, \eta) \in {\rm supp}\: \psi_{\mu}\Xi_{\mu}$ do not meet a neighborhood of $\bigcup_{\ell \not= j} K_{\ell}.$
Consider boundary data 
$$\widetilde{m}_{\mu}(y; k, \eta) = (2\pi)^{-N+1} G(y)\psi_{\mu}(y)e ^{\ii k \la y, \eta \ra},\:\eta \in {\rm supp}\: \Xi_{\mu}.$$
Following Proposition 4.7 in \cite{kn:I3} (see also Proposition 7.5 in \cite{kn:I1}), for every $M \geq 1$ there exists a function
$Z_{\mu, M}(x, s; k, \eta)$ which satisfies condition (S) in $(\Omega_j, \cd1)$ and
\begin{eqnarray} \label{eq:8.4}
\|Z_{\mu, M}(., s; k, \eta)\|_{C^p(\Omega_j \cap \{|x| \leq R\})} \leq C_{R,p} k^p,\: \forall p \in \N,\\ \nonumber
Z_{\mu, M}(y, s; k, \eta) = \widetilde{m}_{\mu}(y; k, \eta) + r^{-M}D_{\mu, M}(y, s; k, \eta),\: y \in \Gamma\;,
\end{eqnarray} 
with
$$\|D_{\mu, M}(., s; k, \eta)\|_{\Gamma, p} \leq C_p k^p,\: \forall p \in \N.$$
The constants in the above estimates are uniform with respect to $\eta$ and $\mu$ and they depend only on the geometry of $K.$

The construction of $Z_{\mu}$ in \cite{kn:I1} is long and technical. We sketch below the main points. The starting point is to introduce oscillatory boundary data
$$(2\pi)^{-N+1} G(y)\psi_{\mu}(y)h(t)e ^{\ii k (\la y, \eta \ra -t)},\:\eta \in {\rm supp}\: \Xi_{\mu}\;,$$
depending on $y$ and $t$ with $h \in C_0^{\infty}(\R^+),$ supp $h \subset (T,T + 1),\: T > 1$ and to construct an asymptotic solution $w_{\mu}(x, t; k, \eta)$ of the wave 
equation $(\partial_t^2 - \Delta_x)u = 0$ for $t \geq 0$ with supp $w_{\mu}(x, t; ., .) \subset \{(x,t):\: t \geq 0\}$ and {\it big parameter} $k$. We omit in the notation here and 
below the dependence on $M$.  In the glancing region we have two phase functions $\varphi_{\pm} =\theta(y, \eta) \pm \frac{2}{3}\rho^{3/2}(y, \eta)$ (see  \cite{kn:I1}, 
\cite{kn:G}, \cite{kn:SV}) and $\varphi_{\pm}$ are constructed so that their traces on ${\rm supp}\: G \cap \gj$ coincide with $\la y, \eta \ra$. The outgoing rays are propagating 
with directions $\nabla \varphi_{+},$ while the incoming rays are propagating with directions $\nabla \varphi_{-}$. The proofs in \cite{kn:I1} and \cite{kn:I3} work assuming $N$ 
odd and one considers the Laplace transform
$$\hat{w}_{\mu}(x, s; k, \eta) = \int_{-\infty}^{\infty} e^{-st} w_{\mu}(x, t; k,\eta) dt,\: s \in \cd1.$$
The assumption $N$ odd is used only by applying the strong Huygens principle to guarantee that for every fixed $x \in \Omega_j$ the support of $w_{\mu}$ with respect to 
$t$ is compact, hence the integral is convergent. For $N$ even we apply the finite speed of propagations and the fact that the supports of the solutions of the transport equations 
are propagating along the rays $\{y + t\nabla \varphi_{+}(y, \eta):\: t \geq 0\}$ to show that for $|x| \leq \rho_0$ the solution $w_{\mu}(x, t; k, \eta)$ vanishes for $t$ large. This 
justifies the existence of $\hat{w}_{\mu}(x, s; k, \eta)$ for $|x| \leq \rho_0$. Next, using the notation of Sect. 6, consider
\begin{equation} \label{eq:8.5}
Z_{\mu}(x, s; k, \eta) =\frac{1}{\hat{h}(s + \ii k)}\Bigl[\Phi \hat{w}_{\mu} - S_0(s)\Bigl(\Phi(\Delta_x - s^2)\hat{w}_{\mu} + [\Delta, \Phi]\hat{w}_{\mu}\Bigr)\Bigr],
\end{equation}
where $h$ is chosen so that $\hat{h}(s + \ii k) \not= 0$ for $|s + \ii k| \leq |\sigma_0| + c.$
Now let $\mu$ be such that the rays $\{y +t \zeta(y, \eta):\: t \leq 0\}, \:(y, \eta) \in {\rm supp}\: \psi_{\mu}\Xi_{\mu},$ do not meet a neighborhood of $\bigcup_{\ell \not= j} K_{\ell}.$ 
In this case we repeat the procedure in Section 7 in \cite{kn:I1} and Section 4 in \cite{kn:I3} to construct an asymptotic solution $w_{\mu}(x, t; k, \eta)$ of the wave equation for 
$t \leq 0$ with supp $w_{\mu}(x, t; ., .) \subset \{(x, t);\: t \leq 0\}$ starting with oscillatory boundary data 
$$(2\pi)^{-N+1} G(y)\psi_{\mu}(y)h(-t)e ^{-\ii k (-\la y, \eta \ra - t)},\:\eta \in {\rm supp}\: \Xi_{\mu}.$$
We express $\la y, \eta \ra$ by the trace of the phase function $\varphi_{-}\bv_{\gj}$ related to the incoming directions and we consider for $|x| \leq \rho_0$ the Laplace transform
$$\hat{w}_{\mu}(x, s; k, \eta) = \int_{-\infty}^{\infty} e^{st} w_{\mu}(x, t; k,\eta) dt,\: s \in \cd1.$$
Next, we define $Z_{\mu}(x, s; k, \eta)$ by (\ref{eq:8.5}) using the outgoing parametrix $S_0(s)$ and deduce the estimates (\ref{eq:8.4}). Finally, we introduce
$$U_1(s; k)F = \sum_{\mu = 1}^{M_1} \int Z_{\mu}(x, s; k, \eta) \Xi_{\mu}(\eta)\chi_1(\eta) \hat{F}(k\eta)k^{N-1} d\eta$$
and conclude that $U_1(s;k)F$ is analytic for $s \in \cd1$ and satisfies 
$$\begin{cases} (\Delta_x - s^2) U_1(s; k)F = 0,\: x \in \Omega_j,\cr
U_1(s; k)F\bv_{\Gamma} = F_1 + L_1(s; k)F.\end{cases}$$
As above, exploiting the estimates (\ref{eq:8.4}), we obtain 
$$\|L_1(s; k)F\|_{L^2({\Gamma})} \leq C_M k^{-M}\|F\|_{L^2(\Gamma)},\: s\in \cd1,$$
and
\begin{equation} \label{eq:8.6}
\| U_1(s; k)F\|_{L^2(\ovo\: \cap \{|x| \leq R\})}\leq C_{1,R} k^{\frac{N-1}{2}}\|F\|_{L^2}.
\end{equation}
\def\tp{\tilde{\varphi}}

Now we pass to the analysis of the term $F_2$ in the elliptic region. Let $\tilde{{\mathcal U}}_j$ be a small neighborhood of $K_j$ and let $\uj = \tilde{{\mathcal U}}_j \setminus K_j.$ Following \cite{kn:SV}, Appendix A. 4, we construct a parametrix $H_e: \tc({\rm supp}\: G) \longrightarrow \tc(\uj)$ as a Fourier integral operator with complex phase function $\tilde{\varphi}(x, \eta)$ and big parameter $k$ having the form
$$(H_e u)(x, s) = \Bigl(\frac{s}{2\pi}\Bigr)^{N-1} \int e^{\ii k (\tp (x, \eta) - \la y, \eta \ra)} \tilde{a}(x, \eta, k) u(y) dy d\eta,$$
so that
$$\begin{cases} (\Delta_s - s^2) H_e u = K_e u,\: x \in \uj,\cr
H_e u\bv_{\gj} = Op(G\chi_2) u,\end{cases},$$
where 
$$Op(G\chi_2) u = \Bigl(\frac{k}{2\pi}\Bigr)^{N -1}\int e^{\ii k \la x - y, \eta \ra}G(x) \chi_2(\eta) u(y) dy d\eta.$$
The last operator is defined for $u \in C^{\infty}(\gj)$ but it can be prolonged to $F \in L^2(\gj)$ since the symbol
$\chi_2(\eta) \in S^{(0,0)}_{0,0}$ (see Proposition A.I. 6 in \cite{kn:G}).

Assume that locally the boundary $\gj$ is given by the equation $x_N = 0$ and let locally $\uj \subset \{x_N \geq 0\}$. To satisfy the equation $(\Delta_x - s^2)H_e u = 0$ modulo negligible terms, we must choose $\tp$ so that
\begin{equation} \label{eq:8.7}
|\nabla \tp|^2 = -\Bigl(\frac{s}{ k}\Bigr)^2 = \gamma^2,\: \tp\bv_{\gj} = \la x, \eta \ra.
\end{equation}
For $|s + \ii k| \leq |\sigma_0| + c$ we see that $\gamma = 1 + {\mathcal O}(k^{-1})$ is a complex parameter close to 1 and we may repeat the argument in Appendix A. 4 in \cite{kn:SV} and Appendix A. II. 4 in \cite{kn:G} to construct $\tp$ with the properties
$$\Im\: \tp(x, \eta) \geq c_0 x_N (1 + |\eta|),\: c_0 > 0,\: |\Re \:\tilde{\varphi}(x, \eta)| \leq c_0'(1 + |\eta|).$$
The phase $\tp$ satisfies the eikonal equation modulo ${\mathcal O}(x_N^{\infty})$, the amplitudes satisfy the corresponding transport equations modulo ${\mathcal O}(x_N^{\infty})$ and $\tilde{a}(x, \eta, k) \in S_{0,0}^{0,0}.$ Notice that the sign of $\Im \: \tp(x, \eta)$ is related to the choice $k > 0.$ We have
$$\Re \Bigl(ik(\tp(x, \eta) - \la y, \eta \ra )\Bigr) = -k \Im\: \tp(x, \eta) \leq -c_0 k x_N (1 + |\eta|)$$
and the integral $H_e F$ is convergent for $x_N > 0$ and $F \in L^2(Y)$. Moreover, we have
$$\sup_{x_N \geq 0} x_N^m e^{-c_0 x_N(1 + |\eta|)} \leq c_m(1 +|\eta|)^{-m} k^{-m},\:\forall m \in \N$$
and this implies that the kernel of $K_e$ is in $\tc(\uj \times {\rm supp}\: G)$ and we obtain $K_e = {\mathcal O}(|k|^{-\infty})$ uniformly with respect to $x_N \in [0, \epsilon].$ 

Next, let $\Psi(x) \in C_0^{\infty}({\mathcal U}_j)$ be a cut-off function such that $\Psi(x) = 1$ in a small neighborhood of $K_j$. Define
$$U_2(s; k) F = \Bigr[\Psi H_e - S_0(s)\Bigl(\Psi K_e + [\Delta, \Psi] H_e \Bigr)\Bigr] F.$$ 
Then $U_2(s; k) F$ satisfies 
$$\begin{cases} (\Delta_x - s^2) U_2(s; k)F = 0,\: x \in \ovo,\: s \in \cd1,\cr
U_2(s; k)F\bv_{\Gamma} = F_2 + L_2(s; k)F,\end{cases}$$
but $U_2(s; k)F$ is not analytic with respect to $s$ which will not be important for the proof of Theorem 2 below.
On the other hand, the trace on $\Gamma$ of $S_0(s) [\Delta, \Psi] H_e F$ is negligible and the same is true for the trace of $S_0(s) \Psi K_e F.$ Thus,
$\|L_2(s; k)F\|_{L^2(\Gamma)} \leq C_M k^{-M} \|F\|_{L^2(\Gamma)},\: \forall M \in \N.$
Moreover, we have the estimate
\begin{equation} \label{eq:8.8}
\|U_2(s; k)F\|_{L^2(\Omega_j \cap \{|x| \leq R\})} \leq  C_{2,R}\|F\|_{L^2(\Gamma)}
\end{equation}
which is a consequence of $L^2$ estimates of $\Psi H_e F$ and $[\Delta, \Psi] H_e F$. In fact, the estimate of $\|[\Delta, \Psi]H_e F\|_{L^2}$ is easy since $\Psi = 1$ in a 
neighborhood of $\Omega_j$ and the kernel of $[\Delta, \Psi]H_e$ is in $\tc(\uj \times {\rm supp}\: G).$ To estimate $\|\Psi H_e F\|_{L^2}$, we observe that for small $x_N \geq 0,$ 
$H_e$ is a Fourier integral operator with non-degenerate phase function of positive type $\phi(x, y, \eta) = \tp(x, \eta) - \la y, \eta \ra$ (see Definition 25.4.3 in \cite{kn:H}). Thus, we 
can estimate 
$$\|(H_e F)(x_N,.,s; k)\|_{L^2(\uj \cap \{x_N = z\})} \leq B\|F\|_{L^2(\Gamma)}$$ 
uniformly with respect to $z \in [0, \epsilon]$ (see Theorem 25.5.6 in \cite{kn:H}) and this leads to (\ref{eq:8.8}).
Finally, introduce
$$L_{Y}(s; k) F =U_0(s; k)F + U_1(s; k)F + U_2(s; k)F$$
and conclude that
$$L_Y(s; k)F\vert_{\Gamma} = F + \sum_{i=0}^2L_i(s; k)F = F + Q_Y(s; k)F$$
with
$$\|Q_Y(s; k)F\|_{L^2(\Gamma)} \leq B_{Y} k^{-M + (N-1)/2 + \alpha} \|F\|_{L^2(\Omega)}.$$
By using a partition of unity on $\Gamma$, we define an operator
$$L(s; k): L^2(\Gamma) \ni f \longrightarrow L(s ; k)f \in C^{\infty}(\ovo)$$
and deduce that $L(s; k)f$ satisfies
$$\begin{cases} (\Delta_x - s^2) L(s, k)f = 0,\: x \in \ovo,\cr
L(s, k)f \in L^2(\Omega),\: {\rm for}\: \Re \: (s) > 0,\cr
L(s, k)F\bv_{\Gamma} = f + Q(s; k)f\end{cases}$$
with
$$\|Q(s; k)f\|_{L^2(\Gamma)} \leq B k^{-M+ (N-1)/2 + \alpha} \|f\|_{L^2(\Gamma)}.$$
Choosing $k_1$ sufficiently large, the operator $I + Q(s; k): L^2(\Gamma) \to L^2(\Gamma)$ is invertible for $s \in \cd1$ and $k \geq k_1$. We define
$$R(s, k) f = L(s; k)(I + Q(s; k))^{-1} f: L^2(\Gamma) \longrightarrow C^{\infty}(\ovo)$$
and it is clear that $R(s, k)f$ for $s \in \cd1$ satisfies (\ref{eq:8.1}).\\

{\it Proof of Theorem $2$.} Given $g \in L^2(\ovo)$ and $\chi \in C_0^{\infty}(\ovo)$ with supp $\chi \subset \{|x| \leq \rho\},\: \rho \geq \rho_0 $, by (\ref{eq:6.3}) we obtain 
$S_0(s)(\chi g) \in H^1(|x| \leq \rho)$ and this yields $[S_0(s) (\chi g)]\bv_{\Gamma} \in H^{1/2} (\Gamma).$ Setting $s = \ii z$, consider for $\Im\: (z) < 0$
\begin{equation} \label{eq:8.9}
 v = S_0(\ii z) (\chi g) - R(\ii z, k)\Bigl([S_0( \ii z)(\chi g)]\bv_{\Gamma}\Bigr).
\end{equation}
Then for the cut-off resolvent $R_{\chi}(z)$ introduced in Sect. 1 we get 
$$R_{\chi}(z) (\chi g) = \chi v,\: \Im \:(z) < 0.$$
The operators $\chi S_0(\ii z)\chi$ and $R_{\chi}(z)$  admit respectively analytic and meromorphic continuation from 
$\Im\: (z) < 0$ to $\{ z  \in \C: \Im\: (z) \leq -\sigma_1,\: \Re\: (z) < -J_1\}$, 
where $-J_1 = \min\{-J, |\sigma_0| + c - k_1\}.$ Thus, $\chi R(\ii z, k)\Bigl([S_0(\ii z)(\chi g)]\bv_{\Gamma}\Bigr)$ is also meromorphic in this domain and to show that it is 
analytic for  $\ii z \in \cd1$ it suffices to prove that this operator is bounded. For $\ii z \in \cd1$ this follows from the estimates (\ref{eq:8.2}), (\ref{eq:8.6}), (\ref{eq:8.8}) above and 
we obtain a polynomial bound for $\|\chi R(\ii z, k)\|_{L^2(\Gamma)\to L^2(\ovo)}$. Consequently, $R_{\chi}(z)$ admits an analytic continuation and we get (\ref{eq:1.7}) for 
$\Re \: (z) < - J_1 < 0.$ Next to cover  the case $\Re \: (z) > J_1 > 0$, we can use the fact that the poles of $R_{\chi}(z)$ are symmetric with respect to $\ii \R^{+}$ or repeat the 
argument with $k << 0$. \endofproof

\bs

To obtain Corollary 1 we establish the estimate
$$\|R_{\chi}(z)\|_{H^L(\ovo) \to L^2(\ovo)} \leq C(1 + |z|)^{m - L},\: z \in {\mathcal S,}$$
where $m \in \N$ is the integer in \eqref{eq:1.7} and $L \in \N,\: L > m$. The proof goes repeating that in the non-trapping case (see Theorem 1 in \cite{kn:TZ}) and 
we omit the details. \endofproof

\section{Appendix A : Stable and instable manifolds for open billiards}
\renewcommand{\theequation}{\arabic{section}.\arabic{equation}}
\setcounter{equation}{0}

Let $z_0 = (x_0,u_0)\in S^*(\Omega)$. For convenience we will assume that $x_0 \notin K$. Assume that
the {\it backward trajectory} $\gamma_-(z_0)$ determined by $z_0$ is bounded, and let $\eta  \in \san$
be its itinerary.

Given $x\in \R^N$ and $\epsilon > 0$, by $B(x,\epsilon)$ we denote
the {\it open ball} with center $x$ and radius $\epsilon$ in $\R^N$.

In this section we use some tools from \cite{kn:I3} to construct the {\it local unstable 
manifold}\footnote{Notice that   $\wuloc (z_0)$ and $W^u_\ep(z_0)$ (see Appendix C) coincide in a neighborhood of $z_0$.} 
$\wuloc (z_0)$ of  $z_0$ in $S^*(\Omega)$ and show that it is Lipschitz in $z_0$ (and $\eta$). 
In a similar way one deals with local stable manifolds. 

Notice that if the boundary $\Gamma$ of $K$ is only $C^k$ ($k \geq 2$) the $C^\infty$ smoothness below should 
be replaced by $C^k$.

\begin{prop} There exists a constant $\epsilon_0 > 0$ such that  for any $z_0 = (x_0,u_0)\in S_{\delta_0}^*(\Omega \cap B_0)$ 
whose backward trajectory $\gamma_-(z_0)$ has an infinite number of 
reflection points $X_j = X_j(z_0)$ ($j\leq 0$) and  $\eta\in \san$ is its itinerary, the following hold:

(a) There exists  a smooth ($C^\infty$) phase function $\psi = \psi_\eta$ 
satisfying part {\rm (i)} of the condition {\rm (${\mathcal P}$)} on $\uu = B(x_0, \epsilon_0)\cap \Omega$  such that
$\psi(x_0) = 0$, $u_0 = \nabla \psi(x_0)$, and such that for any $x\in C_\psi(x_0)\cap \uu^+(\psi)$
the billiard trajectory $\gamma_-(x, \nabla\psi(x))$ has an itinerary $\eta$ and therefore 
$d(\phi_t(x, \nabla \psi(x)), \phi_t(z_0) ) \to 0$ as $t \to -\infty\;.$ That is, 
$$\wuloc (z_0) = \{ (x, \nabla\psi(x)) : x\in C_\psi(x_0)\cap \uu^+(\psi) \}$$
is the local unstable manifold of $z_0$. Moreover, for any $p \geq  1$ 
there exists a global constant $\Con_p > 0$ (independent of $z_0$ and $\eta$) such that
\be
\|\nabla \psi_\eta\|_{(p)} (\uu) \leq \Con_p\;.
\ee

\medskip

(b)  If $(y,v)\in S^*(\Omega \cap B_0)$ is such that $y\in C_\psi(x_0)$ and $\gamma_-(y,v)$
has the same itinerary $\eta$, then $v = \nabla \psi(y)$, i.e. $(y,v) \in \wuloc(z_0)$.

\medskip

(c) There exist  a constant $\alpha \in (0,1)$ depending only on the 
obstacle $K$ and for every $p \geq 1$ a  constant $\Con_p  > 0$ such that for any integer $r\geq 1$ and any 
$\zeta,\eta \in \san$ with $\zeta_j = \eta_j$ for $-r \leq j\leq 0$, we have
$\|\nabla \psi_\eta - \nabla \psi_\zeta \|_p(V) \leq \Con_p\, \alpha^r$, where
$V = \uu(\psi_\eta) \cap \uu(\psi_\zeta)$.
\end{prop}

{\it Proof.} (a) Take $\epsilon_0 > 0$ so small that whenever $(x,u) \in S^*_{\delta_0/2}(\Omega \cap B_0)$ and $(y,v) \in S^*(\Omega)$ is
such that $\| x-y\| < \epsilon_0$ and $\|u-v\| < \epsilon_0$ we have $(y,v) \in S^*_{\delta_0}(\Omega)$. Then define $\uu$ as in part (a) above.

Set,$d_{-m} = \|X_{-m+1} - X_{-m}\|$ and  $u_{-m} = \frac{X_{-m+1} - X_{-m}}{\| X_{-m+1} - X_{-m}\|} \in \sn$ ($m \geq 1$). 
Given any integer $m \geq 1$, consider the linear phase function $\psim = \psime$  in $\Omega$ such that $\nabla \psim \equiv u_{-m}$
and $\psim(X_{-m}) = -(d_{-m} + d_{-m+1} + \ldots + d_{-1})$. Then define
$$\psimm = \psimme = \Phi_{\eta_{-1}}^{\eta_{0}} \circ \Phi_{\eta_{-2}}^{\eta_{-1}} \circ \ldots
\Phi_{\eta_{-m+1}}^{\eta_{-m+2}} \circ \Phi_{\eta_{-m}}^{\eta_{-m+1}}\, (\psim) \;.$$
Clearly $\psimm$ is a smooth phase function defined everywhere on $\uu$ (in fact, on a much larger subset of $\Omega$) with
$\psimm(X_0) = 0$. Moreover, it follows from Proposition 2 in Sect. 2 above that
\be
\| \nabla \psim_m - \nabla \psi^{(m+1)}_{m+1} \|_p(\uu) 
\leq \Con_p\, \alpha^m  \quad, \quad m \geq 1\;
\ee
for some global constant $\Con_p > 0$ depending only on $K$ and $p$. Here we use the fact  that
$\| \nabla \psim - \nabla \psi^{(m+1)} \|_{(p)} \leq \Con$, due to the special choice of the phase functions $\psim$ and $ \psi^{(m+1)}$.
Since $\psim_m (X_0) = \psi^{(m+1)}_{m+1}(X_0) = 0$, it  now follows that there exists a constant $\Con_p > 0$ such that
$\|\psim_m (x) - \psi^{(m+1)}_{m+1} (x)\| \leq \Con_p \, \alpha^m$ for  $x \in \uu\cap B_0\;.$ This implies that for every $x \in \uu$ there exists
$\psi(x) = \lim_{m\to\infty} \psim_m (x)$. Now (9.2) shows that $\psi$ is $C^{\infty}$-smooth in $\uu$ and
\be
\| \nabla \psim_m  - \nabla \psi \|_p(\uu) \leq \Con_p\, \alpha^m \quad, \quad m \geq 1\;.
\ee
In particular, $\| \nabla \psi\| \equiv 1$ in $\uu$.  Extending $\psi$  in a trivial way along straight line rays,
we get a phase function $\psi$ satisfying part (i) of the condition (${\mathcal P}$) in $\uu$. 

Let us now show that $W = \{ (x, \nabla\psi(x)) : x\in C_\psi(x_0)\cap \uu^+(\psi) \}$ is the local unstable
manifold of $z_0$. Given $x\in C_\psi(x_0)\cap \uu^+(\psi)$ sufficiently close to $x_0$ and
an arbitrary integer $r \geq 0$, consider the points
$X^{-r}(x,\psim_m) \in \dk_{\eta_{-r}}$ for $m \geq r$. By Proposition 1 in Sect. 2 above,
there exist global constants $\Con > 0$ and $\alpha \in (0,1)$ such that
$\|X^{-r}(x,\psim_m) - X^{-r}(x,  \psi^{(m')}_{m'})\| \leq \Con \, \alpha^{m-r}$ for $m' \geq m > r\;.$
Thus, there exists $X^{-r} = \lim_{m\to\infty} X^{-r}(x,\psim_m) \in \dk_{\eta_{-r}}$ and
\be
\|X^{-r}(x,\psim_m) - X^{-r}\| \leq \Con \, \alpha^{m-r} \quad, \quad m > r\;.
\ee
It is now easy to see that $\{ X^{-j}\}_{j=0}^\infty$ are the successive reflection points of a billiard trajectory in $\Omega$ and this is the trajectory 
$\gamma_-(x, \nabla \psi)$. The backward itinerary of the latter is obviously $\eta$. Moreover, (9.3) implies
$d(\phi_t(x, \nabla \psi(x)), \phi_t(z_0) ) \to 0$ as $t \to -\infty$, so $(x, \nabla \psi(x)) \in \wuloc(z_0)$. 

Finally, by (2.1), $\|\psim_m\|_{(p)}(\uu) \leq \Con_p \, \|\psim \|_{(p)} \leq \Con_p\;,$ and combining this with (9.3) gives (9.1).

\medskip

(b) Let $(y,v)\in S^*(\Omega)$ be such that $y\in C_\psi(x_0)$ and $\gamma_-(y,v)$
has the same itinerary $\eta$. Define the phase functions $\phimm$ and $\phim$ as in part (a)
replacing the point $z_0 = (x_0,u_0)$ by $z = (y,v)$, and let $\varphi(x) = \lim_{m\to\infty} \phimm(x)$.
Then by part (a), we have $\wuloc (z) = \{ (x, \nabla\psi(x)) : x\in C_\varphi(y)\cap \uu^+(\phi) \}$.

On the other hand, it follows from Proposition 2 that there exist constants $\Con > 0$ and $\alpha \in (0,1)$ such that
$\| \nabla \psimm - \nabla \phimm\| \leq \Con\, \alpha^m$ for all $m \geq 0$, which implies $\varphi = \psi$.
Thus, $v = \nabla \varphi(y) = \nabla \psi(y) \in \wuloc(z_0)$.

\medskip

(c) Choose the constants $\alpha \in (0,1)$ and $\Con_p > 0$ ($p = 1,\ldots,k$) as in part (a). Let
$\zeta,\eta \in \san$ be such that $\zeta_j = \eta_j$ for all $-r \leq j\leq 0$ for some $r \geq 1$.
Construct the phase functions $\psime_m$ and $\psi^{(m,\zeta)}_m$ ($m \geq 1$) as in part (a); then
$\psi_\eta  = \lim_{m\to \infty} \psime_m$, $\psi_\zeta = \lim_{m\to \infty} \psi^{(m,\zeta)}_m$.
It follows from Proposition 2 that  
$\| \nabla \psi^{(r,\eta)} - \nabla \psi^{(r,\zeta)}\| \leq \Con_p\, \alpha^r$. 
Combining this with (9.3) with $m = r$ for $\eta$ and then with $\eta$ replaced by $\zeta$, one gets
$$ \| \nabla \psi_\eta - \nabla \psi_\zeta\| \leq 
\| \nabla \psi_\eta - \nabla \psi^{(r,\eta)}\|  + \| \nabla \psi^{(r,\eta)} - \nabla \psi^{(r,\zeta)}\| 
+\| \nabla \psi^{(r,\zeta)} - \nabla \psi_\zeta\|  \leq \Con_p\, \alpha^r\;.$$
This proves the assertion.
\endofproof

\def\e{\varepsilon}
\def\CC{{\mathcal C}}
\def\ra{\rangle}
\def\la{\langle}
\def\eps{\epsilon}
\def\esp{\vspace{8pt}}
\def\eps{\epsilon}
\def\l2{L^2}
\def\Cp{\C_{+}}
\def\ii{{\bf i}}
\def\vx{|\varphi_y|^2}

\section{Appendix B : Construction of a phase function satisfying the condition (${\mathcal P}$)}
\renewcommand{\theequation}{\arabic{section}.\arabic{equation}}
\setcounter{equation}{0}

Consider a local representation $x_N = h(y)$ of the
boundary $\Gamma_j$ with $y = (y_1,...,y_{N-1}) \in W \subset \R^{N-1}.$ We wish to construct a phase function $\varphi(x; \eta)$ such that
$$ \varphi(y, h(y); \eta) = \la y, \eta \ra, \: (y, h(y)) \in U,\: \eta = (\eta_1,...,\eta_{N-1}),$$
$U$ being a small neighborhood of a fixed point $x_0 \in \Gamma_j$ so that $\varphi(x ; \eta)$ satisfies the conditions $(i)-(iii)$ of Sect. 8. 
Assume that $|\eta| \leq 1 - \mu$, where $0 < \mu < 1$. It is convenient to consider a little more general problem with boundary data given 
by a smooth function $\chi(y)$ such that $|\nabla_{y} \chi(y)| \leq 1 -\mu$ for $y \in \ W.$ We will construct a phase function $\varphi(x)$ such that
\begin{equation} \label{eq:9.5}
\varphi(y, h(y)) = \chi(y), \: y \in W\;,
\end{equation}
omitting the dependence on $\eta$ in the notation.
From the boundary condition (\ref{eq:9.5}) we determine the derivatives of $\varphi$ on the boundary $\Gamma_j.$ 
Set $\varphi_y = (\varphi_{y_1},...,\varphi_{y_{N-1}}),\: 
h_{y} = (h_{y_1},...,h_{y_{N-1}}),\: \chi_y = (\chi_{y_1},...,\chi_{y_{N-1}})$. 
We have $\varphi_y + \varphi_{x_N} h_{y} = \chi_{y}$, so setting $\varphi_{x_N} = \sqrt{1 - |\varphi_{y}|^2}$ and solving the system
$$\varphi_y  +\sqrt{1 - |\varphi_y|^2} h_y = \chi_y$$
we get
$$(1 - \vx) |h_y|^2 = |\chi_y|^2 + |\varphi_y|^2 - 2 \la \chi_y, \varphi_y \ra.$$
On the other hand,
$$2 \la \chi_y, \varphi_y \ra + 2 \sqrt{1 -\vx} \la h_y, \chi_y \ra = 2 |\chi_y|^2 \;, $$
which gives
$$(1 + |h_y|^2) (1 - \vx) - 2 \la h_y, \chi_y \ra \sqrt{1 - \vx} + |\chi_y|^2 - 1 = 0.$$
Consequently, for $\varphi_{x_N} = \sqrt{ 1 - \vx}$ we obtain
$$\varphi_{x_N}(y, h(y)) = \frac{1}{1 + |h_y|^2} \Bigl( \la h_y, \chi_y \ra + 
\sqrt{\la h_y, \chi_y \ra^2 + (1 - |\chi_y|^2)(1 + |h_y|^2)}\Bigr).$$
Now it is easy to see that we have the condition
\begin{equation} \label{eq:10.2}
\la \nabla \varphi(x), \nu(x) \ra \geq \delta_0 > 0,\: x = (y, h(y)) \in U.
\end{equation}
In fact in local coordinates $x = (y, h(y))$ the outward normal to $\Gamma_j$ is given by
$$\nu(x) = \frac{1}{\sqrt{1 + |h_y|^2}}(-h_y, 1)$$
and we deduce
$$\la \nabla \varphi(x), \nu(x) \ra   = \frac{1}{\sqrt{1 + |h_y|^2}}\Bigl[(1 + |h_y|^2) \varphi_{x_N} - 
\la h_y, \chi_y \ra\Bigr] \geq \sqrt{ 1 - |\chi_y|^2}  \geq \sqrt{2\mu - \mu^2} > 0.$$
By using (10.2) and a standard argument, we can solve locally the eikonal equation $|\nabla \varphi(x)| = 1$ with initial data
$$\varphi(y, h(y)) = \chi(y),$$
$$ \nabla_x\varphi (y, h(y)) = \Bigl(\varphi_y(y, h(y)), \varphi_{x_N}(y, h(y))\Bigr),\: (y, h(y)) \in U.$$ 
This argument works for local boundary condition $\chi(y) = \la y, \eta \ra,\: |\eta| \leq 1 - \delta_1/2,$ and we obtain a phase function 
$\varphi(x; \eta),\: x = (y, h(y)),\:y  \in W.$ As in \cite{kn:I3},\cite{kn:B}, we show that the principal curvatures of the wave front
$${\mathcal G}_{\varphi}(z) = \{y \in \R^N:\: \varphi(y; \eta) = \varphi(z; \eta)\}$$
are strictly positive for every $z = (y, h(y)) \in U.$\\

In order to satisfy the condition $({\mathcal P})$ on $\gj$, we will construct a suitable continuation of $\varphi(x; \eta)$. 
For this purpose fix a point $x_0 = (y_0, h(y_0)) \in U.$ Without loss of generality, we can assume that $\varphi(x_0; \eta) = 0$. 
Consider a sphere $S_0$ passing through $x_0$ with center O in the interior of $K_j$ so that the unit outward normal $\nu_0$ of $S_0$ at $x_0$ coincides with  $\nabla \varphi(x_0; \eta).$

Choosing local coordinates $(\theta, z(\theta)),\: \theta \in W \subset \R^{N-1}$ on $S_0$, let $\Xi_0 = \{(\theta, z(\theta)):\: |\theta - \theta_0| \leq 2\epsilon\} \subset S_0$ be a small neighborhood of $x_0 = (\theta_0, z(\theta_0))$. Consider the trace $\Phi(\theta) = \varphi(\theta, z(\theta))$ of $\varphi$ on $\Xi_0.$ (We omit again the dependence on $\eta$ in the notation.) Since $\Phi(\theta_0) = 0$ and $\nabla_{\theta} \Phi(\theta_0) = 0$, we have
$$|\Phi(\theta)| \leq C_0 \epsilon^2,\: |\nabla_{\theta}\Phi(\theta)|\leq C_1 \epsilon, \: \theta \in \Xi_0.$$
 Choose a smooth cut-off function $\alpha(\theta),\: 0 \leq \alpha(\theta) \leq 1,$ such that 
$\alpha(\theta) = 1$ for $|\theta - \theta_0| \leq \epsilon/2,\: \alpha (\theta) = 0$ for $|\theta - \theta_0| \geq  \epsilon$ with $|\nabla_{\theta} \alpha| \leq C_2\epsilon^{-1}.$ 
Set $\chi(\theta) = \alpha(\theta) \Phi(\theta).$ Then for small $\epsilon > 0$ we have
$$|\nabla_{\theta} \chi(\theta)| \leq (C_0 C_2 + C_1) \epsilon < 1 - \mu < 1.$$
By the above procedure we construct a phase function $\Psi(x)$ so that $\Psi(\theta, z(\theta)) = \chi(\theta),\: |\theta - \theta_0| \leq 2\epsilon.$ For 
$\Xi' = \{(\theta, z(\theta)):\:\epsilon \leq |\theta - \theta_0| \leq 2\epsilon\} \subset \Xi_0$, it is easy to see that
$\nabla_x \Psi\bv_{\Xi'}$ coincides with the unit normal $\nu_0$ to $S_0$.  Thus if  $x = z + t\nu_0(z),\: t \geq 0$ with $z \in \Xi'$,  we have 
$\Psi(x) = t$ and for such $x$ the phase $\Psi(x)$ coincides with the 
phase function $\widetilde{\Psi}(x)$ defined {\bf globally} in a neighborhood of $S_0$ and having boundary data 
$\widetilde{\Psi}(x) = 0,\: \forall x \in S_0.$ Consequently, we may consider $\widetilde{\Psi}(x)$ as a continuation of $\Psi(x)$, so $\Psi(x)$ is defined 
globally outside a small neighborhood of the center $O$ of $S_0$ lying in the interior of $K_j.$ It is clear that $\Psi$ satisfies the condition 
$({\mathcal P})$ on $S_0$. On the other hand, for $\Xi_1 = \{(\theta, z(\theta)):\: |\theta - \theta_0| \leq \epsilon/2\}$ we have $\Psi\bv_{\Xi_1} = \varphi\bv_{\Xi_1}$ 
and locally in a neighborhood of $x_0$ the phases $\Psi(x)$ and $\varphi(x)$ coincide. Thus, we can consider $\Psi(x)$ as a continuation of $\varphi(x)$.

\def\nexto{\kern -0.54em}

\def\DF{{\rm {I\ \nexto F}}}
\def\DR{{\rm {I\ \nexto R}}}
\def\DN{{\rm {I\ \nexto N}}}
\def\DZ{{\rm {Z \kern -0.45em Z}}}
\def\DC{{\rm\hbox{C \kern-0.83em\raise0.08ex\hbox{\vrule height5.8pt width0.5pt} \kern0.13em}}}
\def\DQ{{\rm\hbox{Q \kern-0.92em\raise0.1ex\hbox{\vrule height5.8pt width0.5pt}\kern0.17em}}}
\def\T{{\mathbb T}}
\def\S{{\mathbb S}}
\def\sn{{\S}^{n-1}}
\def\sN{{\S}^{N-1}}

\def\ep{\epsilon}
\def\e{\emptyset}
\def\di{\displaystyle}
\def\sk{\smallskip}
\def\bs{\bigskip}
\def\ms{\medskip}
\def\dk{\partial K}
\def\kamin{\kappa_{\min}}
\def\kamax{\kappa_{\max}}

\def\saa{\Sigma_A^+}
\def\sa{\Sigma_A}
\def\san{\Sigma^-_A}

\def\Prf{\mbox{\footnotesize\rm Pr}}
\def\Pr{\mbox{\rm Pr}}

\def\be{\begin{equation}}
\def\ee{\end{equation}}
\def\beqn{\begin{eqnarray}}
\def\eeqn{\end{eqnarray}}
\def\beqn*{\begin{eqnarray*}}
\def\eeqn*{\end{eqnarray*}}
\def\endofproof{{\rule{6pt}{6pt}}}

\def\i{{\bf i}}
\def\iii{{\bf \i}}
\def\iu{\underline{i}}
\def\ju{\underline{j}}
\def\ku{\underline{k}}
\def\pu{\underline{p}}
\def\xx{{\bf x}}
\def\yy{{\bf y}}
\def\zz{{\bf z}}

\def\dist{\mbox{\rm dist}}
\def\diam{\mbox{\rm diam}}
\def\pr{\mbox{\rm pr}}
\def\supp{\mbox{\rm supp}}
\def\Arg{\mbox{\rm Arg}}
\def\In{\mbox{\rm Int}}
\def\Im{\mbox{\rm Im}}
\def\Int{\mbox{\rm Int}}
\def\span{\mbox{\rm span}}
\def\con{\mbox{\rm const}}
\def\Con{\mbox{\rm Const}}
\def\spec{\mbox{\rm spec}\,}
\def\Re{\mbox{\rm Re}}
\def\ecc{\mbox{\rm ecc}}
\def\var{\mbox{\rm var}}
\def\conf{\mbox{\footnotesize\rm const}}
\def\Conf{\mbox{\footnotesize\rm Const}}
\def\mt{\Lambda}

\def\ff{{\mathcal F}}
\def\kk{{\mathcal K}}
\def\cc{{\mathcal C}}
\def\mm{{\mathcal M}}
\def\nn{{\mathcal N}}
\def\ll{{\mathcal L}}
\def\uu{{\mathcal U}}
\def\ss{{\mathcal S}}
\def\oo{{\mathcal O}}
\def\rr{{\mathcal R}}
\def\pp{{\mathcal P}}
\def\gg{{\mathcal G}}
\def\vv{{\mathcal V}}
\def\aa{{\mathcal A}}
\def\jj{{\jmath }}
\def\H{{\mathcal H}}

\def\hy{\hat{y}}
\def\hla{\hat{\lambda}}
\def\hx{\hat{x}}
\def\hxi{\hat{\xi}}
\def\hu{\hat{u}}
\def\hf{\hat{f}}
\def\hg{\hat{g}}
\def\ho{\hat{\omega}}
\def\hD{\hat{\Delta}}
\def\hxi{\hat{\xi}}
\def\heta{\hat{\eta}}
\def\hQ{\hat{Q}}
\def\htau{\hat{\tau}}
\def\hpi{\hat{\pi}}
\def\hL{\widehat{L}}
\def\he{\hat{e}}
\def\hchi{\hat{\chi}}
\def\hju{\hat{\ju}}
\def\hh{\hat{h}}
\def\hw{\hat{w}}
\def\hde{\hat{\delta}}

\def\nv{\nabla \varphi}
\def\Od{\mbox{\rm Od}}

\def\h{h_{\mbox{\rm\footnotesize top}}}

\def\tc{\tilde{c}}
\def\tx{\tilde{x}}
\def\ty{\tilde{y}}
\def\tz{\tilde{z}}
\def\tu{\tilde{u}}
\def\tv{\tilde{v}}
\def\ta{\tilde{a}}
\def\td{\tilde{d}}
\def\tf{\tilde{f}}
\def\tg{\tilde{g}}
\def\tl{\tilde{\ell}}
\def\tr{\tilde{r}}
\def\tt{\tilde{t}}
\def\tw{\tilde{w}}
\def\tth{\tilde{\theta}}
\def\tka{\tilde{\kappa}}
\def\tla{\tilde{\lambda}}
\def\tmu{\tilde{\mu}}
\def\tga{\tilde{\gamma}}  
\def\trho{\tilde{\rho}}
\def\tvar{\tilde{\varphi}}
\def\tF{\tilde{F}}
\def\tU{\tilde{U}}
\def\tV{\tilde{V}}
\def\tS{\tilde{S}}
\def\tP{\widetilde{P}}
\def\tQ{\widetilde{Q}}
\def\tR{\widetilde{R}}
\def\ts{\tilde{\sigma}}
\def\tpi{\tilde{\pi}}
\def\tla{\tilde{\lambda}}
\def\tf{\tilde{f}}

\def\v{{\sf v}}
\def\piU{\pi^{(U)}}
\def\mtb{ \mt_{\dk}}

\def\sAA{\Sigma_{\aa}^+}
\def\sA{\Sigma_{\aa}}
\def\sAn{\Sigma_{\aa}^-}

\def\sAA{\Sigma_{\aa}^+}
\def\sA{\Sigma_{\aa}}
\def\sAn{\Sigma_{\aa}^-}

\section {Appendix C: Dolgopayt type estimates for open billiards}
\renewcommand{\theequation}{\arabic{section}.\arabic{equation}}
\setcounter{equation}{0}

Here we first state the assumptions about the billiard flow and the non-wandering set $\mt$ under which the
results in \cite{kn:St4} imply  the Dolgopyat type estimates (3.3). Following \cite{kn:PS2}, we then explain how 
to apply these in the situation described in  Sect. 6 above. Full details of the arguments can be found in  \cite{kn:PS2}.

For $x\in \Lambda$ and a sufficiently small $\epsilon > 0$ let 
$$W_{\ep}^s(x) = \{ y\in S^*(\Omega) : d (\phi_t(x),\phi_t(y)) \leq \epsilon \: \mbox{\rm for all }
\: t \geq 0 \; , \: d (\phi_t(x),\phi_t(y)) \to_{t\to \infty} 0\: \}\; ,$$
$$W^u_{\ep}(x) = \{ y\in S^*(\Omega) : d (\phi_t(x),\phi_t(y)) \leq \epsilon \: \mbox{\rm for all }
\: t \leq 0 \; , \: d (\phi_t(x),\phi_t(y)) \to_{t\to -\infty} 0\: \}$$
be the (strong) {\it stable} and {\it unstable manifolds} of size $\epsilon$. Then
$E^u(x) = T_x W^u_{\ep}(x)$ and $E^s(x) = T_xW^s_{\ep}(x)$.

The following {\it pinching condition}\footnote{It appears that in the proof of the estimates (3.3), in the case of 
open billiard flows (and some geodesic flows), one should be able to replace the condition (P) by 
just assuming Lipschitzness  of the stable and unstable laminations -- this will be the subject of 
some future work.}  is one of the assumptions mentioned above:

\ms

\noindent
{\sc (P)}:  {\it There exist  constants $C > 0$ and $0 < \alpha \leq \beta$ such that for every $x\in \mt$
we have
$$\frac{1}{C} \, e^{\alpha_x \,t}\, \|u\| \leq \| d\phi_t(x)\cdot u\| \leq C\, e^{\beta_x\,t}\, \|u\|
\quad, \quad  u\in E^u(x) \:\:, t > 0 \;,$$
for some constants $\alpha_x, \beta_x > 0$ depending on $x$ but independent of $u$ and $t$ with
$\alpha \leq \alpha_x \leq \beta_x \leq \beta$ and $2\alpha_x - \beta_x \geq \alpha$ for all $x\in \mt$.}

\ms

Notice that when $N = 2$ this condition is always satisfied. For $N \geq 3$, some general conditions
on $K$ that imply (P) are given in \cite{kn:St5}.  According to general regularity results, (P) implies that 
$W^u_\ep(x)$ and $W^s_\ep(x)$ are Lipschitz in $x\in \mt$.  In fact, it follows from 
\cite{kn:Ha2} (see also \cite{kn:Ha1}) that assuming (P), the map 
$\mt \ni x \mapsto E^u(x)$ is $C^{1+\ep}$ with $\ep = 2 \inf_{x\in \mt}(\alpha_x/\beta_x) -1 > 0$, in the 
sense that this map has a linearization at any $x\in \mt$ that depends (uniformly H\"older) continuously on $x$.
The same applies to the map $\mt \ni x \mapsto E^s(x)$.

Next, we need some definitions from \cite{kn:St4}. Given $z \in \mt$, let $\exp^u_z : E^u(z) \longrightarrow W^u_{\ep_0}(z)$  and
$\exp^s_z : E^s(z) \longrightarrow W^s_{\ep_0}(z)$ be the corresponding
{\it exponential maps}. A  vector $b\in E^u(z)\setminus \{ 0\}$ will be called  {\it tangent to $\mt$} at
$z$ if there exist infinite sequences $\{ v^{(m)}\} \subset  E^u(z)$ and $\{ t_m\}\subset \R\setminus \{0\}$
such that $\exp^u_z(t_m\, v^{(m)}) \in \mt \cap W^u_{\ep}(z)$ for all $m$, $v^{(m)} \to b$ and 
$t_m \to 0$ as $m \to \infty$.  It is easy to see that a vector $b\in E^u(z)\setminus \{ 0\}$ is  tangent to $\mt$ at
$z$ if there exists a $C^1$ curve $z(t)$ ($0\leq t \leq a$) in $W^u_{\ep}(z)$ for some $a > 0$ 
with $z(0) = z$ and $\dot{z}(0) = b$, and $z(t) \in \mt$ for arbitrarily small $t >0$. In a similar way one defines tangent vectors to $\mt$ in $E^s(z)$.  

Denote by $d\alpha$ the standard symplectic form on $T^*(\R^N) = \R^N \times \R^N$. 
The following condition says that $d\alpha$ is in some sense  non-degenerate on the `tangent space' of $\mt$ near some its points:

\ms

\noindent
{\sc (ND)}:  {\it There exist $z_0\in \mt$, $\ep > 0$  and  $\mu_0 > 0$ such that
for any $\hz \in \mt \cap W^u_{\ep}(z_0)$ and any unit vector $b \in E^u(\hz)$ tangent to 
$\mt$ at $\hz$ there exist $\tz \in \mt \cap W^u_{\ep}(z_0)$ arbitrarily close to $\hz$ and  
a unit vector $a \in E^s(\tz)$ tangent to $\mt$ at $\tz$ with $|d\alpha(a,b) | \geq \mu_0$.}

\bs

\noindent
{\bf Remark 6.} Clearly the above is always true for $N = 2$. It was shown very recently in  \cite{kn:St5}  that 
for $N \geq 3$ this conditions is always satisfied for open billiard flows satisfying the pinching condition (P).

\bs

It follows from the hyperbolicity of $\mt$  that if  $\epsilon > 0$ is sufficiently small, there exists $\delta > 0$ such that if $x,y\in \mt$ and 
$d (x,y) < \delta$, then  $W^s_\ep(x)$ and $\phi_{[-\ep,\ep]}(W^u_\ep(y))$ intersect at exactly 
one point $[x,y] \in \mt$  (cf. \cite{kn:KH}). That is, there exists a unique $t\in [-\ep, \ep]$ such that $\phi_t([x,y]) \in W^u_{\ep}(y)$. 
Setting $\Delta(x,y) = t$, defines the so called {\it temporal distance function}. 
Given $E\subset \mt$, we will denote by $\Int_{\mt}(E)$ and $\partial_{\mt} E$  the {\it interior} 
and the {\it boundary} of the subset $E$ of $\mt$ in the topology of $\mt$, and by $\diam(E)$ the {\it diameter} of $E$.
Following \cite{kn:D}, a subset $R$ of $\mt$ will be called a {\it rectangle} if it has the form $R = [U,S] = \{ [x,y] : x\in U, y\in S\}\;,$ 
where $U$ and $S$ are subsets of $W^u_\ep(z) \cap \mt$ and $W^s_\ep(z) \cap \mt$, respectively, for some $z\in \mt$
that coincide with the closures of their interiors in $W^u_\ep(z) \cap \mt$ and $W^s_\ep(z) \cap \mt$.

Let $\rr = \{ R_i\}_{i=1}^k$ be a Markov family of rectangles  $R_i = [U_i  , S_i ]$ for $\mt$
(see e.g. \cite{kn:Bo1}, \cite{kn:D} or \cite{kn:St4} for the definition). Set  $R =  \cup_{i=1}^k R_i\; ,$
denote by  $\pp: R \longrightarrow R$ the corresponding {\it Poincar\'e map}, and
by $\tau$  the {\it first return time} associated with $\rr$. Then $\pp(x) = \phi_{\tau(x)}(x) \in R$ for any $x\in R$.
Notice that  $\tau$ is constant on each stable fiber of each $R_i$. 
We will assume that {\it size} $\chi = \max_i \diam(R_i)$ of the Markov
family $\rr = \{ R_i\}_{i=1}^k$ is sufficiently small so that each rectangle $R_i$ is between
two boundary components $\Gamma_{p_i}$ and $\Gamma_{q_i}$ of $K$, that is for any
$x\in R_i$, the first backward reflection point of the billiard trajectory $\gamma$ determined by $x$
belongs to $\Gamma_{p_i}$, while the first forward reflection point of $\gamma$ belongs to $\Gamma_{q_i}$.

Moreover, using the fact that the intersection of $\mt$ with each
cross-section to the flow $\phi_t$ is a Cantor set, we may assume that the Markov family
$\rr$ is chosen in such a way that

\ms

(i)  for any $i = 1, \ldots, k$ we have $\partial_\mt U_i = \e$.

\ms

Finally, partitioning each $R_i$ into finitely many smaller rectangles if necessary and removing
some `unnecessary' rectangles from the family formed in this way, we may assume that

\ms

(ii) for every $x\in R$ the billiard trajectory of $x$ from $x$ to $\pp(x)$ makes exactly one
reflection.

\ms

From now on we will assume that $\rr = \{ R_i\}_{i=1}^k$ is a fixed Markov family for  $\phi_t$
of size $\chi < \ep_0/2$ satisfying the above conditions (i) and (ii). Set  
$$U = \cup_{i=1}^k U_i\;.$$
The  map $\ts : U   \longrightarrow U$ is given by
$\ts  = \piU \circ \pp$, where $\piU : R \longrightarrow U$ is the {\it projection} along stable leaves.

Let $\aa = (\aa_{ij})_{i,j=1}^k$ be the matrix given by $\aa_{ij} = 1$ 
if $\pp(R_i) \cap R_j \neq  \e$ and $\aa_{ij} = 0$ otherwise. 
Consider the  symbol space 
$$\Sigma_\aa = \{  (i_j)_{j=-\infty}^\infty : 1\leq i_j \leq k, \aa_{i_j\; i_{j+1}} = 1
\:\: \mbox{ \rm for all } \: j\; \},$$
with the product topology and the {\it shift map} $\sigma : \Sigma_\aa \longrightarrow \Sigma_\aa$ 
given by $\sigma ( (i_j)) = ( (i'_j))$, where $i'_j = i_{j+1}$ for all $j$.
As in \cite{kn:Bo1} one defines a natural map  $\Psi : \sA \longrightarrow R\;.$
Namely, given any $(i_j)_{j=-\infty}^\infty \in \Sigma_\aa$ there
is exactly one point $x\in R_{i_0}$ such that $\pp^j(x) \in R_{i_j}$ for
all integers $j$. We then set $\Psi( (i_j)) = x$. One checks that
$\Psi \circ \sigma = \pp\circ \Psi$ on $R$.
It follows from the condition (i) above that {\bf the map $\Psi$ is a  bijection}.   

In a similar way one deals with the one-sided subshift 
$$\sAA = \{  (i_j)_{j=0}^\infty : 1\leq i_j \leq k, \aa_{i_j\; i_{j+1}} = 1
\:\: \mbox{ \rm for all } \: j \geq 0\; \},$$
where the {\it shift map} $\sigma : \sAA \longrightarrow \sAA$ is defined in the same way.
There exists a unique map $\psi : \sAA \longrightarrow U$ such that 
$\psi\circ \pi = \pi^{(U)}\circ \Psi$, where $\pi : \sA \longrightarrow \sAA$
is the {\it natural projection}.

Notice that the {\it roof function} $r : \sA \longrightarrow [0,\infty)$ defined by
$r (\xi) = \tau(\Psi(\xi))$ depends only on the
forward coordinates of $\xi\in \Sigma_\aa$. Indeed, if $\xi_+ = \eta_+$, where
$\xi_+ = (\xi_j)_{j=0}^\infty$, then for $x = \Psi(\xi)$ and $y = \Psi(\eta)$
we have $x, y\in R_i$ for $i = \xi_0= \eta_0$ and $\pp^j(x)$ and  $\pp^j(y)$ 
belong to the same $R_{i_j}$ for all $j\geq 0$. This implies that $x$ and $y$ belong to
the same local stable fibre in $R_i$ and by condition (ii), it follows that $\tau(x) = \tau(y)$.
Thus, $r (\xi) = r (\eta)$. So, we can define a {\it roof function} $r : \sAA \longrightarrow [0,\infty)$ 
such that $r\circ \pi =  \tau\circ \Psi$.

Let $B(\sAA)$ be the {\it space of  bounded functions} 
$g : \sAA \longrightarrow \C$ with its standard norm  
$\|g\|_0 = \sup_{\xi\in \sAA} |g(\xi)|$. Given a function $g \in B(\sAA)$, the  {\it Ruelle transfer operator } 
$\ll_g : B(\sAA) \longrightarrow B(\sAA)$ is defined by
$(\ll_gh)(\eta) = \sum_{\sigma(\eta) = \xi} e^{g(\eta)} h(\eta)\;.$
Denote by $\clip (U)$ the {\it space of Lipschitz functions} $h : U \longrightarrow \C$,
and for $h\in  \clip(U)$ let $\Lip(h)$ denote the {\it Lipschitz constant} of $h$.
For $t\in \R$, $|t| \geq 1$, define
$$\|h\|_{\lip,t} = \|h\|_0 + \frac{\Lip(h)}{|t|} \quad , \quad \|h\|_0 = \sup_{x\in U} |h(x)|\;.$$

Given a real-valued function $g$ on $ \sAA$ with $g\circ \psi^{-1} \in \clip(U)$,
there exists a unique number $s(g) \in \R$ such that $\Pr (-s(g)\, r + g) = 0$. Notice that if
$G : \mt \longrightarrow \C$ is a continuous function such that 
$(g\circ \psi^{-1}\circ \piU)(x) = \int_0^{\tau(x)} G(\phi_t(x))\, dt$ ($x \in R$), then
$s(g) = \Pr_{\phi_t}(G)$, the topological pressure of $G$ with respect to the flow
$\phi_t$ on $\mt$ (see e.g. Ch. 6 in \cite{kn:PP}).

The following is an immediate consequence of the main result in \cite{kn:St4}, taking 
into account the particular considerations for open billiard flows in \cite{kn:St5}.

\begin{thm}   Assume that the billiard flow $\phi_t$ over $\mt$ satisfies the conditions
{\rm (P)} and {\rm (ND)}. Let  $g: \sAA \longrightarrow \R$ be such that $g\circ \psi^{-1} \in \clip(U)$. Then 
there exist constants $a > 0$, $t_0 \geq 1,$ $\sigma(g) < s(g),\: C  > 0$ and 
$0 < \rho < 1$ so that for any $s = \tau + \ii t$ with $\tau \geq \sigma(g)$, $|\tau|\leq a$ and $|t| \geq t_0$,
any integer $n\geq 1$ and any function $v : \Sigma_{\mathcal A}^+ \longrightarrow \C$ with $v \circ \psi^{-1} \in \clip(U)$,
writing $n = p[\log |t|] + l,\: p \in \N,\: 0 \leq l \leq [\log |t|] - 1$,   we have
\begin{equation} \label{eq:11.1}
\|\left(\ll_{-sr + g}^n\ v\right)\circ \psi^{-1}\|_{\lip, t} \leq C \rho^{p[\log|t|]} e^{l {\rm Pr} (-\tau r + g)}\|v\circ \psi^{-1} \|_{\lip, t} \;.
\end{equation}
\end{thm}

\ms

\noindent
{\bf Remark 7.} Another way to state the above estimate is the following (\cite{kn:D}, \cite{kn:St4}):
{\it For every $g: \sAA \longrightarrow \R$ with $g\circ \psi^{-1} \in \clip(U)$ and every
$\epsilon > 0$ there exist constants $0 < \rho < 1$, $a_0 > 0$ and  $C > 0$ such that for any integer $m > 0$,
any $s = \tau + \i\, t\in \C$ with $|\tau| \leq a_0$, $|t| \geq 1/a_0$ and any function $v : \Sigma_{\mathcal A}^+ \longrightarrow \C$
with $v \circ \psi^{-1} \in \clip(U)$ we have}
$$\|\left(\ll_{-sr + g}^m\ v\right)\circ \psi^{-1}\|_{\lip, t} \leq C \;\rho^m \;|t|^{\ep}\; \|v\circ \psi^{-1} \|_{\lip, t}\; .$$

\bs

In the remaining part of this section, following \cite{kn:PS2},  we show  how to apply the Dolgopyat type estimates (11.1) to 
obtain the estimates of $\|L_s^n \gg_s \tv_s\|_{\Gamma,0}$ required in Sect. 5. The problem is
that the operator $L_s$ acts on $C(\saa)$, that is, it is related to the coding of billiard trajectories 
by means of the components of $K$, while the Dolgopyat type estimates apply to Ruelle
transfer operators $\mathcal L_{-sr + g}$ defined by means of Markov families and acting on functions $v$ such that
$v\circ \psi^{-1}$ is Lipschitz with respect to the standard metric in the phase space. Here we 
describe how the two types of Ruelle transfer operators relate, and show that the function
$(\gg_s\tv_s) \circ \psi^{-1}$ is Lipschitz. This  makes it possible to apply (11.1).

Apart from the coding described above, we can also use the coding of the flow over $\mt$ by using the boundary 
components of $K$ described in Sect. 3 above. We will use the notation from there, notably $f(\xi)$, $g(\xi)$, $\eta^{(k)}$
for any $k = 1, \ldots, \kappa_0$, $e(\xi)$, $\chi_f = \chi_1$, $\chi_g = \chi_2$, $\tf(\xi)$ and $\tg(\xi)$. Define the map
$\Phi : \sa \longrightarrow \mtb = \mt \cap S^*_\mt(\Omega)$
by $\Phi(\xi) = (P_0(\xi), (P_1(\xi) - P_0(\xi))/ \|P_1(\xi) - P_0(\xi)\|)$. Then $\Phi$ is a bijection
such that $\Phi\circ \sigma = B \circ \Phi$, where $B : \mtb \longrightarrow \mtb$ is the {\it billiard ball map}.
As before, given any function $G \in B(\saa)$, the  {\it Ruelle transfer operator } $L_G : B(\saa) \longrightarrow B(\saa)$ is defined by
$(L_G H)(\xi) = \sum_{\sigma(\eta) = \xi} e^{G(\eta)} H(\eta)\;.$

Let  $\omega : V_0 \longrightarrow S_{\dk}^*(\Omega)$ be the backward shift along the flow
defined in Sect. 3 on some neighborhood $V_0$ of $\mt$ in $S^*(\Omega)$.
Consider the bijection $\ss  = \Phi^{-1}\circ \omega \circ \Psi : \sA \longrightarrow \sa$.
Its restriction to $\sAA$ defines a bijection $\ss  : \sAA \longrightarrow \saa$.
Moreover $\ss \circ \sigma = \sigma \circ \ss$.
Define the function $g' : \sA \longrightarrow \R$ by $g'(\iu) = g(\ss(\iu))$. 

Next, for any $i = 1, \ldots, k$ choose $\hju^{(i)} = (\ldots, j_{-m}^{(i)}, \ldots,  j_{-1}^{(i)})$ such that
$(\hju^{(i)}, i) \in \sA^-$. It is convenient to make this choice in such a way that $\hju^{(i)}$ {\bf corresponds to the
local unstable manifold} $U_i \subset \mt \cap W^u_\ep(z_i)$, i.e.
the backward itinerary of every $z \in U_i$ coincides with $\hju^{(i)}$.
Now for any $\iu = (i_0, i_1, \ldots) \in \sAA$ (or $\iu \in \sA$) set
$\he(\iu) = (\hju^{(i_0)}; i_0,i_1, \ldots) \in \sA\;.$
According to the choice of $\hju^{(i_0)}$, we then have $\Psi(\he(\iu)) = \psi(\iu) \in U_{i_0}$.
(Notice that without the above special choice we would only have that 
$\Psi(\he(\iu))$ and  $\psi(\iu) \in U_{i_0}$ lie on the same stable leaf in $R_{i_0}$.)
Next, define
$\di \hchi_g (\iu) = \sum_{n=0}^\infty \left[ g'(\sigma^n(\iu)) - g'(\sigma^n \, \he(\iu)) \right]$ for $\iu \in \sA$.
As before, the function $\hg : \sA \longrightarrow \R$ given by
$\hg(\iu) = g'(\iu) - \hchi_g(\iu) + \hchi_g(\sigma \,\iu)$
depends on future coordinates only, so it can be regarded as a function on $\sAA$.

We will now describe a natural relationship between the operators 
$\ll_V : B(\sAA) \longrightarrow B(\sAA)$  and $L_v : B(\saa) \longrightarrow B(\saa)$ 
with $v$ appropriately defined by means of $V$.

First define  $\Gamma : B (\sa) \longrightarrow B(\sA)\;$ by
$\Gamma(v) = v\circ \Phi^{-1}\circ \omega \circ \Psi = v\circ \ss$. 
Since by property (ii) of the Markov family, $\omega : R \longrightarrow \mtb$ is a bijectiion, it
follows that $\Gamma$ is a bijection and $\Gamma^{-1}(V) = V\circ \Psi^{-1}\circ \omega^{-1}\circ \Phi$.
Moreover, $\Gamma$ induces a bijection $\Gamma : B (\saa) \longrightarrow B(\sAA)\;$.
Indeed, assume that $v \in B(\sa)$ depends on future coordinates only. Then $v\circ \Phi^{-1}$
is constant on local stable manifolds in $S^*_\mt(\Omega)$. Hence $v\circ \Phi^{-1} \circ \omega$
is constant on local stable manifolds on $R$, and therefore 
$\Gamma(v) = v\circ \Phi^{-1}\circ \omega \circ \Psi$ depends on future coordinates only.

Next, let $v, w\in B(\saa)$ and let $V = \Gamma(v)$, $W = \Gamma(w)$.
Given $\iu,\ju\in \sAA$ with $\sigma (\ju) = \iu$, setting $\xi = \ss(\iu)$ and $\eta = \ss (\ju)$, we
have $\sigma (\eta) = \xi$. Thus,
\begin{eqnarray*}
\ll_W V (\iu) =  \sum_{\sigma(\ju) = \iu} e^{W(\ju)}\, V(\ju) = \sum_{\sigma(\ju) = \iu} e^{w(\ss (\ju))}\, v(\ss (\ju))
= L_w v (\xi)
\end{eqnarray*}
for all $\iu \in \sAA$. This shows that $(L_w v) \circ \ss = \ll_{\Gamma(w)} \Gamma(v)$.

The equality 
\begin{equation} \label{eq:11.2}
\Pr (-\tau r + \hat{g}) = \Pr (-\tau \tilde{f} + \tilde{g})
\end{equation} 
and the following proposition are established in Section 3 in \cite{kn:PS2}.

\begin{prop} Assume that the map $\Lambda \ni x\mapsto W^u_{\ep}(x)$ is Lispchitz. Then
there exist  Lipschitz functions $\delta_1,\delta_2: U \longrightarrow \R$ such that setting
$\di \hde_s(\iu) = e^{s\, \delta_1(\psi(\iu))+ \delta_2(\psi(\iu))}$, we have
\be \label{eq:11.3}
\left(L_{-s\, \tf+\tg}^{n} u\right) (\ss(\iu))  = \frac{1}{\hde_s(\iu)}  \cdot  
\ll_{-s\, r+\hg}^{n} \left(\hde_s \cdot (u\circ \ss) \right)  (\iu) \quad, \quad \iu \in \sAA\:, \: s\in \C\;,
\ee
for any $u \in C(\saa)$ and any integer $n \geq 1$. 
\end{prop}

Combining (\ref{eq:11.1}), (\ref{eq:11.2}), (\ref{eq:11.3}), we deduce the following

\begin{thm}[\cite{kn:PS2}]   Assume that the billiard flow $\phi_t$ over $\mt$ satisfies the conditions
{\rm (P)} and {\rm (ND)}.  Then  there exist constants $a > 0$, $t_0 \geq 1,$ $\sigma_0 < s_0$, $C'  > 0$ and 
$0 < \rho < 1$ so that for any $s = \tau + \i\, t \in \C$ with $\tau \geq \sigma_0$, $|\tau| \leq a,\:|t| \geq t_0$, any
integer $n \geq 1$ and any function  $u: \saa \longrightarrow \R$ with $u \circ \ss\circ \psi^{-1} \in \clip(U)$, 
writing  $n = p[\log |t|] + l,\: p \in \N,\: 0 \leq l \leq [\log |t|] - 1$, we have
\begin{equation} \label{eq:11.4}
\left\|\left(L_{-s \tf + \tg }^n \, u\right)\circ \ss\circ \psi^{-1} \right\|_{\lip, t} 
\leq C' \rho^{p[\log|t|]} e^{l P (-\tau \tf + \tg)}\|u \circ \ss\circ \psi^{-1}\|_{\lip, t}.
\end{equation}
\end{thm}
The estimate (3.3) is a consequence of (\ref{eq:11.4}) and it could hold even if the assumption
(P) is not fulfilled (see Remark 6 above for (ND)).

\ms

Next, for the needs of Sect. 5 above, we have to estimate $\|L^n_{-s\tf + \tg}\gg_s \tv_s\|_{\Gamma,0}$,
where the operator $\gg_s$ is defined in Sect. 3. For any integer $n \geq 0$ we have
\begin{eqnarray*}
L^n_{-s\tf + \tg}\gg_s v (\xi) 
& = & \sum_{\sigma^n \eta = \xi} \sum_{\sigma \zeta = \eta} e^{-s\, \tf_n(\eta) + \tg_n(\eta)}\;
e^{-\phi^+(\zeta, s)  - s \tf(\zeta) + \tg(\zeta)} v(\zeta)\\
& = & \sum_{\sigma^{n+1} \zeta = \xi}  e^{-s\, \tf_{n+1}(\zeta) + \tg_{n+1}(\zeta)}\;
e^{-\phi^+(\zeta, s)} v(\zeta)
=  L^{n+1}_{-s\tf + \tg} (e^{-\phi^+(\cdot, s)} v)\, (\xi)\;.
\end{eqnarray*}
Thus, it is enough to estimate $\|L^{n+1}_{-s\tf + \tg} (e^{-\phi^+(\cdot, s)}\, \tv_s)\|_{\Gamma,0}$.
As in Sects. 3-5 above, we will consider these operators over $\Gamma_1$.

Given $s\in \C$, consider the functions $w_s : U_1 \longrightarrow \R$ and  $\hw_s : \sAA \longrightarrow \R$ defined by
$$w_s(x) = w_s(\psi(\iu)) = \hw_s(\iu) = e^{-\phi^+(\xi, s)} \tv_s(\xi) \;,$$
for $x= \psi(\iu) \in U_1$, $\iu\in \sAA$, $\xi = \ss(\iu)$. In order to use the Dolgopyat type
estimate (3.3), we have to show that $w_s$ is Lispchitz on $U_1$. We will deal in details with
$$w^{(1)}_s(x) =  e^{s\, \sum_{n = 0}^{\infty} [f(\sigma^n e(\xi)) - f_n^{+}(\xi)] - s\, \varphi(Q_0(\xi))} \; h(Q_0(\xi))\; ;$$
in a similar way one can deal with 
$w^{(2)}_s(x) =  e^{-\sum_{n = 0}^{\infty} [g(\sigma^n e(\xi)) - g_n^{+}(\xi)] }$.
It follows from the definitions of $\phi^+(\xi, s)$ and  $\tv_s$ in Sect. 3 that $w_s(x) = w^{(1)}_s(x) \, w^{(2)}_s(x)$.

Fix an arbitrary point $y_1 \in \mt$ such that
$\eta^{(1)} \in \sa^-$ corresponds to the local unstable manifold $\wuloc(y_1)$, 
i.e. the backward itinerary of every $z \in \wuloc(y_1)\cap V_0$ coincides with $\eta^{(1)}$.
It follows from the Lipschitzness of the stable and unstable laminations that the map
$\H_1: U_1 \longrightarrow \wuloc(y_1)$ defined by  $\H_1(x) = \phi_{\Delta(x,y_1)}([x, y_1])$ is Lipschitz.
Here $\Delta$ is the  temporal distance function defined in the beginning of this section. 

Next, consider the $N$-dimensional submanifold
$X = \{(q, q + t \nabla \varphi(q) : q \in \Gamma_1\;, \; 0 < t \}$
of  $S^*(\R^N)$ and the (stable) holonomy map $\H : \wuloc(y_1)\cap \mt \longrightarrow X$
defined by $\H(y) = \wsloc(y) \cap X$. Since $\varphi$ satisfies Ikawa's condition (${\mathcal P}$), it is easy to
see that $\wsloc(y)$ is transversal to $X$, so $\H(y) = \wsloc(y) \cap X$ is well-defined for 
$y \in \wuloc(y_1) \cap \mt$. Moreover, it follows from our assumtions that the stable (and unstable) holonomy maps
for the billiard flow $\phi_t$ are Lipschitz. In particular, $\H$ is Lipschitz.

We can now write down $w_s^{(1)}(x)$ using the maps $\H$ and $\H_1$ as follows.
Given $x \in U_1$, we have $x= \psi(\iu)$ for some $\iu\in \sAA$, with $i_0 = 1$. Setting
$\xi = \ss(\iu)$, we then have $\xi_0 = 1$. For any integer $m > 1$ consider
\begin{eqnarray*}
B_m
& = & \sum_{n=0}^{m-1}  [f(\sigma^n e(\xi)) - f_n^{+}(\xi)] -  \varphi(Q_0(\xi))\;.
\end{eqnarray*}
Setting $y = \H_1(x) \in \wuloc(y_1)$ and $z = \H(y)$, we have that $z\in \wsloc(y)$, and moreover
$\omega(z) = (Q_0(\xi), \nabla \varphi(Q_0(\xi)))$. Thus, 
$Q_0(\xi) = \pr_1(\omega(z)) = \pr_1(\omega(\H(\H_1(x))))$ is  Lipschitz in $x\in U_1$.
Next, set $\ep(u) = \| \pr_1(u) - \pr_1(\omega(u))\|$; then $u = \phi_{\ep(u)}(\omega(y))$ and 
$\ep(u)$ is a smooth function on an open subset of $S^*(\Omega)$  (where
$\omega$ is defined and takes values in $S^*_{\Gamma_1}(\Omega)$).
For $B_m$ we have
$$B_m = O(\theta^m) + \ep(y) - \ep(z) - \varphi(\omega(z))
=  O(\theta^m) + \ep(y) -  \varphi(z)\;,$$
and letting $m \to \infty$ we get
$$w^{(1)}_s(x) = e^{s[\ep(y) -  \varphi(z)]}\, h(\omega(z))
= e^{s[\ep(\H_1(x)) -  \varphi(\H(\H_1(x)))]}\, h(\omega(\H(\H_1(x))))\;,$$
so $w^{(1)}_s(x)$ is Lipschitz in $x\in U_1$. Moreover, for $x \in U_1$ and bounded $\Re (s)$ we obtain an uniform
bound for the Lipschitz norm of $w^{(1)}_s(x)$. The same argument works for $w^{(2)}_s(x)$.

\end{document}